\documentclass{siamonline181217} 
\usepackage{enumitem}
\usepackage{amsmath}
\usepackage{amssymb}
\usepackage{pifont}
\usepackage[usenames,dvipsnames]{pstricks}
\usepackage{pstricks-add}
\usepackage{pst-grad}
\usepackage{pst-plot} 
\usepackage{pst-eucl} 
\usepackage{theorem} 
\usepackage{euscript}
\usepackage{exscale,relsize}
\usepackage{epic,eepic}
\usepackage{multicol}
\usepackage{tikz}
\usepackage{pgfplots}
\usepackage{makeidx}
\usepackage{srcltx}         
\usepackage{url}
\usepackage{charter}
\usepackage{mathrsfs}       
\usepackage{graphicx}
\definecolor{Brown}{rgb}{0.45,0.0,0.05}
\definecolor{dgreen}{rgb}{0.00,0.49,0.00}
\definecolor{dblue}{rgb}{0,0.08,0.75}
\definecolor{lblue}{rgb}{0,0.7,0.75}
\renewcommand{\leq}{\ensuremath{\leqslant}}
\renewcommand{\geq}{\ensuremath{\geqslant}}

\newcommand{\minimize}[2]{\ensuremath{\underset{\substack{{#1}}}%
{\text{minimize}}\;\;#2 }}
\newcommand{\Frac}[2]{\displaystyle{\frac{#1}{#2}}} 
 
\newcommand{\scal}[2]{{\langle{{#1}\mid{#2}}\rangle}}

\newcommand{\abscal}[2]{\left|\left\langle{{#1}\mid{#2}}%
\right\rangle\right|} 
\newcommand{\menge}[2]{\big\{{#1}~\big |~{#2}\big\}}

\newcommand{\Argmin}{\ensuremath{{\text{\rm Argmin}}}}
\newcommand{\HH}{\ensuremath{{\mathcal H}}}
\newcommand{\GG}{\ensuremath{{\mathcal G}}}
\newcommand{\KK}{\ensuremath{{\mathcal K}}}
\newcommand{\Sum}{\ensuremath{\displaystyle\sum}}

\newcommand{\emp}{\ensuremath{{\varnothing}}}
\newcommand{\Id}{\ensuremath{\operatorname{Id}}}

\newcommand{\RR}{\ensuremath{\mathbb{R}}}
\newcommand{\RP}{\ensuremath{\left[0,+\infty\right[}}

\newcommand{\BL}{\ensuremath{\EuScript B}\,}
\newcommand{\RPP}{\ensuremath{\left]0,+\infty\right[}}

\newcommand{\soft}[1]{\ensuremath{{\:\text{\rm soft}}_{{#1}}\,}}

\newcommand{\RX}{\ensuremath{\left]-\infty,+\infty\right]}}

\newcommand{\NN}{\ensuremath{\mathbb N}}

\newcommand{\bK}{\ensuremath{\mathbb K}}

\newcommand{\weakly}{\ensuremath{\:\rightharpoonup\:}}

\newcommand{\moyo}[2]{\ensuremath{\sideset{^{#2}}{}%
{\operatorname{}}\!\!#1}}

\newcommand{\range}{\ensuremath{\text{\rm range}\,}}

\newcommand{\pinf}{\ensuremath{{+\infty}}}

\newcommand{\dom}{\ensuremath{\text{\rm dom}\,}}

\newcommand{\prox}{\ensuremath{\text{\rm prox}}}
\newcommand{\proj}{\ensuremath{\text{\rm proj}}}

\newcommand{\sign}{\ensuremath{\text{\rm sign}}}

\newcommand{\sri}{\ensuremath{\text{\rm sri}\,}}

\newcommand{\infconv}{\ensuremath{\mbox{\small$\,\square\,$}}}

\newcommand{\rzeroun}{\ensuremath{\left]0,1\right]}}   
   
\newtheorem{Algorithm}[theorem]{Algorithm}
\newtheorem{problem}[theorem]{Problem}

\newtheorem{remark}[theorem]{Remark}
\newtheorem{example}[theorem]{Example}

\renewcommand\theenumi{(\roman{enumi})}
\renewcommand\theenumii{\alph{enumii})}

\numberwithin{equation}{section}

\title{Proximal Activation of Smooth Functions in
Splitting Algorithms for Convex Image Recovery\thanks{The work 
of P. L. Combettes was supported by the National Science Foundation
under grant CCF-1715671.}}
\author{Patrick L. Combettes\thanks{North Carolina State
University, Department of Mathematics, 
Raleigh, NC 27695-8205, USA (plc@math.ncsu.edu)} 
\and 
Lilian E. Glaudin\thanks{Sorbonne Universit\'e,
Laboratoire Jacques-Louis Lions, F-75005 Paris, France
(glaudin@ljll.math.upmc.fr)}
}

\begin{document}
{\scriptsize\noindent{\sc SIAM J. Imaging Sci.} \hfill to appear}

\maketitle

\begin{abstract}
Structured convex optimization problems typically involve a mix of
smooth and nonsmooth functions. The common practice is to activate
the smooth functions via their gradient and the nonsmooth ones via
their proximity operator. We show that, although intuitively
natural, this approach is not necessarily the most efficient
numerically and that, in particular, activating all the functions
proximally may be advantageous. To make this viewpoint viable
computationally, we derive a number of new examples of proximity 
operators of smooth convex functions arising in applications.
A novel variational model to relax inconsistent convex feasibility 
problems is also investigated within the proposed framework. 
Several numerical applications to image recovery are presented to 
compare the behavior of fully proximal versus mixed
proximal/gradient implementations of several splitting algorithms.
\end{abstract}


\section{Introduction}
\label{sec:1}

Splitting in convex optimization methods for image recovery can
be traced back to the influential work of Youla 
\cite{Youl78,Youl82}. 
The convex feasibility framework he proposed consists in
formulating the image recovery problem as that of finding an
image in a Hilbert space $\HH$ satisfying $m$ constraints derived
from \emph{a priori} knowledge and the observed data. The
constraints are represented by closed convex sets 
$(C_i)_{1\leq i\leq m}$ and the problem is therefore to 
\begin{equation}
\label{e:feas}
\text{find}\;x\in\bigcap_{i=1}^mC_i.
\end{equation}
Now, for every $i\in\{1,\ldots,m\}$, let $\proj_{C_i}$ be the
projection operator onto $C_i$, which maps each $x\in\HH$
to its unique closest point in $C_i$, that is,
\begin{multline}
\label{e:proj}
\proj_{C_i}\colon\HH\to\HH\colon 
x\mapsto\underset{y\in\HH}{\text{argmin}}\;
\bigg(\iota_{C_i}(y)+\frac{1}{2}\|x-y\|^2\bigg),\\
\;\text{where}\quad\iota_{C_i}\colon y\mapsto
\begin{cases}
0,&\text{if}\;\;y\in C_i;\\
\pinf,&\text{if}\;\;y\notin C_i.
\end{cases}
\end{multline}
The methodology of projection methods is to split the problem of
finding a point in $\bigcap_{i=1}^mC_i$ into a sequence of simpler
problems involving the sets $(C_i)_{1\leq i\leq m}$ individually
\cite{Baus96,Aiep96}. For instance, the POCS (Projection Onto
Convex Sets) algorithm advocated in \cite{Youl82} is governed by
the updating rule
\begin{equation}
(\forall n\in\NN)\quad
x_{n+1}=(\proj_{C_1}\circ\cdots\circ\proj_{C_m})x_n.
\end{equation}
Convex variational formulations arising in modern image recovery 
have complex structures that require sophisticated analysis tools 
and solution methods. Since projection operators are of limited use
beyond feasibility and best approximation problems, 
to solve such formulations, one strategy 
is to use an extended notion of a projection operator.
In \cite{Smms05} it was suggested to use Moreau's 
proximity operator \cite{Mor62b} for this purpose. Recall that the
proximity operator of a proper lower semicontinuous convex function
$\varphi\colon\HH\to\RX$ is
\begin{equation}
\label{e:prox}
\prox_{\varphi}\colon\HH\to\HH\colon 
x\mapsto\underset{y\in\HH}{\text{argmin}}\;
\bigg(\varphi(y)+\frac{1}{2}\|x-y\|^2\bigg),
\end{equation}
and that it reduces to \eqref{e:proj} when $\varphi=\iota_{C_i}$. We 
refer the reader to \cite[Chapter~24]{Livre1} for a detailed
account of the properties of proximity operators with various
examples, to \cite{Banf11} for a tutorial on proximal methods 
in signal processing, and to 
\cite{Aspe16,Beck09,Cham11,Chou15,Jsts07,Ocon14,Papa14,Ragu15} for 
specific applications to image recovery.
Current proximal splitting methods can handle highly structured 
convex minimization problems such as the following, which will be
the focus of our discussion (see below for notation).

\begin{problem}
\label{prob:1}
\rm
Let $\HH$ be a real Hilbert space, let $I$ and $J$ be
disjoint finite subsets of $\NN$ such that $K=I\cup J\neq\emp$, 
let $f\in\Gamma_0(\HH)$, and let $(\GG_k)_{k\in K}$ be a 
family of real Hilbert spaces.
For every $k\in K$, suppose that $L_k\colon\HH\to\GG_k$ is a
nonzero bounded linear operator.
For every $i\in I$, let $g_i\in\Gamma_0(\GG_i)$ and, 
for every $j\in J$, let $\mu_j\in\RPP$ and
let $h_j\colon\GG_j\to\RR$ be convex and differentiable with a 
$\mu_j$-Lipschitzian gradient. Assume that 
(see \cite[Proposition~4.3]{Svva12} for sufficient conditions)
\begin{equation}
\label{e:cq1}
0\in\range\bigg(\partial f+\sum_{i\in I}L_i^*\circ\partial g_i 
\circ L_i+\sum_{j\in J}L_j^*\circ(\nabla h_j)
\circ L_j\bigg).
\end{equation}
The goal is to 
\begin{equation}
\label{e:primal}
\minimize{x\in\HH}{f(x)+\sum_{i\in I}\,g_i(L_ix)
+\sum_{j\in J}\,h_j(L_jx)}. 
\end{equation}
\end{problem}

The principle of a splitting method for solving \eqref{e:primal} is 
to use separately each of the functions $f$, $(g_i)_{i\in I}$, and 
$(h_j)_{j\in J}$, and each of the operators $(L_k)_{k\in K}$, so as
to reduce the execution of the algorithm to a sequence of simple
steps. A prevalent viewpoint in first order convex splitting
algorithms is that to activate each function $\varphi$ appearing in
the model there are two options:
\begin{itemize}
\setlength{\itemsep}{-1pt}
\item
if $\varphi$ is smooth, i.e., real-valued and differentiable
everywhere with a Lipschitzian gradient, then use $\nabla\varphi$;
\item
otherwise, use $\varphi$ proximally, i.e., via its proximity 
operator \eqref{e:prox}.
\end{itemize}

\begin{figure}[ht!]
\scalebox{0.82} 
{
\begin{pspicture}(-9.7,-5.0)(9.2,4.3)
\psline[linewidth=0.04cm,arrowsize=0.07cm 4.0,%
arrowlength=1.4,arrowinset=0.4]{->}(-7,0)(7.5,0)
\psline[linewidth=0.04cm,arrowsize=0.07cm 4.0,%
arrowlength=1.4,arrowinset=0.4]{->}(0,-5)(0,5.0)
\psellipse[linewidth=0.04cm,linecolor=dgreen,rot=45]%
(0,0)(6.402,2.134)
\psellipse[linewidth=0.04cm,linecolor=red,rot=45](0,0)(1.17,0.39)
\psellipse[linewidth=0.04cm,linecolor=blue,rot=45](0,0)(4.95,1.65)
\rput(0.0,-4.0){$-$}
\rput(-0.4,-4.0){$-4$}
\rput(0.0,4.0){$-$}
\rput(-0.3,4.0){$4$}
\rput(0.0,2.0){$-$}
\rput(-0.2,1.8){$2$}
\rput(0.0,-2.0){$-$}
\rput(-0.4,-2.0){$-2$}
\rput(2.0,0.0){$|$}
\rput(2.0,-0.5){$2$}
\rput(4.0,0.0){$|$}
\rput(4.0,-0.4){$4$}
\rput(-6.0,0.0){$|$}
\rput(-6.0,-0.4){$-6$}
\rput(6.0,0.0){$|$}
\rput(6.0,-0.4){$6$}
\rput(-2.0,0.0){$|$}
\rput(-2.0,-0.4){$-2$}
\rput(-4.0,0.0){$|$}
\rput(-4.0,-0.4){$-4$}
\rput(7.9,0){$\xi_1$}
\rput(0,5.4){$\xi_2$}
\psline[linewidth=0.04cm,arrowsize=0.07cm 4.0,linecolor=dgreen,%
arrowlength=1.4,linestyle=dashed,arrowinset=0.4]{->}%
(3.682500,1.122500)(0.88535,2.85655)
\psline[linewidth=0.04cm,arrowsize=0.07cm 4.0,linecolor=dgreen,%
arrowlength=1.4,arrowinset=0.4]{->}%
(3.682500,1.122500)(4.53243,0.59560)
\rput(4.800000,0.520000){\color{dgreen} $d_2$}
\rput(1.500000,4.700000){\color{dgreen} $\bullet$}
\rput(1.200000,4.700000){\color{dgreen} $x_1$}
\rput(3.682500,1.122500){\color{dgreen} $\bullet$}
\rput(4.072500,1.182500){\color{dgreen} $x_2$}
\rput(0.837938,2.885937){\color{dgreen} $\bullet$}
\rput(0.537938,2.985937){\color{dgreen} $x_3$}
\rput(2.262202,0.623802){\color{dgreen} $\bullet$}
\rput(2.462202,0.823802){\color{dgreen} $x_4$}
\rput(0.462966,1.773686){\color{dgreen} $\bullet$}
\rput(0.522966,2.073686){\color{dgreen} $x_5$}
\rput(1.390991,0.342415){\color{dgreen} $\bullet$}
\rput(1.690991,0.272415){\color{dgreen} $x_6$}
\psline[linewidth=0.04cm,arrowsize=0.07cm 4.0,linecolor=blue,%
arrowlength=1.4,arrowinset=0.4]{->}%
(3.5,3.5)(4.2071,4.2071)
\rput(3.500000,3.500000){\blue $\bullet$}
\rput(3.900000,3.400000){\color{blue} $x_1$}
\rput(3.800000,4.350000){\color{blue} $d_1$}
\rput(3.062500,3.062500){\blue $\bullet$}
\rput(3.360000,3.062500){\color{blue} $x_2$}
\rput(2.584223,2.584223){\blue $\bullet$}
\rput(2.884223,2.584223){\color{blue} $x_3$}
\rput(2.094132,2.094132){\blue $\bullet$}
\rput(2.094132,2.394132){\color{blue} $x_4$}
\rput(1.618308,1.618308){\blue $\bullet$}
\rput(1.348308,1.818308){\color{blue} $x_5$}
\rput(1.178447,1.178447){\blue $\bullet$}
\rput(0.878447,1.378447){\color{blue} $x_6$}
\psline[linewidth=0.04cm,arrowsize=0.07cm 4.0,linecolor=red,%
arrowlength=1.4,arrowinset=0.4]{->}%
(0.860606,0.739394)(1.83182,0.97762)
\psline[linewidth=0.04cm,arrowsize=0.07cm 4.0,linecolor=red,%
arrowlength=1.4,linestyle=dashed,arrowinset=0.4]{->}%
(0.860606,0.739394)(5.9229,1.9811)
\rput(6.000000,2.000000){\black $\bullet$}
\rput(6.000000,2.300000){\black $x_0$}
\rput(0.860606,0.739394){\red $\bullet$}
\rput(1.160606,0.599394){\red $x_1$}
\rput(1.860606,1.289394){\red $d_1$}
\rput(0.161837,0.158163){\red $\bullet$}
\rput(0.361837,0.358163){\red $x_2$}
\end{pspicture}
}
\caption{Comparison of the steepest descent method \eqref{e:j2}
(in green) and of the proximal point algorithm \eqref{e:martinet} (in
red) in $\HH=\RR^2$ for 
$\varphi\colon(\xi_1,\xi_2)\mapsto 9\xi_1^2-14\xi_1\xi_2+9\xi_2^2$.
The ellipsoids represent the level lines of $\varphi$.
The steepest descent method is implemented with $\gamma=1.8/\beta$
as this choice gave rise to the fastest convergence. On the other
hand, the proximal point algorithm is implemented with the default
choice $\gamma=1$ (larger values gave even faster convergence). 
The two algorithms behave quite differently, both in terms of
directions of movement and of trajectories. At iteration $n$, call
$d_n=\nabla\varphi(x_n)/\|\nabla\varphi(x_n)\|$ the normalized
gradient at $x_n$. Consider the action of the 
steepest descent, say at iteration $n=2$. The next iterate
$x_3$ is obtained by moving from $x_2$ in the direction opposite 
to the gradient at $x_2$. By contrast, consider the action of the
proximal point algorithm, say at iteration $n=0$. The next iterate
$x_1$ satisfies the implicit equation
$x_0-x_1=\gamma\nabla\varphi(x_1)$, which means that $x_1$ is
obtained by moving from $x_0$ in the direction opposite to the
gradient at $x_1$. Finally, we include the orbit (in blue) of the
inertial version of the steepest descent method obtained by setting
$f=0$ and $h=\varphi$ in Algorithm~\ref{al:cd}, and choosing the
parameters $\alpha=2.01$ and $\gamma=1/\beta$, which gave the
fastest convergence.
}
\label{fig:Jim}
\end{figure}
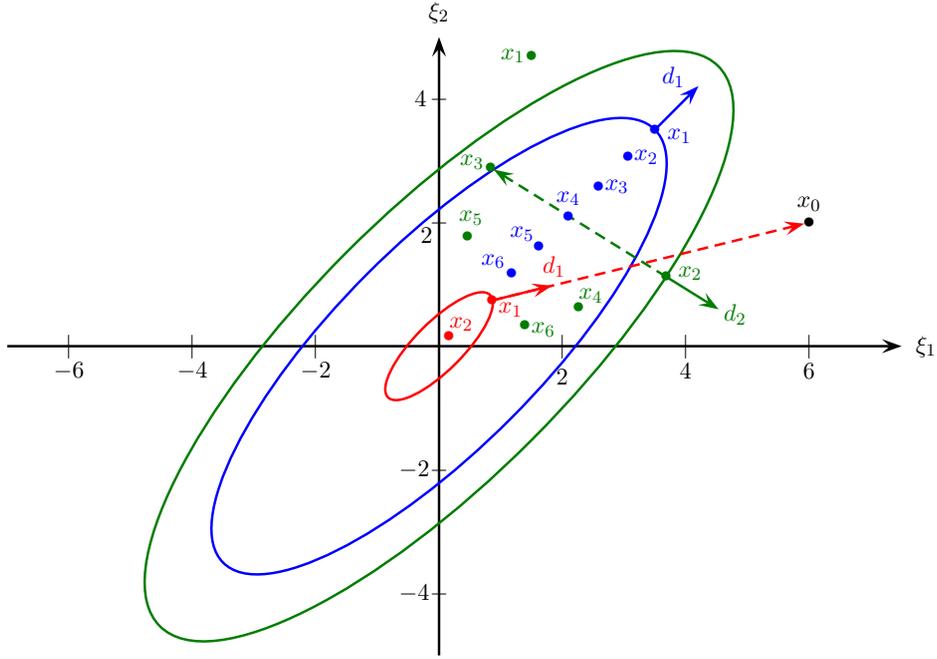 
\noindent
In the present paper we propose a more nuanced viewpoint and submit
that, when $\varphi$ is smooth, it may be computationally 
advantageous to activate it proximally when its proximity
operator can be implemented. To motivate this viewpoint, let us
first observe that a tight Lipschitz constant for the gradient of
$\varphi$ may not be easy to estimate (see, e.g.,
\cite{Abbo17,Nmtm09,Chaa11}), which limits the range of the
proximal parameters and may have a detrimental incidence on the
speed of convergence. Our second observation is that proximal steps
behave numerically quite differently from gradient steps, which may
have a positive impact on the asymptotic performance of algorithms. 
To illustrate this fact, consider the problem of minimizing a
differentiable convex function $\varphi\colon\HH\to\RR$ with a
$\beta$-Lipschitzian gradient (see 
Fig.~\ref{fig:Jim} for a concrete example). The associated
continuous-time gradient dynamics is \cite[Section~3.4]{Aub84a}
\begin{equation}
\label{e:j1}
x(0)=x_0\quad\text{and}\quad -\frac{dx(t)}{dt}=\nabla\varphi(x(t)).
\end{equation}
The forward Euler (explicit) discretization of this equation
with time step $\gamma\in\RPP$ assumes the form
$-(x_{n+1}-x_n)/\gamma=\nabla\varphi(x_n)$, which leads to the
steepest descent algorithm
\begin{equation}
\label{e:j2}
(\forall n\in\NN)\quad
x_{n+1}=x_n-\gamma\nabla\varphi(x_n).
\end{equation}
On the other hand, the backward Euler (implicit) discretization 
of \eqref{e:j1} is
$-(x_{n+1}-x_n)/\gamma=\nabla\varphi(x_{n+1})$, which leads to 
Martinet's proximal point algorithm \cite{Mart72} 
\begin{equation}
\label{e:martinet}
(\forall n\in\NN)\quad
x_{n+1}=\prox_{\gamma\varphi}x_n.
\end{equation}
Alternatively, it follows from 
\cite[Proposition~12.30]{Livre1} that the proximal point algorithm
coincides with the steepest descent method applied to the Moreau
envelope of $\varphi$, namely,
\begin{equation}
\label{e:martinet1}
(\forall n\in\NN)\quad
x_{n+1}=\prox_{\gamma\varphi}x_n=
x_n-\gamma\nabla\big(\moyo{\varphi}{\gamma}\big)(x_n),
\quad\text{where}\quad
\moyo{\varphi}{\gamma}=\varphi\infconv(q/\gamma).
\end{equation}
While the convergence of \eqref{e:martinet} is guaranteed for any
$\gamma\in\RPP$, that of \eqref{e:j2} requires 
$\gamma<2/\beta$ \cite[Chapter~28]{Livre1}, which results in
potentially slow convergence. Historically, the idea of using
proximal steps in smooth minimization problems can be found
in \cite[Section~5.8]{Bell66}. There, the problem under
consideration is the standard least-squares problem of minimizing
the smooth function 
$\varphi\colon\RR^N\to\RR\colon x\mapsto\|Ax-b\|^2/2$ in
connection with the numerical inversion of the Laplace transform.
Given $\gamma\in\RPP$ and $x_0\in\RR^N$, the algorithm proposed in 
\cite[Eq.~(5.8.3)]{Bell66} is
\begin{equation}
\label{e:bellman}
(\forall n\in\NN)\quad
x_{n+1}=\big(\Id+\gamma A^\top A\big)^{-1}\big(x_n+\gamma A^\top
b\big),
\end{equation}
and it is reported to be better than the standard steepest descent
approach. Remarkably, \eqref{e:bellman} is nothing but an early 
instance of the proximal point algorithm \eqref{e:martinet}.

The paper is organized as follows. In Section~\ref{sec:2}, we 
enrich the list of known proximity operators by providing new
closed form expressions for those of various smooth convex functions
commonly encountered in applications. This investigation is of
interest in its own right since some splitting algorithms operate
exclusively with proximal steps; see, e.g., 
\cite{Cham11,Icip14,Svva10,Jsts07}. 
In connection with the numerical
solution of Problem~\ref{prob:1}, we review in Section~\ref{sec:3}
some pertinent proximal splitting methods. Image recovery 
applications are presented in Section~\ref{sec:4}. Numerical
comparisons between splitting algorithms in which smooth functions
are activated via gradient steps and those in which all functions
are activated via their proximity operators are conducted. In
particular, in Section~\ref{sec:44}, we propose a new variational
model, based on Problem~\ref{prob:1} and the results of
Section~\ref{sec:2}, to relax inconsistent convex feasibility
problems.
While no universal conclusion may be drawn from these experiments,
they suggest that fully proximal splitting algorithms deserve to be
given serious consideration in applications.

{\bfseries Notation.} 
The notation follows that of \cite{Livre1}. 
Throughout, $\HH$ is a real Hilbert space with scalar product 
$\scal{\cdot}{\cdot}$, associated norm $\|\cdot\|$, and identity 
operator $\Id$. Weak convergence is denoted by $\weakly$.
Given a real Hilbert space $\GG$, we denote by
$\BL(\HH,\GG)$ the space of continuous linear operators from
$\HH$ to $\GG$. We set $q=\|\cdot\|^2/2$ and denote by 
$\Gamma_0(\HH)$ the class of lower semicontinuous 
convex functions $f\colon\HH\to\RX$ such that 
$\dom f=\menge{x\in\HH}{f(x)<\pinf}\neq\emp$.
Let $f\in\Gamma_0(\HH)$. Then $f^*$ denotes the
conjugate of $f$, $\partial f$ the 
subdifferential of $f$, and $f\infconv g$ the inf-convolution of
$f$ and $g\in\Gamma_0(\HH)$.
Let $C$ be a convex subset of $\HH$. The strong relative interior
of $C$ is denoted by $\sri C$, the indicator function of
$C$ by $\iota_C$, the distance function to $C$ by
$d_C$, the support function of $C$ by $\sigma_C$ and, 
if $C$ is nonempty and closed,
the projection operator onto $C$ by $\proj_C$.
The Hilbert direct sum of family of real
Hilbert spaces $(\HH_i)_{i\in I}$ is denoted by
$\bigoplus_{i\in I}\HH_i$; in addition if, for every $i\in I$,
$f_i\colon\HH_i\to[0,\pinf]$, then 
\begin{equation}
\bigoplus_{i\in I} f_i\colon \bigoplus_{i\in I}\HH_i
\to[0,\pinf]\colon (x_i)_{i\in I}\mapsto \sum_{i\in I}f_i(x_i).
\end{equation}
The standard Euclidean norm on $\RR^N$ is denoted by $\|\cdot\|_2$.

\section{Proximity operators of smooth convex functions}
\label{sec:2}

Let $\beta\in\RPP$, let $\gamma\in\RPP$, and let $h\colon\HH\to\RR$ 
be a convex function with a $\beta$-Lipschitzian gradient. Then
there exists a function $g\in\Gamma_0(\HH)$ such that 
$h=g^*\infconv(\beta q)$ \cite[Corollary~18.19]{Livre1}. In 
this case, we derive from
\cite[Propositions~12.30 and~24.8(vii)]{Livre1} that
\begin{equation}
\label{e:jjm69}
\nabla h=\beta\big(\Id-\prox_{g^*/\beta}\big)
\quad\text{and}\quad
\prox_{\gamma h}x=\Id+\dfrac{\gamma\beta}{\gamma\beta+1}
\big(\prox_{(\gamma\beta+1)g^*/\beta}-\Id\big).
\end{equation}
The closed form expression for $\prox_{\gamma h}x$ above
is however of limited 
use since $g^*$ and its proximity operator are usually not
available explicitly. Even when $\HH=\RR$, computing the
proximity operator of a smooth convex function may be involved:
for instance the derivative of $h\colon
x\mapsto\sqrt[3]{1+x^6}/3$ is $\sqrt[3]{2}$-Lipschitzian but
evaluating $\prox_{\gamma h}$ requires solving a high degree
polynomial equation. Nonetheless, as we now show, a variety of smooth
convex functions encountered in applications have readily
computable proximity operators.

\subsection{Functions involving distances}

The following fact will be needed.
\begin{lemma}{\rm\cite[Proposition~2.1]{Nmtm09}}
\label{l:3}
Let $C$ be a nonempty closed convex subset of $\HH$, let
$\phi\in\Gamma_0(\RR)$ be even, and set $\varphi=\phi\circ d_C$.
Then $\varphi\in\Gamma_0(\HH)$. Moreover, $\prox_{\varphi}=\proj_C$
if $\dom\phi={\{0\}}$ and, otherwise, for every $x\in\HH$,
\begin{equation}
\label{e:26juin2009}
\prox_{\varphi}x=
\begin{cases}
x+\displaystyle{\frac{\prox_{\phi^*}d_C(x)}{d_C(x)}}
\big(\proj_Cx-x\big),&\text{if}\;\;d_C(x)>\max\partial\phi(0);\\
\proj_Cx,&\text{if}\;\;x\notin C\;\text{and}\;\; d_C(x)
\leq\max\partial\phi(0);\\
x,&\text{if}\;\;x\in C.
\end{cases}
\end{equation}
\end{lemma}

We start with an example which leads to affine gradient and
proximal operators.

\begin{example}
\label{ex:-1}
\rm
Let $I$ be a nonempty finite set. For every $i\in I$, 
let $\GG_i$ be a real Hilbert space, let $V_i$ be a closed vector
subspace of $\GG_i$, let $r_i\in\GG_i$, let 
$L_i\in\BL(\HH,\GG_i)$, and let $\alpha_i\in\RPP$. Set 
$h\colon\HH\to\RR\colon x\mapsto(1/2)
\sum_{i\in I}\alpha_id^2_{V_i}(L_ix-r_i)$ and 
$Q=(\Id+\gamma\sum_{i\in I}\alpha_iL_i^*\proj_{V_i^\bot}L_i)^{-1}$.
Let $\gamma\in\RPP$, set 
$\beta=\sum_{i\in I}\alpha_i\|L_i\|^2$, and let $x\in\HH$.
Then $h\colon\HH\to\RR$ is convex and (Fr\'echet) differentiable
with a $\beta$-Lipschitzian gradient, 
\begin{equation}
\label{e:kj50}
\nabla h(x)=\sum_{i\in I}\alpha_iL_i^*\Big(\proj_{V_i^\bot}
\big(L_ix-r_i\big)\Big),
\quad\text{and}\quad
\prox_{\gamma h}x=Q
\bigg(x+\gamma\sum_{i\in I}\alpha_iL_i^*
\Big(\proj_{V_i^\bot}r_i\Big)\bigg).
\end{equation}
\end{example}
\begin{proof}
The convexity of $h$ is clear. We have
$h(x)=(1/2)\sum_{i\in I}\alpha_i\|\proj_{V_i^\bot}(L_ix-r_i)\|^2$ 
and $\nabla h(x)$ is given by \eqref{e:kj50} since 
$(\forall i\in I)$ 
$\nabla d_{V_i}^2/2=\Id-\proj_{V_i}=\proj_{V_i^\bot}$.
Moreover, 
\begin{equation}
(\forall i\in I)\quad\|L_i^*\proj_{V_i^\bot}L_i\|
\leq\|L_i^*\|\,\|\proj_{V_i^\bot}\|\,\|L_i\|\leq\|L_i\|^2.
\end{equation}
Hence, for every $y\in\HH$,
\begin{equation}
\|\nabla h(x)-\nabla h(y)\|=
\bigg\|\sum_{i\in I}\alpha_iL_i^*\Big(\proj_{V_i^\bot}L_i(x-y)\Big)
\bigg\|\leq\sum_{i\in I}\alpha_i\|L_i\|^2\|x-y\|=\beta\|x-y\|.
\end{equation}
Now set $p=\prox_{\gamma h}x$. Then we derive from \eqref{e:kj50}
that
\begin{equation}
\label{e:kj51}
x-p=\gamma\nabla h(p)=
\gamma\Bigg(\sum_{i\in I}\alpha_iL_i^*\proj_{V_i^\bot}L_i\Bigg)p
-\gamma\sum_{i\in I}\alpha_iL_i^*\Big(\proj_{V_i^\bot}r_i\Big),
\end{equation}
which yields the expression for $\prox_{\gamma h}x$.
\end{proof}

The next construction, which involves the distance function $d_C$
to a convex set $C$, will be seen to capture a broad range of
functions of interest.

\begin{example}
\label{ex:2}
\rm
Let $C$ be a nonempty closed convex subset of $\HH$, let
$\beta\in\RPP$, let $\phi\colon\RR\to\RR$ be even, convex, and
differentiable with a $\beta$-Lipschitzian derivative, and
set $h=\phi\circ d_C$. Let $\gamma\in\RPP$ and $x\in\HH$.
Then $h\colon\HH\to\RR$ is convex and Fr\'echet differentiable
with a $\beta$-Lipschitzian gradient, 
\begin{equation}
\label{e:kj80}
\nabla h(x)=
\begin{cases}
\dfrac{\phi'\big(d_C(x)\big)}{d_C(x)}\big(x-\proj_Cx\big),
&\text{if}\;\;x\notin C;\\
0,&\text{if}\;\;x\in C,
\end{cases}
\end{equation}
and 
\begin{equation}
\label{e:kj52}
\prox_{\gamma h}x=
\begin{cases}
\proj_Cx+\Frac{\prox_{\gamma\phi}d_C(x)}{d_C(x)}
(x-\proj_Cx),&\text{if}\;\;x\notin C;\\
x,&\text{if}\;\;x\in C.
\end{cases}
\end{equation}
\end{example}
\begin{proof}
Since $\phi$ and $d_C$ are convex and $\phi$ is increasing on
$\RP$, $h$ is convex \cite[Proposition~11.7(ii)]{Livre1}. In
addition, since \cite[Proposition~11.7(i)]{Livre1} asserts that $0$ 
is a minimizer of $\phi$, 
\begin{equation}
\label{e:sh12}
\partial\phi(0)=\{\phi'(0)\}=\{0\}.
\end{equation}
First, we derive \eqref{e:kj80} from \cite[Proposition
17.33(ii)]{Livre1} and \eqref{e:sh12}.
Next, we infer from \cite[Proposition 13.26]{Livre1} 
that $h^*=\sigma_C+\phi^*\circ\|\cdot\|$. We invoke
\cite[Theorem~18.15]{Livre1} to deduce that
$\phi^*$ is $(1/\beta)$-strongly convex and that $h^*$
is therefore likewise, and then to conclude that $\nabla h$ is
$\beta$-Lipschitzian. On the other hand, we derive from 
Lemma~\ref{l:3}, \eqref{e:sh12}, and 
Moreau's decomposition \cite[Remark~14.4]{Livre1} that
\begin{align}
\prox_{\gamma h}x
&=
\begin{cases}
x+\Frac{\prox_{\phi^*}d_C(x)}{d_C(x)}
\big(\proj_Cx-x\big),&\text{if}\;\;d_C(x)>\max\partial\phi(0);\\
\proj_Cx,&\text{if}\;\;d_C(x)\leq\max\partial\phi(0)
\end{cases}
\label{e:7}\\
&=
\begin{cases}
x+\Frac{d_C(x)-\prox_{\gamma\phi}d_C(x)}{d_C(x)}
(\proj_Cx-x),&\text{if}\;\;d_C(x)>0;\\
\proj_Cx,&\text{if}\;\;d_C(x)\leq 0
\end{cases}
\nonumber\\
&=
\begin{cases}
\proj_Cx+\Frac{\prox_{\gamma\phi}d_C(x)}{d_C(x)}
(x-\proj_Cx),&\text{if}\;\;x\notin C;\\
x,&\text{if}\;\;x\in C,
\end{cases}
\label{e:kj53}
\end{align}
as announced.
\end{proof}

We now investigate a generalization of the 
Vapnik $\varepsilon$-insensitive loss function \cite{Vapn00}.

\begin{example}[abstract smooth Vapnik loss function]
\label{ex:6bis}
\rm
Let $C$ be a nonempty closed convex subset of $\HH$,
let $\varepsilon\in\RPP$, let
$\beta\in\RPP$, let $\psi\colon\RR\to\RR$ be even, convex, and
differentiable with a $\beta$-Lipschitzian derivative, and
set $h=\psi\circ\max(d_C-\varepsilon,0)$. 
Let $\gamma\in\RPP$ and $x\in\HH$.
Then $h\colon\HH\to\RR$ is convex and Fr\'echet differentiable
with a $\beta$-Lipschitzian gradient, 
\begin{equation}
\label{e:kj81}
\nabla h(x)=
\begin{cases}
\dfrac{\psi'\big(d_C(x)-\varepsilon\big)}{d_C(x)}
\big(x-\proj_Cx\big),
&\text{if}\;\;d_C(x)>\varepsilon;\\
0,&\text{if}\;\;d_C(x)\leq\varepsilon,
\end{cases}
\end{equation}
and 
\begin{equation}
\label{e:kj5-5}
\prox_{\gamma h}x=
\begin{cases}
\proj_Cx+\dfrac{\varepsilon+\prox_{\gamma\psi}
\big(d_C(x)-\varepsilon\big)}{d_C(x)}
(x-\proj_C x),&\text{if}\;\;d_C(x)>\varepsilon;\\
x,&\text{if}\;\;d_C(x)\leq\varepsilon.
\end{cases}
\end{equation}
\end{example}
\begin{proof}
Let $\vartheta=\max(|\cdot|-\varepsilon,0)$ be the standard Vapnik 
loss function, set $\phi=\psi\circ\vartheta$, and let $\xi\in\RR$. 
Upon applying Example~\ref{ex:2} in $\RR$ with
$C=[-\varepsilon,\varepsilon]$, we obtain that 
$\phi\colon\RR\to\RR$ is convex and Fr\'echet differentiable
with a $\beta$-Lipschitzian derivative, that
\begin{equation}
\phi'(\xi)=
\begin{cases}
\psi'(|\xi|-\varepsilon)\sign(\xi),
&\text{if}\;\;|\xi|>\varepsilon;\\
0,&\text{if}\;\;|\xi|\leq\varepsilon,
\end{cases}
\end{equation}
and that
\begin{equation}
\label{e:kj55}
\prox_{\gamma\phi}\xi=
\begin{cases}
\big(\varepsilon+\prox_{\gamma\psi}(|\xi|-\varepsilon)\big)
\sign(\xi),&\text{if}\;\;|\xi|>\varepsilon;\\
\xi,&\text{if}\;\;|\xi|\leq\varepsilon.
\end{cases}
\end{equation}
Since $h=\phi\circ d_C$ and $\phi$ is even, we apply 
Example~\ref{ex:2} to conclude.
\end{proof}

The following is an extension of the Huber loss function 
\cite{Hube81}.

\begin{example}[abstract Huber function]
\label{ex:9}
\rm
Let $C$ be a nonempty closed convex subset of $\HH$, let
$\rho\in\RPP$, and set 
\begin{equation}
h\colon\HH\to\RR\colon x\mapsto
\begin{cases}
\rho d_C(x)-\Frac{\rho^2}{2}, &\text{if}\;\; 
d_C(x)>\rho;\\[+4mm]
\Frac{d_C(x)^2}{2},&\text{if}\;\; d_C(x)\leq\rho.
\end{cases}
\end{equation}
Let $\gamma\in\RPP$ and $x\in\HH$.
Then $h\colon\HH\to\RR$ is convex and Fr\'echet differentiable
with a nonexpansive gradient,
\begin{equation}
\label{e:kj801}
\nabla h(x)=
\begin{cases}
\dfrac{\rho}{d_C(x)}\big(x-\proj_Cx\big),
&\text{if}\;\;d_C(x)>\rho;\\
x-\proj_Cx,&\text{if}\;\;d_C(x)\leq\rho,
\end{cases}
\end{equation}
and 
\begin{equation}
\label{e:kj56}
\prox_{\gamma h}x=
\begin{cases}
x+\dfrac{\gamma\rho}{d_C(x)}(\proj_Cx-x),&\text{if}\;\;
d_C(x)>(\gamma+1)\rho;\\[4mm]
\dfrac{1}{\gamma+1}\big(x+\gamma\proj_Cx\big),
&\text{if}\;\;d_C(x)\leq(\gamma+1)\rho.
\end{cases}
\end{equation}
\end{example}
\begin{proof}
Let 
\begin{equation}
\label{e:huber}
\mathfrak{h}_{\rho}\colon\RR\to\RR\colon\xi\mapsto
\begin{cases}
\rho|\xi|-\Frac{\rho^2}{2},&\text{if}\;\;|\xi|>\rho;\\[+4mm]
\Frac{|\xi|^2}{2},&\text{if}\;\;|\xi|\leq\rho
\end{cases}
\end{equation}
be the standard Huber function with parameter $\rho$. Then
$\mathfrak{h}_{\rho}'$ is $1$-Lipschitzian. In addition, using the
expression of $\prox_{\gamma\mathfrak{h}_{\rho}}$ from
\cite[Example~24.9]{Livre1} and then Example~\ref{ex:2}, we obtain
\eqref{e:kj56}.
\end{proof}

\begin{example}
\label{ex:8}
\rm
Let $C$ be a nonempty closed convex subset of $\HH$, let
$\omega\in\RPP$, and set $\beta=\omega^2$ and
$h=\omega d_C-\ln(1+\omega d_C)$. Let $\gamma\in\RPP$ and 
let $x\in\HH$. Then $h\colon\HH\to\RR$ is convex and Fr\'echet
differentiable with a $\beta$-Lipschitzian gradient,
\begin{equation}
\label{e:kj800}
\nabla h(x)=
\begin{cases}
\dfrac{\omega^2}{1+\omega d_C(x)}\big(x-\proj_Cx\big),
&\text{if}\;\;x\notin C;\\
0,&\text{if}\;\;x\in C,
\end{cases}
\end{equation}
and 
\begin{multline}
\label{e:kj58}
\hskip -0mm
\prox_{\gamma h}x=\\
\hskip -12mm
\begin{cases}
\proj_Cx+\dfrac{\gamma\omega^2+1-\omega d_C(x)-
\sqrt{|\omega d_C(x)-\gamma\omega^2-1|^2+4\omega d_C(x)}}
{2\omega d_C(x)}(\proj_Cx-x),&\text{if}\;\:x\notin C;\\[4mm]
x,&\text{if}\;\:x\in C.
\end{cases}
\end{multline}
\end{example}
\begin{proof}
We apply Example~\ref{ex:2} with 
$\phi=\omega|\cdot|-\ln(1+\omega|\cdot|)$. Note that
$\phi'\colon\xi\mapsto\omega^2\xi/(1+\omega|\xi|)$ is 
$\omega^2$-Lipschitzian. Furthermore, we derive
$\prox_{\gamma\phi}$ by arguing as in 
\cite[Example~24.42]{Livre1} (where $\gamma=1$) and we then 
invoke \eqref{e:kj52} to get \eqref{e:kj58}.
\end{proof}

The following extension of Example~\ref{ex:2} involves a
composition with a linear operator. 

\begin{example}
\label{ex:7}
\rm
Let $\GG$ be a real Hilbert space and let $M\in\BL(\HH,\GG)$ be 
such that $MM^*=\theta\Id$ for some $\theta\in\RPP$.
Let $D$ be a nonempty closed convex subset of $\GG$, let
$\mu\in\RPP$, let $\phi\colon\RR\to\RR$ be even, convex, and
differentiable with a $\mu$-Lipschitzian derivative, and
set $h=\phi\circ d_D\circ M$ and $\beta=\mu\|M\|^2$. 
Let $\gamma\in\RPP$ and $x\in\HH$.
Then $h\colon\HH\to\RR$ is convex and Fr\'echet differentiable
with a $\beta$-Lipschitzian gradient,
\begin{equation}
\label{e:kj802}
\nabla h(x)=
\begin{cases}
\dfrac{\phi'\big(d_D(Mx)\big)}{d_D(Mx)}M^*\big(Mx-\proj_D(Mx)\big),
&\text{if}\;\;Mx\notin D;\\
0,&\text{if}\;\;Mx\in D,
\end{cases}
\end{equation}
and
\begin{equation}
\label{e:kj60}
\prox_{\gamma h}x=
\begin{cases}
x+\Frac{\theta^{-1}\big(d_D(Mx)-
\prox_{\gamma\theta\phi}d_D(Mx)\big)}{d_D(Mx)}
M^*\big(\proj_D(Mx)-Mx\big),&\text{if}\;\;Mx\notin D;\\
x,&\text{if}\;\;Mx\in D.
\end{cases}
\end{equation}
\end{example}
\begin{proof}
We have $h=(\phi\circ d_D)\circ M$ and therefore 
$\nabla h=M^*\circ\nabla(\phi\circ d_D)\circ M$. In turn, the
Lipschitz constant of $\nabla h$ is $\|M^*\|\mu\|M\|=\beta$, 
and we derive \eqref{e:kj802} from \eqref{e:kj80}.
Now set $g=\gamma\phi\circ d_D$. We derive from 
\cite[Proposition~24.14]{Livre1} that 
\begin{equation}
\label{e:kj61}
\prox_{\gamma h}x=\prox_{g\circ M}x=
x+\theta^{-1}M^*\big(\prox_{\theta g}(Mx)-Mx\big).
\end{equation}
We then obtain the expression for 
$\prox_{\theta g}=\prox_{\gamma\theta\phi\circ d_D}$ 
from \eqref{e:kj52}, which yields \eqref{e:kj60}.
\end{proof}

\begin{remark}
\label{r:7}
\rm
The condition $MM^*=\theta\Id$ used in
Example~\ref{ex:7} arises in particular in problems involving 
tight frame representations \cite{Chau07}. When it is not
satisfied, one can still deal with smooth functions of the type
$\phi\circ d_C\circ M$ in modern structured proximal splitting
techniques by activating $\prox_{\phi\circ d_C}$ and $M$
separately; see Propositions~\ref{p:2012}, \ref{p:fbp}, and 
\ref{p:kt} below and \cite{Livre1}. 
One can then invoke Example~\ref{ex:2} to compute the former. 
\end{remark}

\subsection{Integral functions}

\begin{example}
\label{ex:int}
\rm
Let $(\Omega,\mathcal{F},\mu)$ be a complete $\sigma$-finite
measure space, let $(\mathsf{H},\scal{\cdot}{\cdot}_\mathsf{H})$ be
a separable real Hilbert space, let $\mathsf{C}$ be a
closed convex subset of $\mathsf{H}$ such that $\mathsf{0}\in
\mathsf{C}$,
and let $\phi\colon\RR\to\RR$ be even, convex, and
differentiable with a $\beta$-Lipschitzian derivative.
Suppose that $\HH=L^2((\Omega,\mathcal{F},\mu);\mathsf{H})$,
and that $\mu(\Omega)<\pinf$ or $\phi(0)=0$. 
Set $h\colon\HH\rightarrow\RR\colon
x\mapsto\int_\Omega\phi(d_\mathsf{C}(x(\omega)))\mu(d\omega)$.
Let $\gamma\in\RPP$ and $x\in\HH$. 
Then $h\colon\HH\to\RR$ is convex and Fr\'echet differentiable
with a $\beta$-Lipschitzian gradient, and for 
$\mu$-almost every $w\in\Omega$,
\begin{equation}
\label{e:sonne}
\big(\nabla h(x)\big)(\omega)=
\begin{cases}
\Frac{\phi'(d_\mathsf{C}(x(\omega)))}{d_\mathsf{C}(x(\omega))}
\big(x(\omega)-\proj_\mathsf{C}x(\omega)\big),
&\text{if}\;\;x(\omega)\notin \mathsf{C};\\
0,&\text{if}\;\;x(\omega)\in\mathsf{C},
\end{cases}
\end{equation}
and
\begin{equation}
\label{e:proxint}
\big(\prox_{\gamma h}x\big)(\omega)=
\begin{cases}
\proj_\mathsf{C}x(\omega)+\Frac{\prox_{\gamma\phi}
d_\mathsf{C}(x(\omega))}{d_\mathsf{C}(x(\omega))}
\big(x(\omega)-\proj_\mathsf{C}x(\omega)\big),
&\text{if}\;\;x(\omega)\notin \mathsf{C};\\
x(\omega),&\text{if}\;\;x(\omega)\in \mathsf{C}.
\end{cases}
\end{equation}
\end{example}
\begin{proof}
Set $\varphi=\phi\circ d_\mathsf{C}$. As seen in
Example~\ref{ex:2}, $\varphi\colon\mathsf{H}\to\RR$ is convex and
Fr\'echet differentiable with a $\beta$-Lipschitzian gradient, 
\begin{equation}
\label{e:kj8-0}
(\forall\mathsf{x}\in\mathsf{H})\quad
\nabla \varphi(\mathsf{x})=
\begin{cases}
\dfrac{\phi'\big(d_\mathsf{C}(\mathsf{x})\big)}
{d_\mathsf{C}(\mathsf{x})}
\big(\mathsf{x}-\proj_\mathsf{C}\mathsf{x}\big),
&\text{if}\;\;\mathsf{x}\notin \mathsf{C};\\
0,&\text{if}\;\;\mathsf{x}\in \mathsf{C},
\end{cases}
\end{equation}
and
\begin{equation}
\label{e:kj85}
(\forall\mathsf{x}\in\mathsf{H})\quad
\prox_{\gamma\varphi}\mathsf{x}=
\begin{cases}
\proj_\mathsf{C}\mathsf{x}+
\Frac{\prox_{\gamma\phi}d_\mathsf{C}(\mathsf{x})}
{d_\mathsf{C}(\mathsf{x})}(\mathsf{x}-\proj_\mathsf{C}\mathsf{x}),
&\text{if}\;\;\mathsf{x}\notin \mathsf{C};\\
\mathsf{x},&\text{if}\;\;\mathsf{x}\in\mathsf{C}.
\end{cases}
\end{equation}
Since $\phi$ is convex and even, it is minimized by $0$
\cite[Proposition~11.7(i)]{Livre1}. Hence, 
$\phi'(0)=0$ and, since $\mathsf{0}\in\mathsf{C}$,
$\varphi(\mathsf{0})=\phi(0)$,
while \eqref{e:kj8-0} yields
$\nabla\varphi(\mathsf{0})=\mathsf{0}$. 
Consequently, by virtue of the descent lemma 
\cite[Theorem~18.15(iii)]{Livre1}, 
\begin{equation}
\label{e:brennt}
(\forall\mathsf{x}\in\mathsf{H})\quad
\varphi(\mathsf{x})\leq 
\varphi(\mathsf{0})+\scal{\mathsf{x}}
{\nabla\varphi(\mathsf{0})}_\mathsf{H}
+\frac{\beta}{2}\|\mathsf{x}\|_\mathsf{H}^2
=\varphi(\mathsf{0})+\frac{\beta}{2}\|\mathsf{x}\|_\mathsf{H}^2.
\end{equation}
In turn,
\begin{equation}
\label{e:87fr}
h(x)=\int_\Omega\varphi(x(\omega))\mu(d\omega)\leq
\varphi(\mathsf{0})\mu(\Omega)+\dfrac{\beta}{2}\|x\|^2<\pinf.
\end{equation}
On the other hand, 
\cite[Proposition~16.63(ii)]{Livre1} asserts that
$\nabla h(x)=(\nabla\varphi)\circ x$ $\mu$-a.e. which, 
combined with \eqref{e:kj8-0}, yields \eqref{e:sonne}.
Now let $y$ be in $\HH$. Then
\begin{align}
\|\nabla h(x)-\nabla h(y)\|^2&=
\int_\Omega\|(\nabla h(x))(\omega)-(\nabla
h(y))(\omega)\|_\mathsf{H}^2\mu(d\omega)\nonumber\\
&=\int_\Omega\|\nabla\varphi(x(\omega))-\nabla
\varphi(y(\omega))\|_\mathsf{H}^2\mu(d\omega)\nonumber\\
&\leq\beta^2\int_\Omega\|x(\omega)-y(\omega)\|_\mathsf{H}^2
\mu(d\omega)\nonumber\\
&=\beta^2\|x-y\|^2,
\end{align}
which shows that $\nabla h$ is $\beta$-Lipschitzian. Finally, we
apply \cite[Proposition~24.13]{Livre1} to derive \eqref{e:proxint}
from \eqref{e:kj85}.
\end{proof}

\begin{remark}
\label{r:gaussTV}
\rm
Let $\Omega$ be a nonempty bounded smooth open subset of $\RR^2$, 
let $\mathsf{H}=\RR^2$, let $\mu$ be the Lebesgue measure, and
suppose that $\mathsf{C}=\{\mathsf{0}\}$ in Example~\ref{ex:int}. 
Furthermore, let $\mathfrak{h}_\rho$ be the Huber function of 
\eqref{e:huber} and, for every $x\in H^1_0(\Omega)$, let $Dx$ be
the gradient of $x$. Then the function 
\begin{equation}
h\circ D\colon H^1_0(\Omega)\to\RR\colon x\mapsto
\int_\Omega \mathfrak{h}_\rho(\|D x(\omega)\|_2)d\omega
\end{equation}
can be found in \cite{Hint06,Keel03,Louc13,Schn94} and it is called 
the Gauss-TV (or TV-Huber) function.
\end{remark}

\subsection{Functionals involving orthonormal decompositions}
We first revisit a construction proposed in
\cite{Smms05}; see also \cite{Save18,Daub04,Demo02} for
special cases.

\begin{example}
\label{ex:12}
\rm
Suppose that $\HH$ is separable and that $\emp\neq\bK\subset\NN$, 
and let $(e_k)_{k\in\bK}$ be an orthonormal basis of $\HH$.
For every $k\in\bK$, let $\beta_k\in\RPP$ and let 
$\phi_k\colon\RR\to\RR$ be a differentiable convex function such 
that $\phi_k\geq\phi_k(0)=0$ and $\phi_k'$ is
$\beta_k$-Lipschitzian. Suppose that
$\beta=\sup_{k\in\bK}\beta_k<\pinf$ and define
$(\forall x\in\HH)$ $h(x)=\sum_{k\in\bK}\phi_k(\scal{x}{e_k})$. 
Let $\gamma\in\RPP$.
Then $h\colon\HH\to\RR$ is convex and Fr\'echet differentiable
with a $\beta$-Lipschitzian gradient,
\begin{equation}
(\forall x\in\HH)\quad\nabla h(x)=
\Sum_{k\in\bK}\phi_k'(\scal{x}{e_k})e_k,
\end{equation}
and 
\begin{equation}
\label{e:fin}
(\forall x\in\HH)\quad\prox_{\gamma h}x=
\sum_{k\in\bK}\big(\prox_{\gamma\phi_k}\scal{x}{e_k}\big)e_k.
\end{equation}
\end{example}
\begin{proof}
The identity \eqref{e:fin} is established in \cite{Smms05}.
The frame analysis operator is
\begin{equation}
\label{e:K}
F\colon\HH\to\ell^2(\bK)\colon x\mapsto
(\scal{x}{e_k})_{k\in\bK}
\end{equation}
and its adjoint is the frame synthesis operator
\begin{equation}
\label{e:Ka}
F^*\colon\ell^2(\bK)\to\HH\colon(\xi_k)_{k\in\bK}
\mapsto\sum_{k\in\bK}\xi_ke_k.
\end{equation}
Now denote by $\mathsf{x}=(\xi_k)_{k\in\bK}$ a generic element in
$\ell^2(\bK)$ and define
\begin{equation}
\varphi\colon\ell^2(\bK)\to\RX\colon \mathsf{x}\mapsto
\sum_{k\in\bK}\phi_k(\xi_k). 
\end{equation}
Then $h=\varphi\circ F$. Since all the functions
$(\phi_k)_{k\in\bK}$ are minimized at
$0$, we have $(\forall k\in\bK)$ $\phi_k'(0)=0=\phi_k(0)$. In
turn, we derive from the descent lemma
\cite[Theorem~18.15(iii)]{Livre1} that
\begin{equation}
\label{e:caba}
(\forall k\in\bK)(\forall\xi_k\in\RR)\quad
\phi_k(\xi_k)\leq\phi_k(0)+(\xi_k-0)\phi_k'(0)
+\frac{\beta_k}{2}|0-\xi_k|^2=\frac{\beta}{2}|\xi_k|^2.
\end{equation}
As a result, 
\begin{equation}
\label{e:vage}
\big(\forall\mathsf{x}\in\ell^2(\bK)\big)\quad
\varphi(\mathsf{x})=\sum_{k\in\bK}\phi_k(\xi_k)\leq
\frac{\beta}{2}\sum_{k\in\bK}|\xi_k|^2=\frac{\beta}{2}
\|\mathsf{x}\|^2
\end{equation}
and, therefore, $\varphi\colon\ell^2(\bK)\to\RR$. In addition,
\begin{equation}
\label{e:ret}
(\forall k\in\bK)(\forall\xi_k\in\RR)\quad
|\phi_k'(\xi_k)|^2=|\phi_k'(\xi_k)-\phi_k'(0)|^2\leq
\beta_k^2|\xi_k-0|^2\leq\beta^2|\xi_k|^2,
\end{equation}
which yields 
\begin{equation}
\label{e:sau}
\big(\forall\mathsf{x}\in\ell^2(\bK)\big)\quad
\sum_{k\in\bK}|\phi_k'(\xi_k)|^2\leq
\beta^2\sum_{k\in\bK}|\xi_k|^2=\beta^2
\|\mathsf{x}\|^2<\pinf,
\end{equation}
and allows us to conclude that $\varphi$ is differentiable with 
$(\forall\mathsf{x}\in\ell^2(\bK))$ 
$\nabla\varphi(\mathsf{x})=(\phi_k'(\xi_k))_{k\in\bK}$. Moreover,
\begin{align}
\big(\forall\mathsf{x}\in\ell^2(\bK)\big)
\big(\forall\mathsf{y}\in\ell^2(\bK)\big)\quad
\|\nabla\varphi(\mathsf{x})-\nabla\varphi(\mathsf{y})\|^2
&=\sum_{k\in\bK}|\phi_k'(\xi_k)-\phi_k'(\eta_k)|^2
\nonumber\\
&\leq\sum_{k\in\bK}\beta_k^2|\xi_k-\eta_k|^2\nonumber\\
&\leq\beta^2\|\mathsf{x}-\mathsf{y}\|^2.
\end{align}
Altogether, $\varphi\colon\ell^2(\bK)\to\RR$ is convex and
differentiable with a $\beta$-Lipschitzian gradient. Hence
$h=\varphi\circ F\colon\HH\to\RR$ is convex and differentiable,
with $\nabla h=F^*\circ\nabla\varphi\circ F$. Furthermore, 
since $F$ and $F^*$ are isometries, 
\begin{align}
(\forall x\in\HH)(\forall y\in\HH)\quad
\|\nabla h(x)-\nabla h(y)\|
&=\|F^*(\nabla\varphi(Fx))-F^*(\nabla\varphi(Fy))\|\nonumber\\
&=\|\nabla\varphi(Fx)-\nabla\varphi(Fy)\|\nonumber\\
&\leq\beta\|Fx-Fy\|\nonumber\\
&=\beta\|x-y\|,
\end{align}
which shows that $\nabla h$ is $\beta$-Lipschitzian.
\end{proof}

\begin{remark}
\label{r:nd}
\rm
It is clear from the proof of Example~\ref{ex:12} that if the 
condition
\begin{equation}
\label{e:c6}
\phi_k\geq\phi_k(0)=0
\end{equation}
is not satisfied for a finite number of indices $k\in\bK$, the
results remain valid. In particular, if $\HH$ is
finite-dimensional, \eqref{e:c6} is not required. An example of
a smooth convex function $\phi_k\colon\RR\to\RR$ with an explicit
proximity operator is $\phi_k=d_{C_k}^2/2$, where $C_k$ is a
nonempty closed interval in $\RR$ (set $\phi=|\cdot|^2/2$ in
Example~\ref{ex:2}). In
this case, \eqref{e:c6} holds if and only if $0\in C_k$. If
$C_k=\left[1,\pinf\right[$, $\phi_k$ is the squared hinge loss 
\cite{Stei05}; if $C_k=\left[-\varepsilon,\varepsilon\right]$ for
some $\varepsilon\in\RPP$, $\phi_k$ is the smooth Vapnik
insensitive loss \cite{Arav13} (see also Example~\ref{ex:6bis}).
Specializing Examples~\ref{ex:2}, \ref{ex:9}, and \ref{ex:8}, as
well as \eqref{e:kj63} to $\HH=\RR$ provides further examples of
functions $\phi_k=\psi_k\circ d_{C_k}$ of interest. For instance, 
taking $\psi_k$ to be the Huber function \eqref{e:huber} and
$C_k=\left[1,\pinf\right[$ yields the modified Huber function
of \cite{Zhan04}.
\end{remark}

\begin{example}
\label{ex:npsi2}
\rm
Let $M$ be a strictly positive integer and let 
$\emp\neq\bK\subset\NN$. For every
$i\in\{1,\ldots,M\}$, suppose that $\HH_i$ is a separable real
Hilbert space with orthonormal basis $(e_{i,k})_{k\in\bK}$.
Let $\beta\in\RPP$, and let $\phi\colon\RR\to\RR$ be an even 
differentiable convex function such 
that $\phi\geq\phi(0)=0$ and $\phi'$ is
$\beta$-Lipschitzian. For every 
$(x_1,\ldots,x_M)$ in $\HH_1\oplus\cdots\oplus\HH_M$, set
$h(x_1,\ldots,x_M)=\sum_{k\in\bK}
\phi(\|(\scal{x_1}{e_{1,k}},\ldots,\scal{x_M}{e_{M,k}})\|_2)$.
Let $\gamma\in\RPP$, let 
$(x_1,\ldots,x_M)\in\HH_1\oplus\cdots\oplus\HH_M$, and set,
for every $k\in\bK$,
\begin{equation}
\alpha_k=
\begin{cases}
\dfrac{\phi'(\|(\scal{x_1}{e_{1,k}},\ldots,\scal{x_M}{e_{M,k}})\|_2)}
{\|(\scal{x_1}{e_{1,k}},\ldots,\scal{x_M}{e_{M,k}})\|_2},
&\text{if}\;\;\displaystyle\max_{1\leq i\leq M}
\abscal{x_i}{e_{i,k}}>0;\\ 
0,&\text{otherwise,}
\end{cases}
\end{equation}
and
\begin{equation}
\delta_k=
\begin{cases}
\dfrac{\prox_{\gamma\phi}
(\|(\scal{x_1}{e_{1,k}},\ldots,\scal{x_M}{e_{M,k}})\|_2)}
{\|(\scal{x_1}{e_{1,k}},\ldots,\scal{x_M}{e_{M,k}})\|_2},
&\text{if}\;\;\displaystyle\max_{1\leq i\leq M}
\abscal{x_i}{e_{i,k}}>0;\\[3mm]
1,&\text{otherwise.}
\end{cases}
\end{equation}
Then $h\colon\HH_1\oplus\cdots\oplus\HH_M\to\RR$ is convex and
Fr\'echet 
differentiable with a $\beta$-Lipschitzian gradient,
\begin{equation}
\label{e:13455}
\nabla h(x_1,\ldots,x_M)=\Bigg(
\Sum_{k\in\bK}\alpha_k\scal{x_1}{e_{1,k}}e_{1,k},\ldots,
\Sum_{k\in\bK}\alpha_k\scal{x_M}{e_{M,k}}e_{M,k}\Bigg),
\end{equation}
and 
\begin{equation} 
\label{e:ex140} 
\prox_{\gamma h}(x_1,\ldots,x_M)= 
\Bigg( \Sum_{k\in\bK}\delta_k\scal{x_1}{e_{1,k}}e_{1,k},\ldots,
\Sum_{k\in\bK}\delta_k\scal{x_M}{e_{M,k}}e_{M,k}\Bigg).
\end{equation}
\end{example}
\begin{proof}
Arguing as in \eqref{e:caba}--\eqref{e:vage}, we obtain that $h$
is real-valued and, arguing as in \eqref{e:ret}, we obtain
$\sup_{k\in\bK}\alpha_k\leq\beta$. Hence
$(\forall i\in\{1,\ldots,M\})$
$(\alpha_k\scal{x_i}{e_{i,k}})_{k\in\bK}\in\ell^2(\bK)$ and
$\sum_{k\in\bK}\alpha_k\scal{x_i}{e_{i,k}}e_{i,k}\in\HH_i$. 
Likewise, since $\phi$ is even, 
$\phi\geq\phi(0)=0$ \cite[Proposition~11.7(i)]{Livre1}, 
hence $\prox_{\gamma\phi}0=0$,
and thus $\sup_{k\in\bK}\delta_k\leq 1$. It therefore follows that
$(\forall i\in\{1,\ldots,M\})$
$(\delta_k\scal{x_i}{e_{i,k}})_{k\in\bK}\in\ell^2(\bK)$ and
$\sum_{k\in\bK}\delta_k\scal{x_i}{e_{i,k}}e_{i,k}\in\HH_i$. 
Now let $\ell^2(\bK;\RR^M)=\bigoplus_{k\in\bK}\RR^M$ be the space of
square summable sequences with entries in $\RR^M$, and set
\begin{equation}
\begin{cases}
\varphi\colon\RR^M\rightarrow\RR\colon(\xi_1,\ldots,\xi_M)\mapsto
\phi(\|(\xi_1,\ldots,\xi_M)\|_2)\\
U\colon\HH_1\oplus\cdots\oplus\HH_M
\to\ell^2(\bK;\RR^M)\colon(x_1,\ldots,x_M)\mapsto
\big((\scal{x_1}{e_{1,k}},\ldots,
\scal{x_M}{e_{M,k}})\big)_{k\in\bK}.
\end{cases}
\end{equation}
Then $U$ is a bijective isometry,
\begin{equation}
\begin{array}{lcl}
U^{-1}=U^*\colon&\ell^2(\bK;\RR^M)&\to
\;\HH_1\oplus\cdots\oplus\HH_M\\
&\big((\mu_{1,k},\ldots,\mu_{M,k})\big)_{k\in\bK}&\mapsto
\;\Bigg(\Sum_{k\in\bK}\mu_{1,k}e_{1,k},\ldots,
\Sum_{k\in\bK}\mu_{M,k}e_{M,k}\Bigg),
\end{array}
\end{equation}
and $h=\big(\bigoplus_{k\in\bK}\varphi\big)\circ U$. In turn,
\begin{align}
\nabla h
&=U^*\circ\nabla\Bigg(\bigoplus_{k\in\bK}\varphi\Bigg)\circ U
\nonumber\\
&=U^*\circ\nabla\Bigg(\bigoplus_{k\in\bK}\big(\phi\circ
\|\cdot\|_2\big)\Bigg)\circ U
\nonumber\\
&=U^*\circ\Bigg(\underset{k\in\bK}{\bigotimes}\nabla\big(\phi\circ
\|\cdot\|_2\big)\Bigg)\circ U,
\end{align}
which yields \eqref{e:13455} by applying Example~\ref{ex:2} with
$C=\{0\}$. On the other 
hand, using \cite[Proposition~24.14]{Livre1}, we obtain
\begin{equation}
\prox_{\gamma h}
=\prox_{\gamma(\bigoplus_{k\in\bK}\varphi)\circ U}
=\big(U^*\circ\prox_{\gamma\bigoplus_{k\in\bK}\varphi}
\circ U\big)
=U^*\circ\Bigg(\underset{k\in\bK}{\bigotimes}
\prox_{\gamma\phi\circ\|\cdot\|_2}\Bigg)\circ U,
\end{equation}
which yields \eqref{e:ex140} by applying Example~\ref{ex:2} 
with $C=\{0\}$.
\end{proof}

\subsection{Function involving explicit infimal convolutions}

If $h$ is explicitly constructed in terms of a function
$g\in\Gamma_0(\HH)$ as $h=g\infconv(\beta q)$ for some
$\beta\in\RPP$, then, as in \eqref{e:jjm69}, we obtain
\begin{equation}
\label{e:miles}
\nabla h=\beta\big(\Id-\prox_{g/\beta}\big)
\quad\text{and}\quad
\prox_{\gamma h}x=\Id+\dfrac{\gamma\beta}{\gamma\beta+1}
\big(\prox_{(\gamma\beta+1)g/\beta}-\Id\big).
\end{equation}
In the same spirit, we have the following construction.

\begin{example}
\label{ex:22}
\rm
Let $\varphi\in\Gamma_0(\HH)$, let $\beta\in\RPP$, and 
define $h=\beta q-(\beta q)\infconv\varphi$.
Let $\gamma\in\RPP$ and $x\in\HH$.
Then $h\colon\HH\to\RR$ is convex and Fr\'echet differentiable
with a $\beta$-Lipschitzian gradient, 
\begin{equation}
\label{e:kj62}
\nabla h(x)=\beta\prox_{\varphi/\beta}x,
\quad\text{and}\quad
\prox_{\gamma h}x=
x-\beta\gamma\prox_{\frac{\varphi}{\beta(1+\beta\gamma)}}
\bigg(\frac{1}{1+\beta\gamma}x\bigg).
\end{equation}
\end{example}
\begin{proof}
It follows from \cite[Proposition~12.30]{Livre1} that $\nabla h=
\beta(\Id-(\Id-\prox_{\varphi/\beta}))=\beta\prox_{\varphi/\beta}$,
which is $\beta$-Lipschitzian since $\prox_{\varphi/\beta}$ is 
nonexpansive \cite[Proposition~12.28]{Livre1}. Finally 
the expression of $\prox_{\gamma h}x$ follows from
\cite[Proposition~24.8(viii)]{Livre1}.
\end{proof}

\begin{remark}
\label{r:23}
\rm
Let $C$ be a nonempty closed convex subset of $\HH$, let 
$\beta\in\RPP$, let $\gamma\in\RPP$, and let $x\in\HH$. Set
$\varphi=\iota_C$ in Example~\ref{ex:22}.
Then we derive that $h=\beta(q-d_C^2/2)\colon\HH\to\RR$ is convex
and Fr\'echet differentiable with a $\beta$-Lipschitzian gradient,
and that
\begin{equation}
\label{e:kj63}
\nabla h(x)=\beta\proj_{C}x
\quad\text{and}\quad
\prox_{\gamma h}x=x-\beta\gamma\proj_{C}
\bigg(\frac{1}{1+\beta\gamma}x\bigg).
\end{equation}
For $\beta=1$, $h$ is the generalized Huber function of
\cite[Example~3.2]{Svva17}, that is,
\begin{equation}
(\forall x\in\HH)\quad h(x)=
\begin{cases}
\scal{x}{\proj_Cx}-\dfrac{\|\proj_Cx\|^2}{2},
&\text{if}\;\:x\notin C;\\[3mm]
\dfrac{\|x\|^2}{2},
&\text{if}\;\:x\in C.
\end{cases}
\end{equation}
If $\HH=\RR$ and $C=[-\rho,\rho]$ for some $\rho\in\RPP$, $h$
reduces to the standard Huber function \eqref{e:huber}.
\end{remark}

\section{Splitting algorithms}
\label{sec:3}

In this section, we review several proximal splitting algorithms
which are relevant to our discussion and will be used in the 
numerical experiments (see \cite{Livre1} for the supporting
theory). We start with the forward-backward splitting algorithm.

\begin{Algorithm}[forward-backward]
\label{al:fb}
Let $\beta\in\left]0,\pinf\right[\,$, 
let $f\in\Gamma_0(\HH)$, let $h\colon\HH\to\RR$ 
be convex and differentiable with 
a $\beta$-Lipschitzian gradient, let $(\gamma_n)_{n\in\NN}$ be a
sequence in $\left]0,2/\beta\right[$ such that
$0<\inf_{n\in\NN}\gamma_n\leq\sup_{n\in\NN}\gamma_n<2/\beta$.
Let $x_0\in\HH$ and iterate 
\begin{equation}
\label{e:fb}
\begin{array}{l}
\text{for}\;n=0,1,\ldots\\
\left\lfloor
\begin{array}{l}
y_n=x_n-\gamma_n\nabla h(x_n)\\
x_{n+1}=\prox_{\gamma_n f}y_n.
\end{array}
\right.\\
\end{array}
\end{equation}
\end{Algorithm}

\begin{proposition}{\rm\cite{Save18}}
\label{p:fb}
Let $(x_n)_{n\in\NN}$ be a sequence generated by
Algorithm~\ref{al:fb} and suppose that $\Argmin(f+h)\neq\emp$. 
Then there exists $x\in\Argmin(f+h)$ such that
$x_n\weakly x$. In addition $(f+h)(x_n)-\inf (f+h)(\HH)=o(1/n)$.
\end{proposition}

Inertial variants of the above method have been popularized by
\cite{Beck09}. They require additional storage capabilities but
have been shown to be advantageous in terms of convergence speed
in certain situations. The following implementation proposed in
\cite{Cham15} guarantees convergence of the iterates.

\begin{Algorithm}[inertial forward-backward]
\label{al:cd}
Let $\beta\in\left]0,\pinf\right[\,$, 
let $f\in\Gamma_0(\HH)$, let $h\colon\HH\to\RR$ 
be convex and differentiable with 
a $\beta$-Lipschitzian gradient, let 
$\gamma\in\left]0,1/\beta\right]$, 
and let $\alpha\in\left]2,\pinf\right[$.
Let $x_0=x_{-1}\in\HH$ and iterate 
\begin{equation}
\label{e:cd1}
\begin{array}{l}
\text{for}\;n=0,1,\ldots\\
\left\lfloor
\begin{array}{l}
z_n=x_n+\dfrac{n-1}{n+\alpha}(x_n-x_{n-1})\\[3mm]
y_n=z_n-\gamma\nabla h(z_n)\\
x_{n+1}=\prox_{\gamma f}y_n.
\end{array}
\right.\\
\end{array}
\end{equation}
\end{Algorithm}

\begin{proposition}{\rm\cite{Cham15}}
\label{p:cd}
Let $(x_n)_{n\in\NN}$ be a sequence generated by
Algorithm~\ref{al:cd} and suppose that $\Argmin(f+h)\neq\emp$. 
Then there exists $x\in\Argmin(f+h)$ such that
$x_n\weakly x$.
In addition $(f+h)(x_n)-\inf (f+h)(\HH)=O(1/n^2)$.
\end{proposition}

The next algorithm does not require smoothness of any of the
functions.

\begin{Algorithm}[Douglas-Rachford]
\label{al:dr}
Let $f$ and $g$ be functions in $\Gamma_0(\HH)$
such that $0\in\sri(\dom g-\dom f)$, 
let $\gamma\in\left]0,\pinf\right[$,
and let $(\lambda_n)_{n\in\NN}$ be sequence in $\left[0,2\right]$
such that $\sum_{n\in\NN}\lambda_n(2-\lambda_n)=\pinf$.
Let $y_0\in\HH$ and iterate 
\begin{equation}
\label{e:dr1}
\begin{array}{l}
\text{for}\;n=0,1,\ldots\\
\left\lfloor
\begin{array}{l}
z_n=\prox_{\gamma g}y_n\\
x_n=\prox_{\gamma f}(2z_n-y_n)\\
y_{n+1}=y_n+\lambda_n(x_n-z_n).
\end{array}
\right.\\
\end{array}
\end{equation}
\end{Algorithm}

\begin{proposition}
{\rm\cite{Jsts07}}
\label{p:dr}
Let $(x_n)_{n\in\NN}$ be a sequence generated by
Algorithm~\ref{al:dr} and suppose that $\Argmin(f+g)\neq\emp$. 
Then there exists $x\in\Argmin(f+g)$ such that
$x_n\weakly x$.
\end{proposition}

Although this feature will not be used in Section~\ref{sec:4}, it
should be noted that the forward-backward \cite{Save18},
inertial forward-backward \cite{Aujo15}, and 
Douglas-Rachford \cite{Opti04} algorithms tolerate errors in 
the implementations of the proximity operators. 
The next three algorithms are specifically tailored to handle
Problem~\ref{prob:1}. Although they also compute dual solutions,
for brevity, we present only the primal convergence result for
the error-free, unrelaxed formulations of these algorithms.
The first one is known as the primal-dual
forward-backward-forward algorithm; see \cite{Svva12} for details
and \cite{Bcvu13} for a variable metric version.

\begin{Algorithm}
\label{al:2012}
Consider the setting of Problem~\ref{prob:1}. Set 
$\beta=\sqrt{\sum_{i\in I}\|L_i\|^2}+
\sum_{j\in J}\mu_j\|L_j\|^2$, 
let $\varepsilon\in\left]0,1/(\beta+1)\right[$, 
let $(\gamma_n)_{n\in\NN}$ be a sequence in 
$\left[\varepsilon,(1-\varepsilon)/\beta\right]$, and 
let $(\forall i\in I)$ $v^*_{i,0}\in\GG_i$. Let $x_0\in\HH$ and
iterate 
\begin{equation}
\label{e:blackpagepart3}
\hskip -1mm
\begin{array}{l}
\text{for}\;n=0,1,\ldots\\
\left\lfloor
\begin{array}{l}
y_{1,n}=x_n-\gamma_n\Big(
\sum_{i\in I}L_i^*v^*_{i,n}+
\sum_{j\in J}L_j^*\big(\nabla h_j(L_jx_n)\big)\Big)\\
p_{1,n}=\prox_{\gamma_n f}\,y_{1,n}\\
\text{for every}\;i\in I\\
\left\lfloor
\begin{array}{l}
y_{2,i,n}=v^*_{i,n}+\gamma_nL_ix_n\\
p_{2,i,n}=y_{2,i,n}-\gamma_n
\prox_{g_i/\gamma_n}(y_{2,i,n}/\gamma_n)\\
q_{2,i,n}=p_{2,i,n}+\gamma_nL_ip_{1,n}\\
v^*_{i,n+1}=v^*_{i,n}-y_{2,i,n}+q_{2,i,n}
\end{array}
\right.\\[1mm]
q_{1,n}=p_{1,n}-\gamma_n\Big(
\sum_{i\in I}L_i^*p_{2,i,n}+
\sum_{j\in J}L_j^*\big(\nabla h_j(L_jp_{1,n})\big)\Big)\\
x_{n+1}=x_n-y_{1,n}+q_{1,n}.
\end{array}
\right.\\
\end{array}
\end{equation}
\end{Algorithm}

\begin{proposition}{\rm\cite{Svva12}}
\label{p:2012}
Let $(x_n)_{n\in\NN}$ be a sequence generated by
Algorithm~\ref{al:2012}.
Then there exists a solution $x$ to
\eqref{e:primal} such that $x_n\weakly x$.
\end{proposition}

The following algorithm is an implementation of the
forward-backward algorithm in a renormed primal-dual space (see
\cite{Cham11,Icip14,Cond13,Xiao12,Bang13} for special cases and
variants, and \cite{Opti14} for a more general variable metric
version).

\begin{Algorithm}
\label{al:fbp}
Consider the setting of Problem~\ref{prob:1}
and let $(\tau_n)_{n\in\NN}$ be a sequence in
$\RPP$ such that $(\forall n\in\NN)$ $\tau_{n+1}\geq\tau_n$.
For every $i\in I$, let $v^*_{i,0}\in\GG_i$, let 
$(\sigma_{i,n})_{n\in \NN}$ be a sequence in $\RPP$ such
that $(\forall n\in\NN)$ $\sigma_{i,n+1}\geq\sigma_{i,n}$.
Suppose that
\begin{equation}
\label{e:z12}
\sup_{n\in\NN}\Bigg(\sqrt{\tau_n\sum_{i\in I}
\sigma_{i,n}\|L_i\|^2}+
\frac{1}{2}\max\Big\{\tau_n,\max_{i\in I}\sigma_{i,n}\Big\}
\sum_{j\in J}\mu_j\|L_j\|^2\,\Bigg)<1.
\end{equation}
Let $x_0\in\HH$ and iterate 
\begin{equation}
\label{e:bv}
\begin{array}{l}
\text{for}\;n=0,1,\ldots\\
\left\lfloor
\begin{array}{l}
y_{1,n}=x_n-\tau_n\Big(\sum_{i\in I}L_i^*v^*_{i,n}+
\sum_{j\in J}L_j^*\big(\nabla h_j(L_jx_n)\big)\Big)\\
x_{n+1}=\prox_{\tau_n f}\,y_{1,n}\\
z_n=2x_{n+1}-x_n\\
\text{for every}\;i\in I\\
\left\lfloor
\begin{array}{l}
y_{2,i,n}=v^*_{i,n}+\sigma_{i,n}(L_iz_n)\\
v^*_{i,n+1}=y_{2,i,n}-\sigma_{i,n}\prox_{g_i/\sigma_{i,n}}
(y_{2,i,n}/\sigma_{i,n}).
\end{array}
\right.\\[2mm]
\end{array}
\right.\\[2mm]
\end{array}
\end{equation}
\end{Algorithm}

\begin{proposition}{\rm\cite{Opti14}}
\label{p:fbp}
Let $(x_n)_{n\in\NN}$ be a sequence generated by
Algorithm~\ref{al:fbp}.
Then there exists a solution $x$ to
\eqref{e:primal} such that $x_n\weakly x$.
\end{proposition}

The next algorithm, which was first proposed in \cite{Siop14}
in the case when $J=\emp$, was extended in \cite{MaPr18} 
to a block-coordinate and block-iterative asynchronous method.
The following version, which explicitly exploits smooth functions,
is proposed in \cite{Ecks18}.

\begin{Algorithm}
\label{al:kt}
Consider the setting of Problem~\ref{prob:1}
and let $(\gamma_n)_{n\in\NN}$ be a sequence in $\RPP$ such that 
$0<\inf_{n\in\NN}\gamma_n\leq\sup_{n\in\NN}\gamma_n<\pinf$.
For every $k\in I\cup J$, let $v^*_{k,0}\in\GG_k$, let 
$(\mu_{k,n})_{n\in \NN}$ be a sequence in $\RPP$ such that 
$0<\inf_{n\in\NN}\mu_{k,n}\leq\sup_{n\in\NN}\mu_{k,n}<\pinf$, and
let $(\lambda_n)_{n\in\NN}$ be a sequence in $]0,2[$ such that
$0<\inf_{n\in\NN}\lambda_n\leq\sup_{n\in\NN}\lambda_n<2$.
Let $x_0\in\HH$ and iterate
\begin{equation}
\label{e:kt}
\begin{array}{l}
\text{for}\;n=0,1,\ldots\\
\left\lfloor
\begin{array}{l}
l^*_n=\sum_{k\in I\cup J}L_{k}^*v_{k,n}^*\\
a_n=\prox_{\gamma_n f}(x_{n}-\gamma_nl^*_n)\\
a^*_n=\gamma_n^{-1}(x_n-a_n)-l^*_n\\
\text{for}\;i\in I\\
\left\lfloor
\begin{array}{l}
l_{i,n}=L_{i}x_{n}\\
b_{i,n}=\prox_{\mu_{i,n} g_i}\big(l_{i,n}+\mu_{i,n}v_{i,n}^*
\big)\\
b^*_{i,n}= v^*_{i,n}+\mu_{i,n}^{-1}(l_{i,n}-b_{i,n})\\
t_{i,n}=b_{i,n}-L_ia_n
\end{array}
\right.\\[1mm]
\text{for}\;j\in J\\
\left\lfloor
\begin{array}{l}
l_{j,n}=L_{j}x_{n}\\
b_{j,n}=l_{j,n}-\mu_{j,n}\big(\nabla h_j(l_{j,n})-v^*_{j,n}\big)\\
b^*_{j,n}=\nabla h_j(b_{j,n})\\
t_{j,n}=b_{j,n}-L_ja_n
\end{array}
\right.\\[1mm]
t^*_n=a^*_n+\sum_{k\in I\cup J}L^*_{k}b^*_{k,n}\\
\tau_n=\|t^*_n\|^2+\sum_{k\in I\cup J}\|t_{k,n}\|^2\\
\text{if}\;\tau_n=0\\
\left\lfloor
\begin{array}{l}
x=a_{n}\\
\text{terminate}.
\end{array}
\right.\\
\text{if}\;\tau_n>0\\
\left\lfloor
\begin{array}{l}
\theta_n=\dfrac{\lambda_n}{\tau_n}
\max\Big\{0,\scal{x_n}{t^*_n}-\scal{a_n}{a_n^*}+\sum_{k\in I\cup J}
\big(\scal{t_{k,n}}{v^*_{k,n}}-\scal{b_{k,n}}
{b^*_{k,n}}\big)\Big\}
\\
x_{n+1}=x_{n}-\theta_n t^*_{n}\\
\text{for}\;k\in I\cup J\\
\left\lfloor
\begin{array}{l}
v^*_{k,n+1}=v^*_{k,n}-\theta_n t_{k,n}.
\end{array}
\right.\\
\end{array}
\right.\\
\end{array}
\right.\\[8mm]
\end{array}
\end{equation}
\end{Algorithm}

\begin{proposition}{\rm\cite{Ecks18}}
\label{p:kt}
Either Algorithm~\ref{al:kt} terminates at a solution $x$ to
\eqref{e:primal} in a finite number of iterations, or it 
generates an infinite sequence $(x_n)_{n\in\NN}$ which converges
weakly to a solution to \eqref{e:primal}.
\end{proposition}

\section{Applications and numerical illustrations}
\label{sec:4}
We illustrate the viewpoint formulated in the Introduction, which
suggests that it may be computationally advantageous to activate
smooth functions proximally in certain instances. 
\begin{figure}[ht!]
\centering
\begin{tabular}{c@{}c@{}c@{}}
\includegraphics[width=4.7cm]{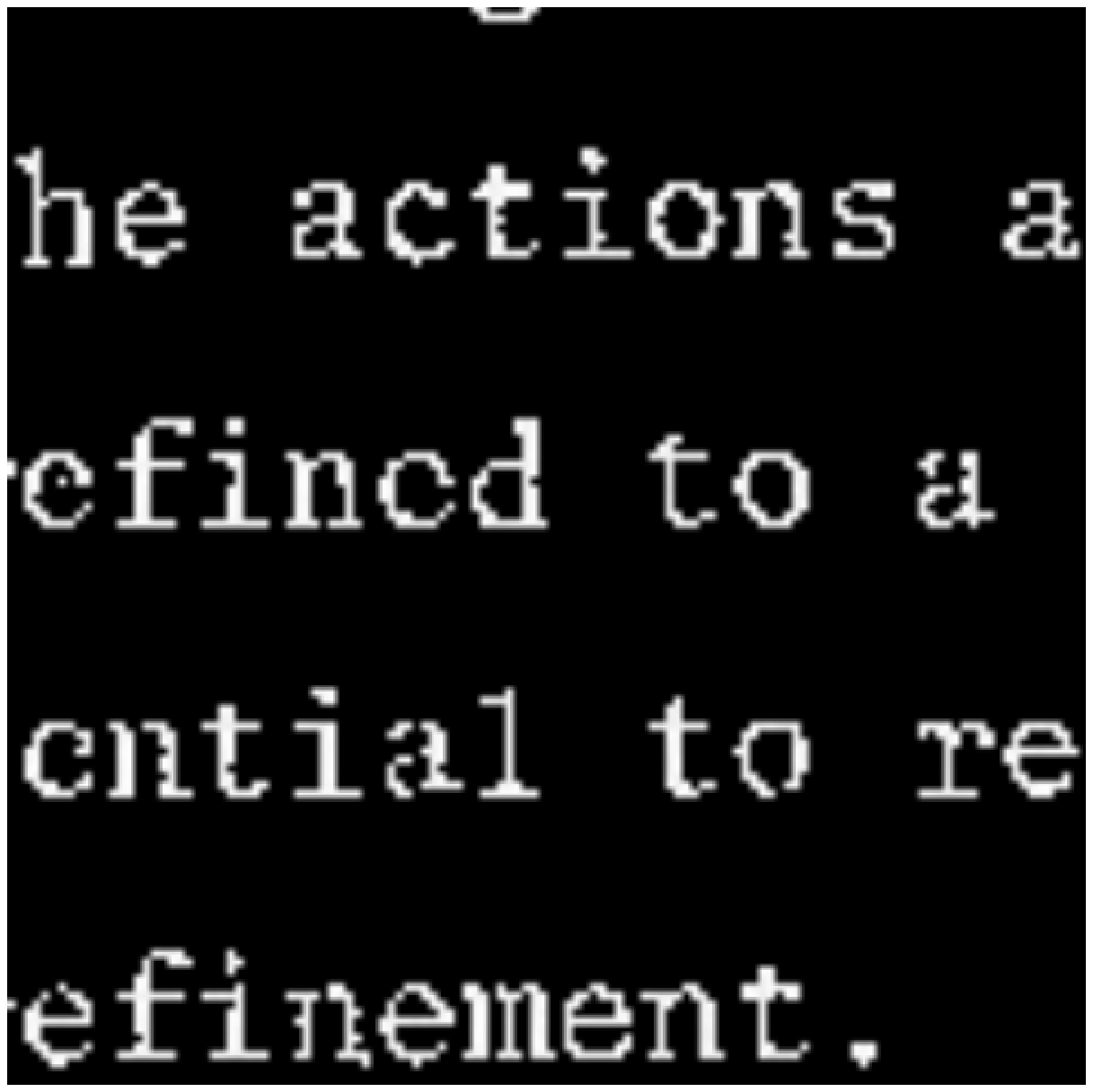}&
\hspace{0.4cm}
\includegraphics[width=4.7cm]{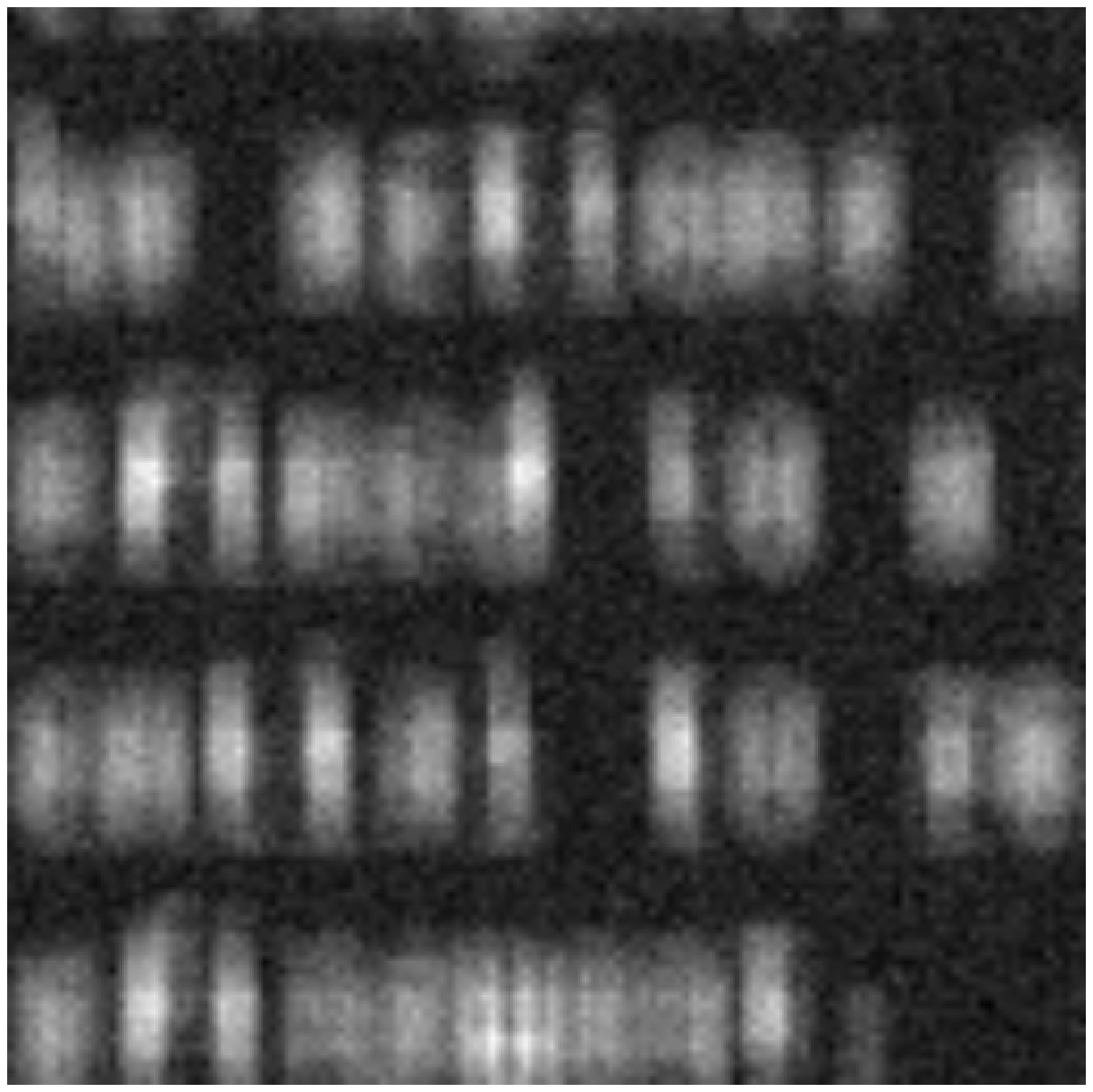}&
\hspace{0.4cm}
\includegraphics[width=4.7cm]{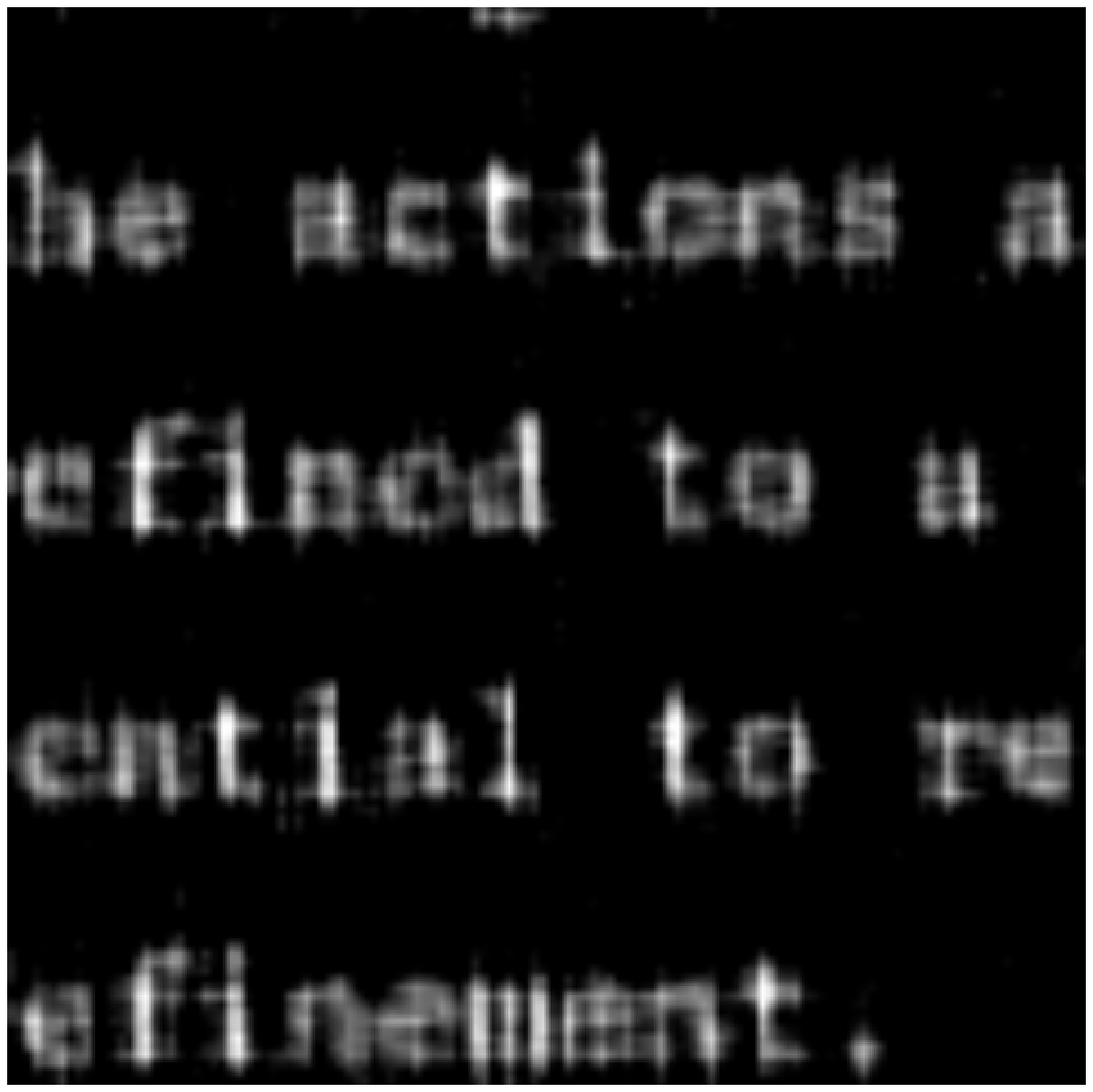}\\
\small{(a)} & \small{(b)} & \small{(c)}\\
\includegraphics[width=4.7cm]{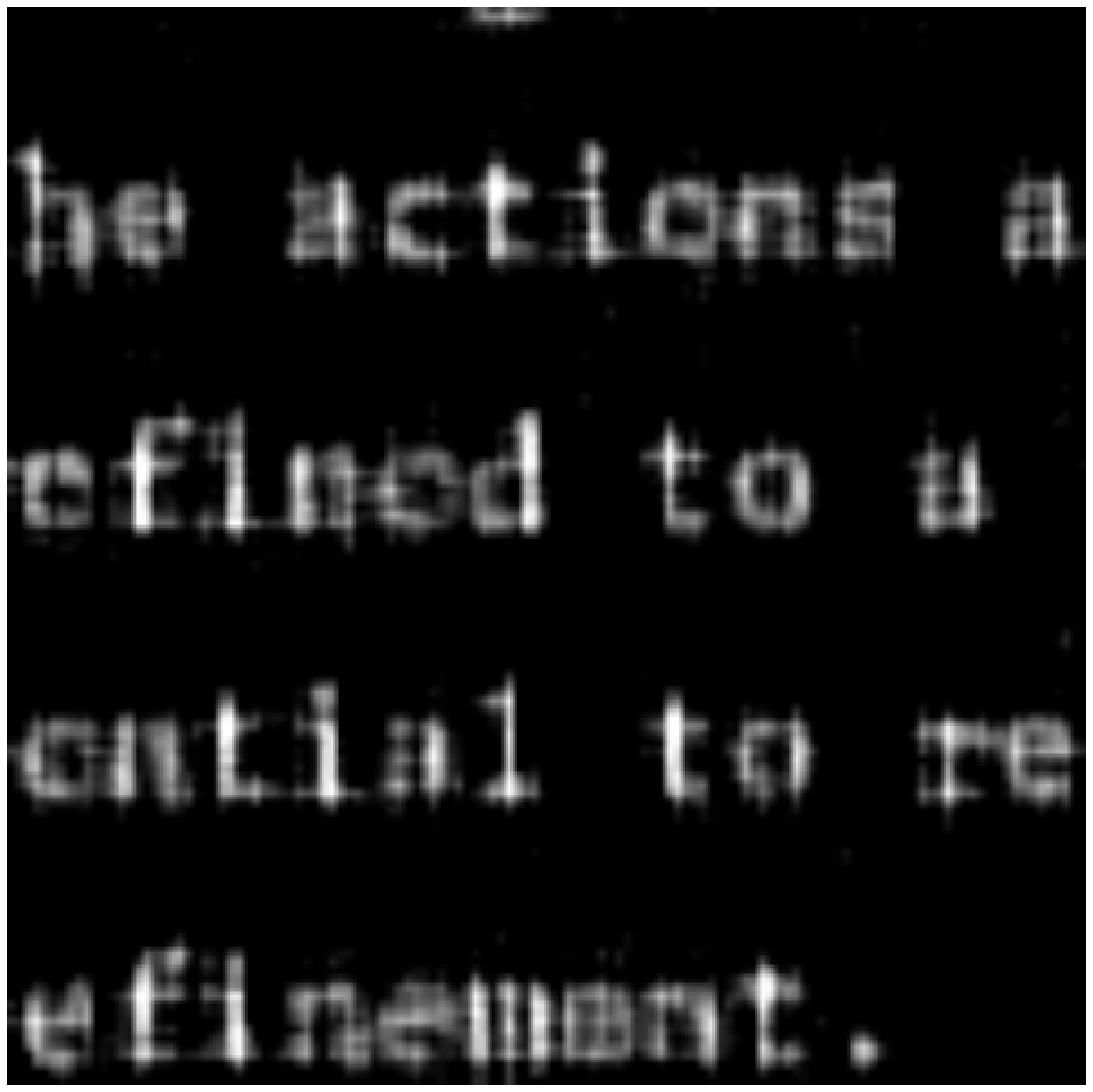}&
\hspace{0.4cm}
\includegraphics[width=4.7cm]{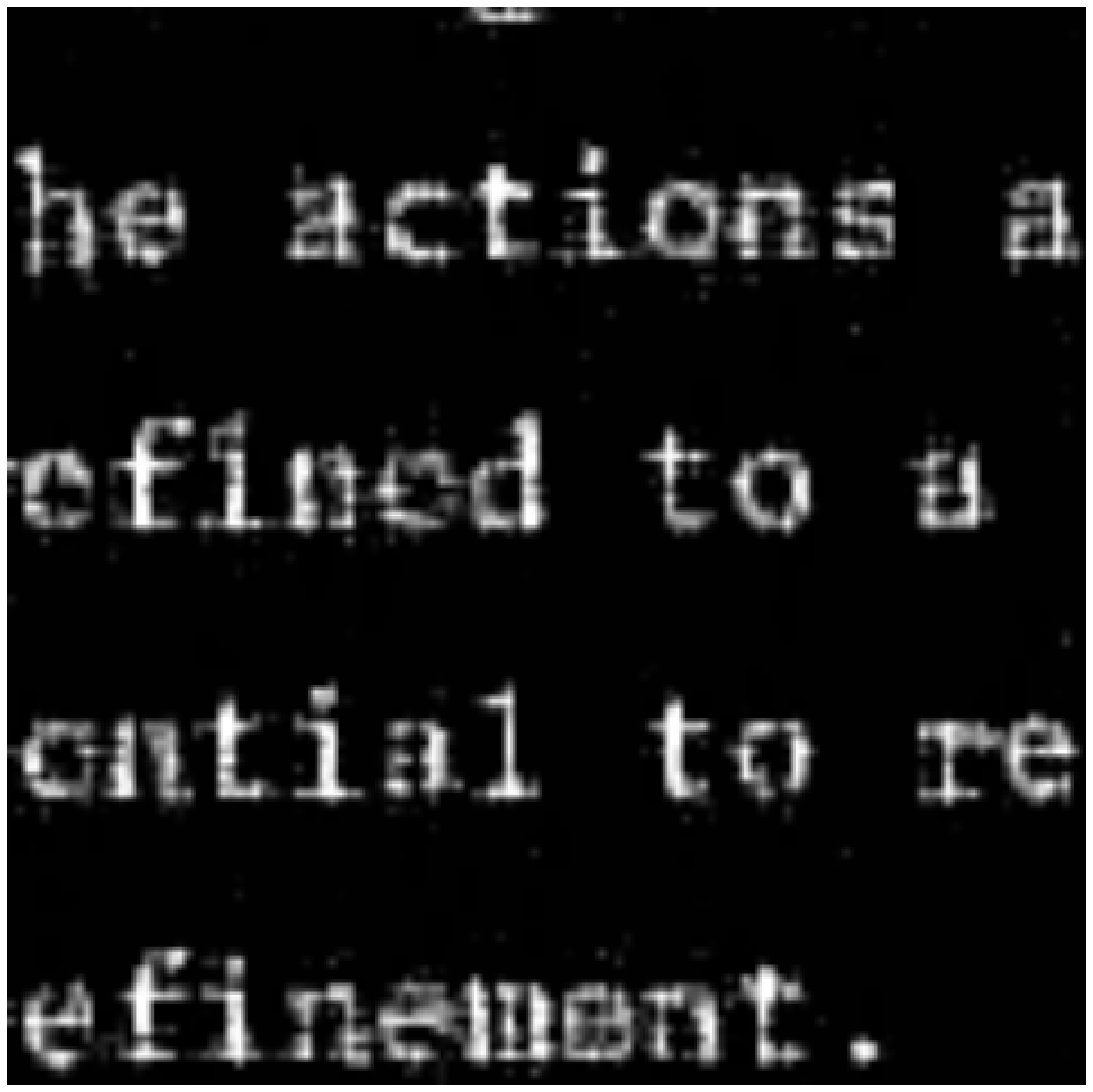}&
\hspace{0.4cm}
\includegraphics[width=4.7cm]{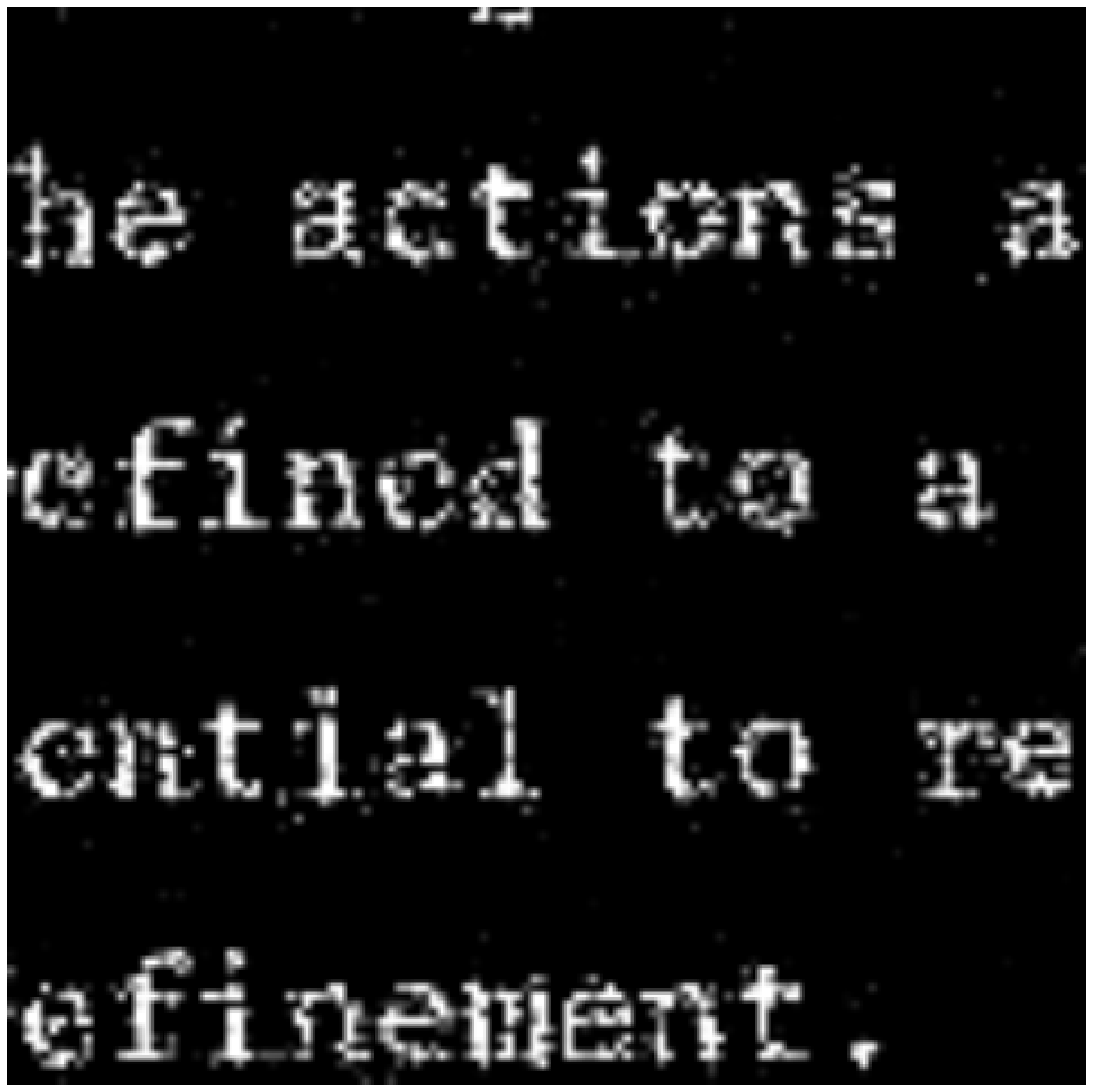}\\
\small{(d)} & \small{(e)} & \small{(f)}\\
\end{tabular} 
\caption{(a) Original image $\overline{x}$. 
(b) Degraded image $y$.
(c) Image restored by the forward-backward algorithm
(Algorithm~\ref{al:fb}) after 50 iterations.
(d) Image restored by the inertial forward-backward algorithm
(Algorithm~\ref{al:cd}) after 50 iterations.
(e) Image restored by the Douglas-Rachford algorithm 
(Algorithm~\ref{al:dr}) after 50 iterations.
(f) Restored image (all algorithms yield visually equivalent
images).
}
\label{fig:ex1im}
\end{figure}

\begin{figure}[ht!]
\begin{tikzpicture}[scale=0.55]
\begin{axis}[height=10cm,width=14cm,xmin =0, xmax=40, ymax=0,
legend cell align={left}, legend style={at={(0.02,0.02)},
anchor=south west}]
\addplot [dashed, very thick, mark=none, color=dgreen] 
table[x={time_FB}, 
y={FB}] {figures/ex1/dist_time_log.txt};
\addlegendentry{Alg.~\ref{al:fb}-Forward-backward}
\addplot [dashed, very thick, mark=none, color=blue] 
table[x={time_CD}, 
y={CD}] {figures/ex1/dist_time_log.txt};
\addlegendentry{Alg.~\ref{al:cd}-Inertial forward-backward}
\addplot [very thick, mark=none, color=red] table[x={time_DR}, 
y={DR}] {figures/ex1/dist_time_log.txt};
\addlegendentry{Alg.~\ref{al:dr}-Douglas-Rachford}
\end{axis}
\end{tikzpicture}
\hfill
\begin{tikzpicture}[scale=0.55]
\begin{axis}[height=10cm,width=14cm,xmin =0, xmax=500, ymax=0,
legend cell align={left}, legend style={at={(0.02,0.02)},
anchor=south west}]
\addplot [dashed, very thick, mark=none, color=dgreen] table[x={n}, 
y={FB}] {figures/ex1/dist_log.txt};
\addlegendentry{Alg.~\ref{al:fb}-Forward-backward}
\addplot [dashed, very thick, mark=none, color=blue] table[x={n}, 
y={CD}] {figures/ex1/dist_log.txt};
\addlegendentry{Alg.~\ref{al:cd}-Inertial forward-backward}
\addplot [very thick, mark=none, color=red] table[x={n}, 
y={DR}] {figures/ex1/dist_log.txt};
\addlegendentry{Alg.~\ref{al:dr}-Douglas-Rachford}
\end{axis}
\end{tikzpicture}\\
\begin{tikzpicture}[scale=0.55]
\hskip -0.0mm
\begin{axis}[height=10cm,width=14cm, legend cell align={left},
xmin =0, xmax=6.5, ymax=0]
\addplot [dashed, very thick, mark=none, color=dgreen] 
table[x={time_FB},
y={FB}] {figures/ex1/cost_time_log.txt};
\addlegendentry{Alg.~\ref{al:fb}-Forward-backward}
\addplot [dashed, very thick, mark=none, color=blue] 
table[x={time_CD},
y={CD}] {figures/ex1/cost_time_log.txt};
\addlegendentry{Alg.~\ref{al:cd}-Inertial forward-backward }
\addplot [very thick, mark=none, color=red] table[x={time_DR},
y={DR}] {figures/ex1/cost_time_log.txt};
\addlegendentry{Alg.~\ref{al:dr}-Douglas-Rachford}
\end{axis}
\end{tikzpicture}
\hfill
\begin{tikzpicture}[scale=0.55]
\begin{axis}[height=10cm,width=14cm, legend cell align={left},
xmin =0, xmax=100, ymax=0]
\addplot [dashed, very thick, mark=none, color=dgreen] table[x={n},
y={FB}] {figures/ex1/cost_log.txt};
\addlegendentry{Alg.~\ref{al:fb}-Forward-backward}
\addplot [dashed, very thick, mark=none, color=blue] table[x={n},
y={CD}] {figures/ex1/cost_log.txt};
\addlegendentry{Alg.~\ref{al:cd}-Inertial forward-backward }
\addplot [very thick, mark=none, color=red] table[x={n},
y={DR}] {figures/ex1/cost_log.txt};
\addlegendentry{Alg.~\ref{al:dr}-Douglas-Rachford}
\end{axis}
\end{tikzpicture}
\caption{Top left: Normalized distance in dB to the asymptotic image
produced by each algorithm versus execution time in seconds. 
Top right: Normalized distance in dB to the asymptotic image
produced by each algorithm versus iteration number.
In this experiment the objective function values
remain finite and can therefore be displayed.
Bottom left: Normalized objective function of \eqref{e:p1} in dB
versus execution time in seconds.
Bottom right: Normalized objective function of \eqref{e:p1} in dB
versus iteration number.}
\label{fig:ex1cost}
\end{figure}
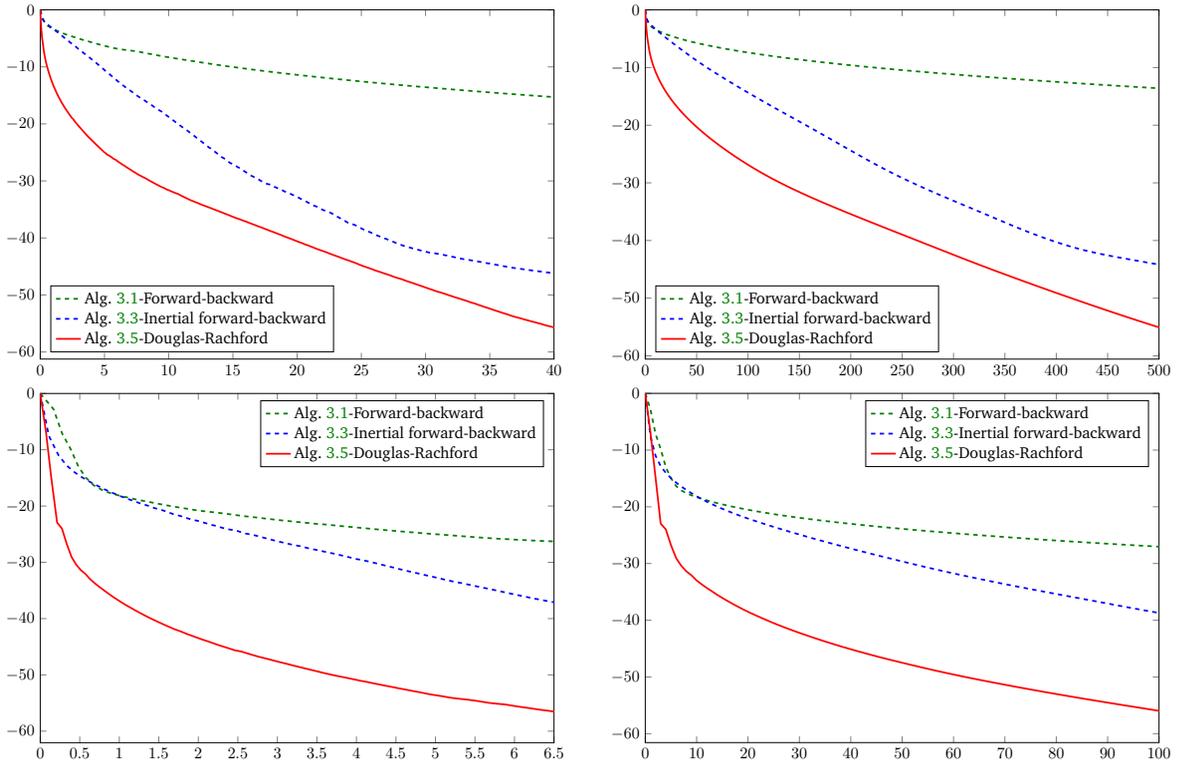 
We compare the splitting methods reviewed in Section~\ref{sec:3} 
on various digital image restoration and reconstruction problems. 
All images have $\sqrt{N}\times\sqrt{N}$ pixels and therefore 
the underlying Hilbert space is $\HH=\RR^N$ 
($N\in\{96^2,128^2,512^2\}$) equipped
with the standard Euclidean norm $\|\cdot\|_2$. 
All the algorithms guarantee the
convergence of their iterates $(x_n)_{n\in\NN}$ to a solution $x$
of the underlying optimization problem. Let us also note that this
set of experiments constitutes
the first implementation of Algorithm~\ref{al:kt} to image recovery. 
The simulations are run in Python on a personal computer running 
Linux Ubuntu version 18.04 with a 2.60GHz dual-core processor 
and 8Gb of RAM. Finally, the normalization used to plot the
decibel value of the squared distance of an iterate $x_n$ 
to a solution $x$ is
\begin{equation}
20\log_{10}\dfrac{\|x_n-x\|_2}{\|x_0-x\|_2},
\end{equation}
and that used for the objective value $\varphi(x_n)$ at iteration
$n$ is 
\begin{equation}
10\log_{10}\dfrac{\varphi(x_n)-\varphi(x)}{\varphi(x_0)-\varphi(x)}.
\end{equation}

\begin{remark}
\label{r:2012}
\rm
The applications to be considered below are instances of
Problem~\ref{prob:1} in which $f$ has bounded domain and the
remaining functions have full domain, which ensures existence of at
least one solution \cite[Corollary~11.16(i)]{Livre1}. In addition, 
the inclusion \eqref{e:cq1} holds by virtue of 
\cite[Proposition~4.3(ii)]{Svva12}.
\end{remark}

\subsection{Sparse image deconvolution}
\label{sec:41}
We consider a very basic instance of
Problem~\ref{prob:1} with only two functions. More specifically, 
we compare the numerical behavior of the forward-backward
algorithms of Propositions~\ref{p:fb} and \ref{p:cd} with that 
of the Douglas-Rachford algorithm of Proposition~\ref{p:dr},
which is a fully proximal method. Note that an elementary
comparison of these three algorithms to minimize $f+h$ already 
appears in Fig.~\ref{fig:Jim}, where $f=0$ and $h=\varphi$.
The images have size $128\times 128$.

The original image is $\overline{x}$ and the degraded image is 
$y=H\overline{x}+w$,
where $H$ models a convolution with a uniform rectangular kernel 
of size $15\times 5$ and $w$ is a Gaussian white noise
realization (see Fig.~\ref{fig:ex1im}(a)--(b)). 
The blurred image-to-noise-ratio is 15.5 dB.
Since each pixel value is known to be in $[0,255]$, we use
the hard constraint set $C=\left[0,255\right]^N$. As is
customary, the natural sparsity of $\overline{x}$ is promoted
using the function $\|\cdot\|_1$. Altogether, the problem is to 
\begin{equation}
\label{e:p1}
\minimize{x\in C}\|x\|_1+\frac{1}{2}\|Hx-y\|_2^2.
\end{equation}
Now set $f=\|\cdot\|_1+\iota_{C}$, set $h=\|H\cdot-y\|_2^2/2$, and
let $\gamma\in\RPP$.
Then $f\in\Gamma_0(\HH)$ and
$\prox_{\gamma f}=\proj_C\circ\soft{\gamma}$
\cite[Propositions~24.12(ii) and 24.47]{Livre1}, where 
$\soft{\gamma}$ is the soft thresholder on $[-\gamma,\gamma]$; 
the projector $\proj_C$ is implemented by
setting to $0$ the pixel values less than $0$, and to 
$255$ those larger than $255$.
On the other hand, $h$ is smooth, with gradient and proximity 
operator provided by Example~\ref{ex:-1} as
\begin{equation}
\label{e:ex1}
\nabla h\colon x\mapsto H^*(Hx-y)
\quad\text{and}\quad
\prox_{\gamma h}\colon x\mapsto
(\Id+\gamma H^*H)^{-1}(x+\gamma H^*y).
\end{equation}
The linear operator $H$ models a convolution and it is
representable by a block-circulant matrix. The computation of the
inverse in \eqref{e:ex1} is therefore straightforward via the 
fast Fourier transform \cite{Andr77}. 

We solve \eqref{e:p1} with the forward-backward algorithms
(Algorithm~\ref{al:fb} and Algorithm~\ref{al:cd}), as well as with
the fully proximal Douglas-Rachford algorithm
(Algorithm~\ref{al:dr}).
The algorithms are initialized at zero and implemented with
the parameters for which they seem to perform best, that is,
$\gamma_n\equiv 1.99/\beta$ for Algorithm~\ref{al:fb},
$\gamma=1/\beta$ and $\alpha=3$ for Algorithm~\ref{al:cd}, and 
$\gamma=30$ and $\lambda_n\equiv 1.9$ for Algorithm~\ref{al:dr}. 
The results of
Figs.~\ref{fig:ex1im}--\ref{fig:ex1cost} show a superior
performance for Algorithm~\ref{al:dr}, which is
fully proximal.

\subsection{Multiview image reconstruction from partial diffraction
data}
\label{sec:42}

\begin{figure}[t]
\centering
\begin{tabular}{c@{}c@{}}
\includegraphics[width=4.7cm]{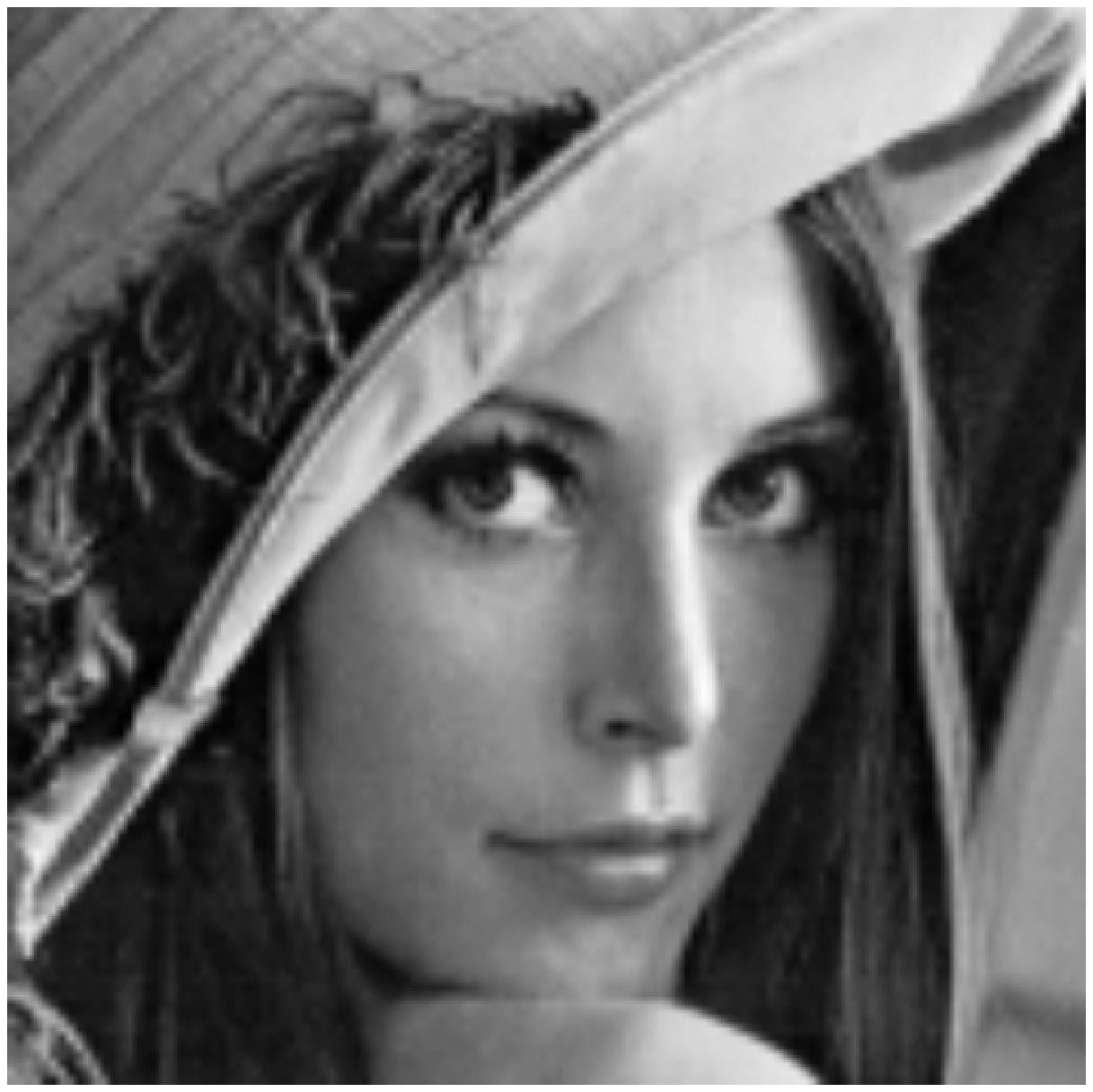}&
\hspace{0.4cm}
\includegraphics[width=4.7cm]{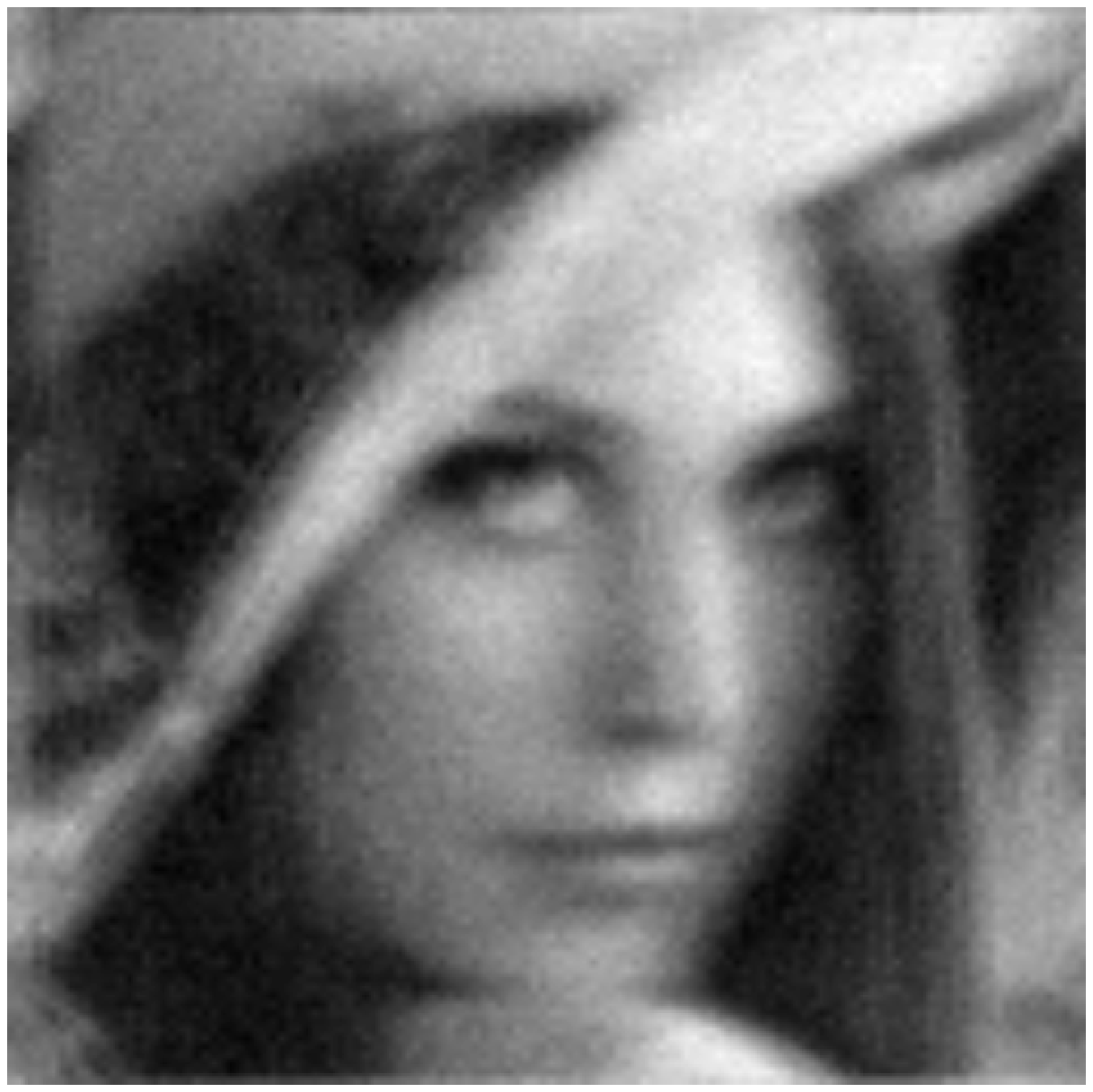}\\
\small{(a)} & \small{(b)}\\
\includegraphics[width=4.7cm]{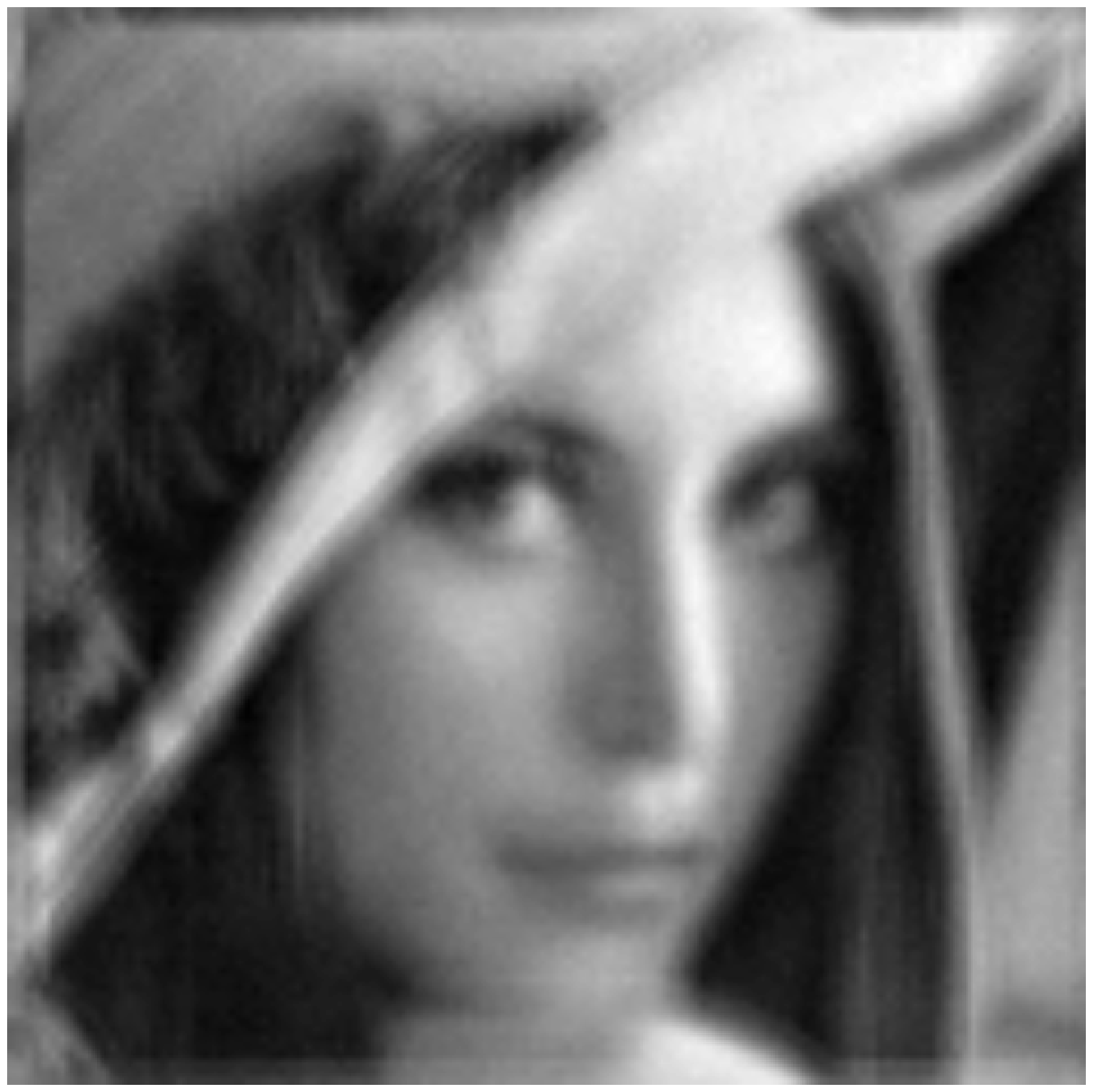}&
\hspace{0.4cm}
\includegraphics[width=4.7cm]{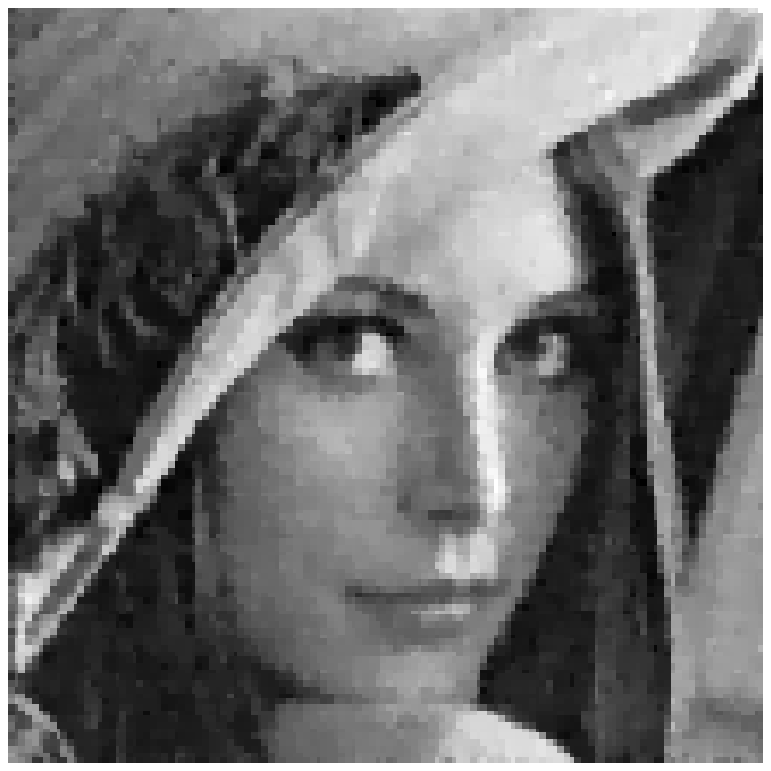}\\
\small{(c)} & \small{(d)}\\
\end{tabular} 
\caption{(a) Original image $\overline{x}$. 
(b) Degraded image $y_1$.
(c) Degraded image $y_2$. 
(d) Reconstructed image (all algorithms yield visually equivalent
images).}
\label{fig:ex2im}
\end{figure}

\begin{figure}[t]
\centering
\begin{tabular}{c@{}c@{}c@{}}
\includegraphics[width=4.7cm]{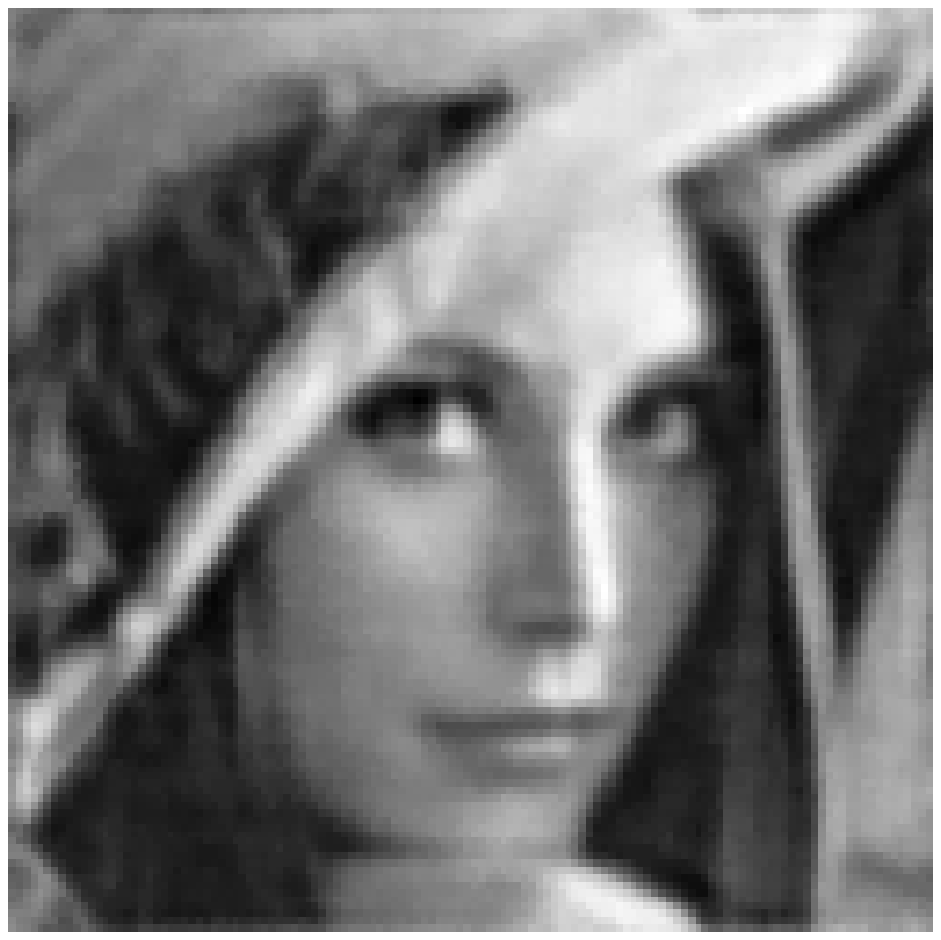}&
\hspace{0.4cm}
\includegraphics[width=4.7cm]{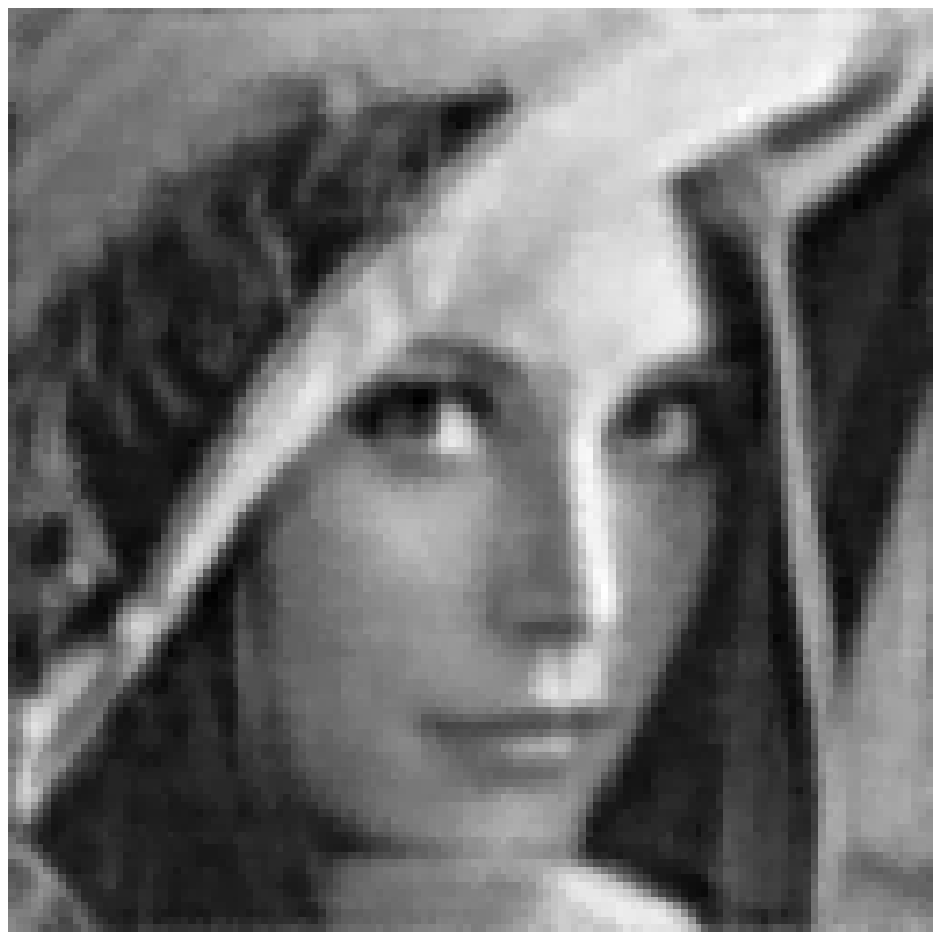}&
\hspace{0.4cm}
\includegraphics[width=4.7cm]{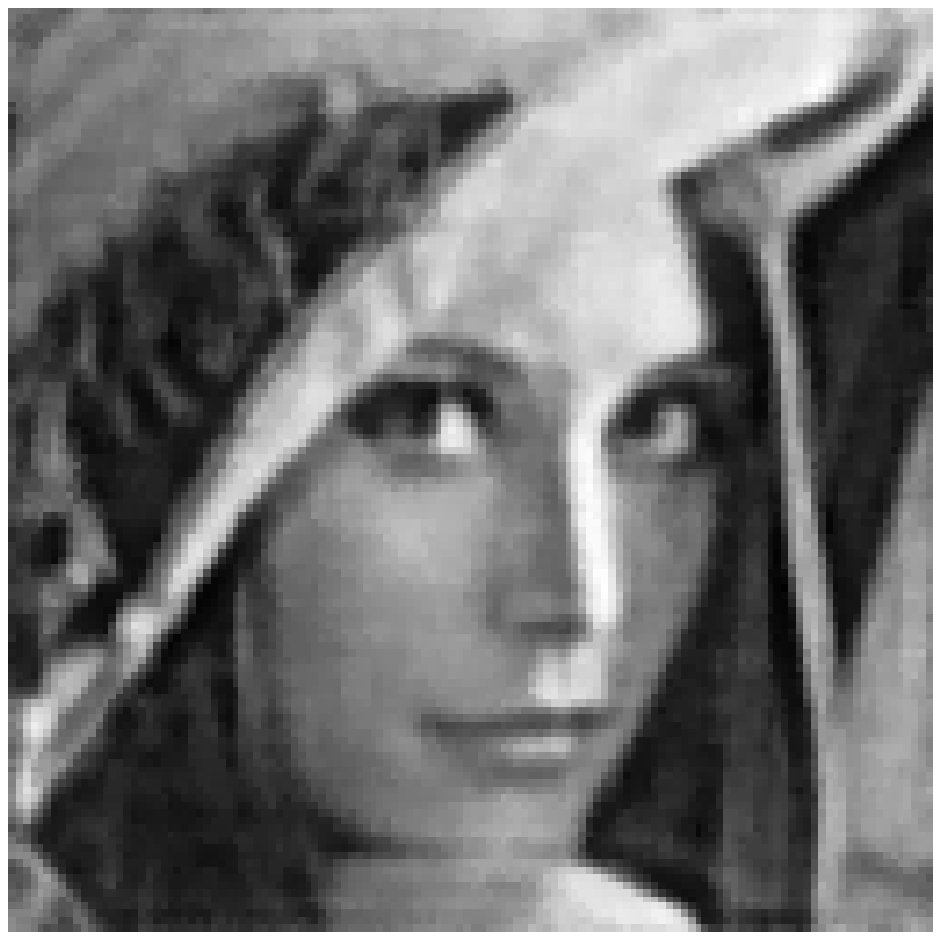}\\
\small{(a)} & \small{(b)} & \small{(c)}\\
\includegraphics[width=4.7cm]{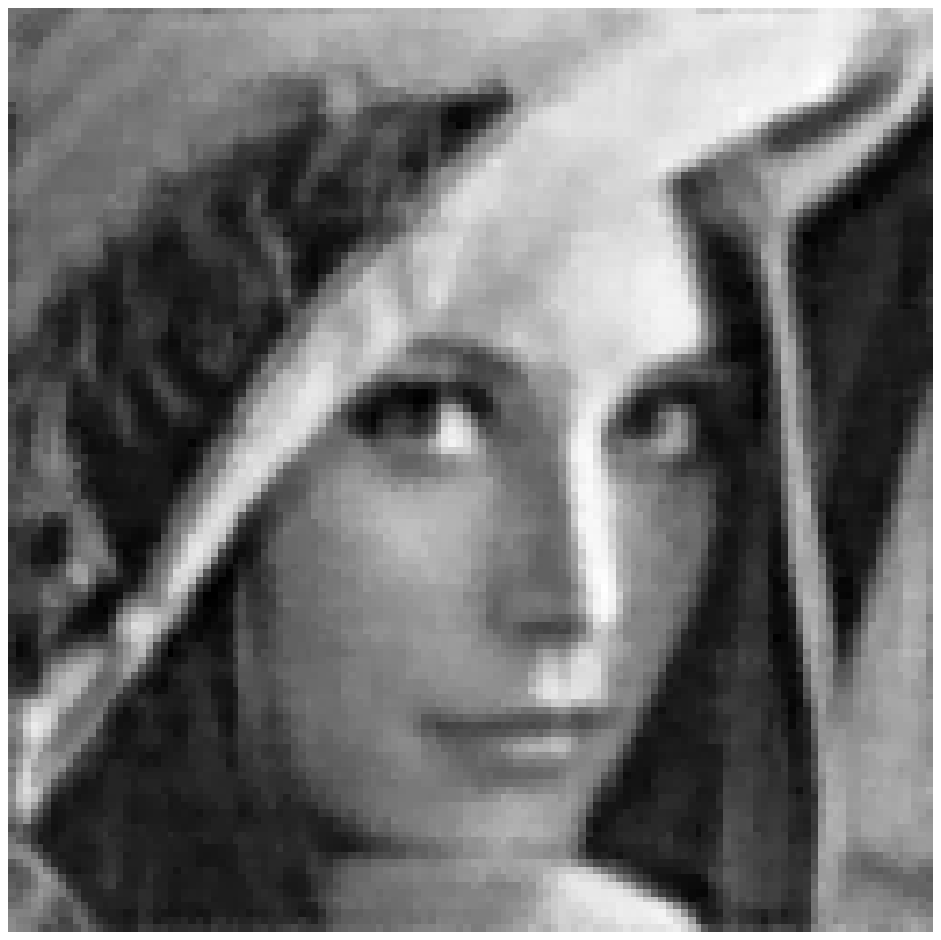}&
\hspace{0.4cm}
\includegraphics[width=4.7cm]{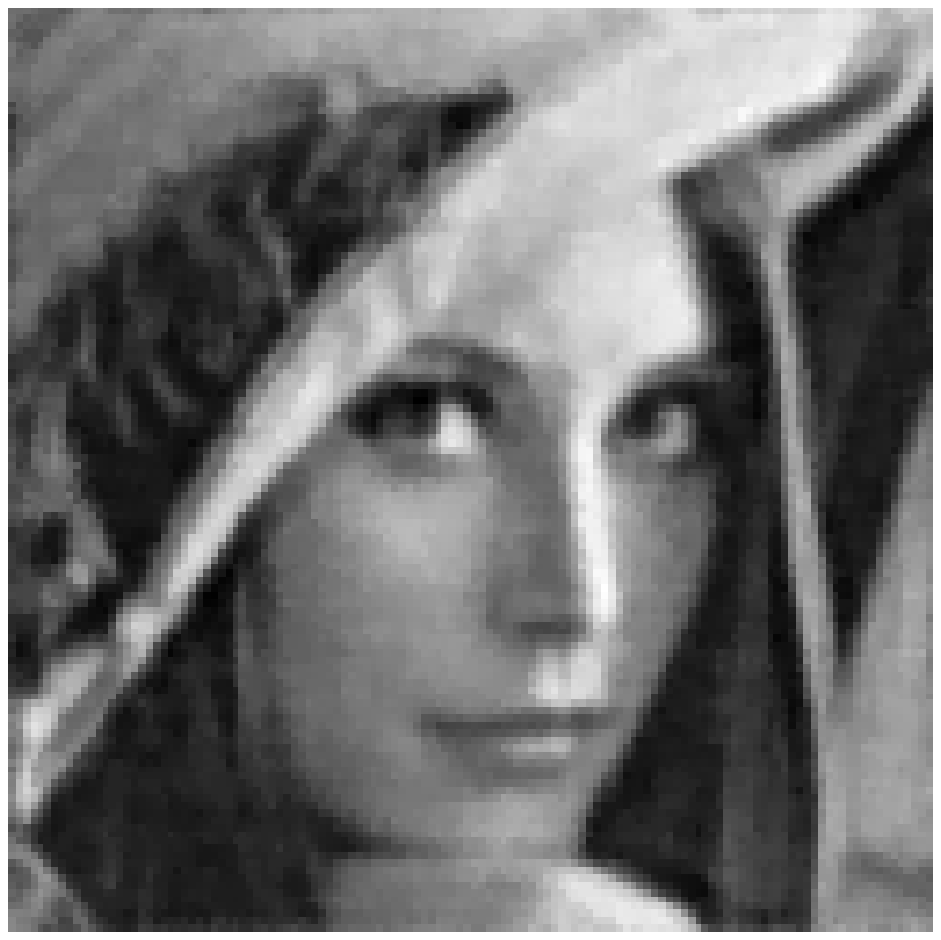}&
\hspace{0.4cm}
\includegraphics[width=4.7cm]{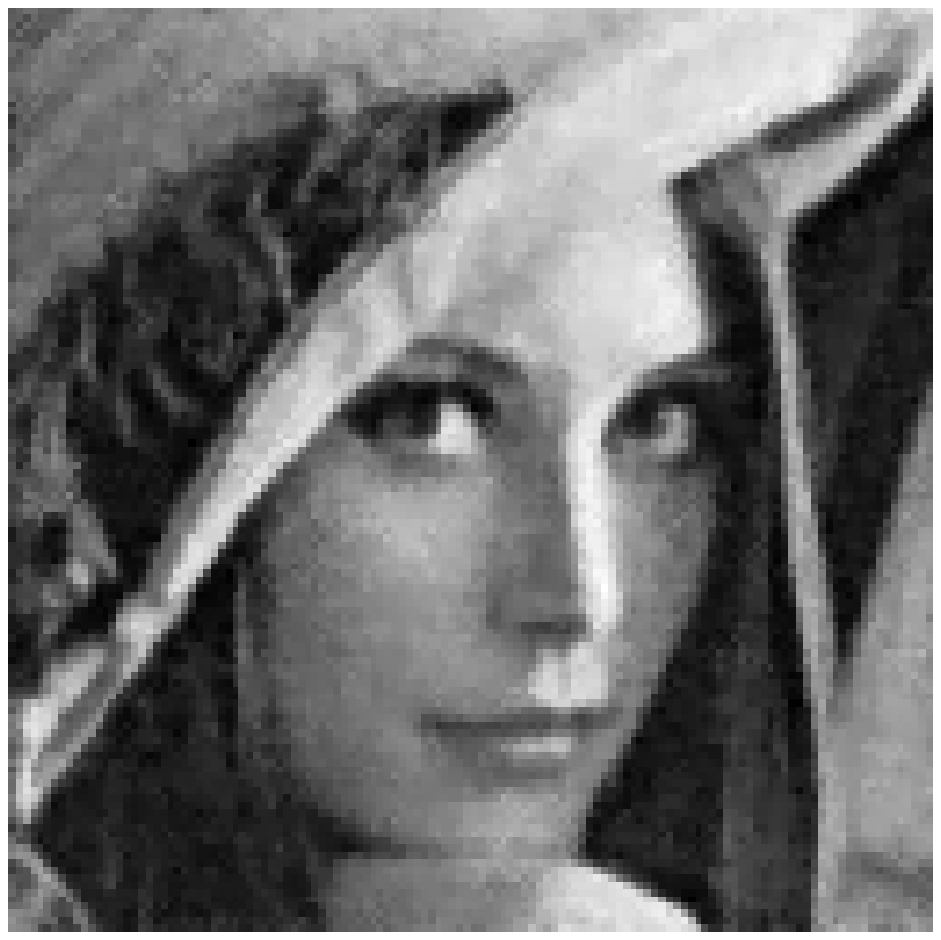}\\
\small{(d)} & \small{(e)} & \small{(f)}\\
\end{tabular} 
\caption{ 
(a) Image restored by Pb.~\ref{prob:50}/Alg.~\ref{al:2012}-S 
after 20 iterations.
(b) Image restored by Pb.~\ref{prob:50}/Alg.~\ref{al:fbp}-S 
after 20 iterations.
(c) Image restored by Pb.~\ref{prob:50}/Alg.~\ref{al:kt}-S 
after 20 iterations.
(d) Image restored by Pb.~\ref{prob:51}/Alg.~\ref{al:2012}-P 
after 20 iterations.
(e) Image restored by Pb.~\ref{prob:51}/Alg.~\ref{al:fbp}-P 
after 20 iterations.
(f) Image restored by Pb.~\ref{prob:51}/Alg.~\ref{al:kt}-P 
after 20 iterations.}
\label{fig:ex2im2}
\end{figure}
We consider the problem of reconstructing a $128\times 128$
image $\overline{x}$ from a partial observation of its diffraction 
over some frequency range $R$, possibly with 
measurement errors \cite{Seza83}. To exploit this information we 
use the soft constraint penalty $d_E$ associated with the set
\begin{equation}
\label{e:E}
E=\menge{x\in\RR^N}{(\forall k\in R)\;
\widehat{x}(k)=\widehat{\overline{x}}(k)},
\end{equation}
where $\widehat{x}$ denotes the two-dimensional discrete Fourier 
transform of $x$. The set $R$ contains the frequencies in 
$\{0,\ldots,15\}^2$ as well as those resulting from the 
symmetry properties of the discrete Fourier transform. In
addition, we have at our disposal two blurred noisy observations of 
$\overline{x}$, namely (see Fig.~\ref{fig:ex2im}(a)--(c))
$y_1=H_1\overline{x}+w_1$ and $y_2=H_2\overline{x}+w_2$. 
Here, $H_1$ and $H_2$ model convolutional blurs with kernels 
of size $3\times 11$ and of $7\times 5$, respectively, and $w_1$
and $w_2$ are Gaussian white noise realizations. The blurred
image-to-noise-ratios are 27.3 dB and 35.4 dB.
\begin{figure}[ht!]
\begin{tikzpicture}[scale=0.55]
\hskip -2mm
\begin{axis}[height=10cm,width=14cm,legend style={at={(0.02,0.02)},
legend cell align={left},
anchor=south west}, xmin =0, xmax=35, ymax=0]
\addplot [dashed, very thick, mark=none, color=blue] 
table[x={time_FBF_Expl}, y={FBF_Expl}]
{figures/ex2/dist_time_log.txt};
\addlegendentry{Pb.~\ref{prob:50}/Alg.~\ref{al:2012}-S}
\addplot [very thick, mark=none, color=blue] table[x={time_FBF_Imp}, 
y={FBF_Imp}]{figures/ex2/dist_time_log.txt};
\addlegendentry{Pb.~\ref{prob:51}/Alg.~\ref{al:2012}-P}
\addplot [dashed, very thick, mark=none, color=dgreen] 
table[x={time_FB_Expl}, y={FB_Expl}] {figures/ex2/dist_time_log.txt};
\addlegendentry{Pb.~\ref{prob:50}/Alg.~\ref{al:fbp}-S}
\addplot [very thick, mark=none, color=dgreen] 
table[x={time_FB_Imp}, 
y={FB_Imp}]{figures/ex2/dist_time_log.txt};
\addlegendentry{Pb.~\ref{prob:51}/Alg.~\ref{al:fbp}-P}
\addplot [dashed, very thick, mark=none, color=red] 
table[x={time_KT_Expl}, y={KT_Expl}]{figures/ex2/dist_time_log.txt};
\addlegendentry{Pb.~\ref{prob:50}/Alg.~\ref{al:kt}-S}
\addplot [very thick, mark=none, color=red] table[x={time_KT_Imp}, 
y={KT_Imp}] {figures/ex2/dist_time_log.txt};
\addlegendentry{Pb.~\ref{prob:51}/Alg.~\ref{al:kt}-P}
\end{axis}
\hfill
\end{tikzpicture}
\begin{tikzpicture}[scale=0.55]
\begin{axis}[height=10cm,width=14cm,legend style={at={(0.02,0.02)},
legend cell align={left},
anchor=south west}, xmin =0, xmax=160, ymax=0]
\addplot [dashed, very thick, mark=none, color=blue] 
table[x={n}, y={FBF_Expl}]{figures/ex2/dist_log.txt};
\addlegendentry{Pb.~\ref{prob:50}/Alg.~\ref{al:2012}-S}
\addplot [very thick, mark=none, color=blue] table[x={n}, 
y={FBF_Imp}]{figures/ex2/dist_log.txt};
\addlegendentry{Pb.~\ref{prob:51}/Alg.~\ref{al:2012}-P}
\addplot [dashed, very thick, mark=none, color=dgreen] 
table[x={n}, y={FB_Expl}] {figures/ex2/dist_log.txt};
\addlegendentry{Pb.~\ref{prob:50}/Alg.~\ref{al:fbp}-S}
\addplot [very thick, mark=none, color=dgreen] table[x={n}, 
y={FB_Imp}]{figures/ex2/dist_log.txt};
\addlegendentry{Pb.~\ref{prob:51}/Alg.~\ref{al:fbp}-P}
\addplot [dashed, very thick, mark=none, color=red] 
table[x={n}, y={KT_Expl}]{figures/ex2/dist_log.txt};
\addlegendentry{Pb.~\ref{prob:50}/Alg.~\ref{al:kt}-S}
\addplot [very thick, mark=none, color=red] table[x={n}, y={KT_Imp}] 
{figures/ex2/dist_log.txt};
\addlegendentry{Pb.~\ref{prob:51}/Alg.~\ref{al:kt}-P}
\end{axis}
\hfill
\end{tikzpicture}
\centering
\caption{Left: Normalized distance in dB to the asymptotic image
produced by each algorithm versus execution time in seconds. 
Right: Normalized distance in dB to the asymptotic image
produced by each algorithm versus iteration number.
Each algorithm is represented by a given color, the solid line 
corresponds to the fully proximal implementation and the dashed 
line to the implementation with gradient steps for the smooth
functions.}
\label{fig:ex2costdist}
\end{figure}
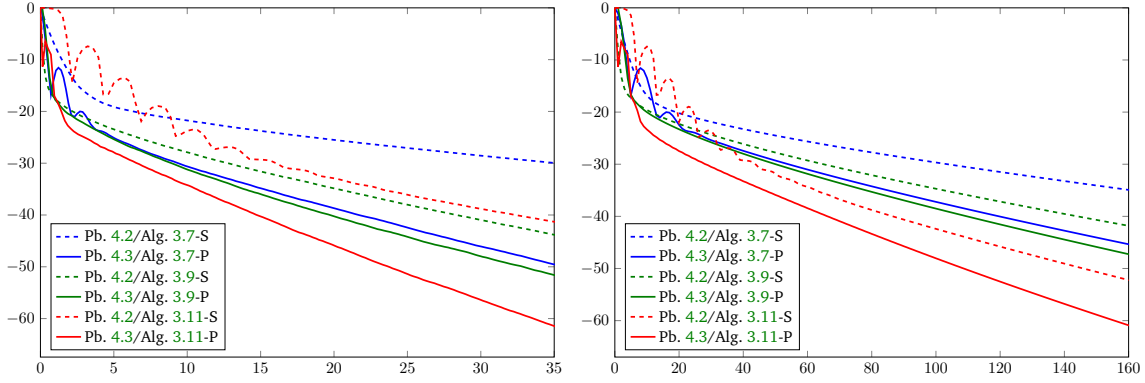
We use $C=\left[0,255\right]^N$ as a hard constraint
set. Finally, we utilize a discrete version of the Gauss-TV 
penalty of Remark~\ref{r:gaussTV} to control
oscillations in the reconstructed image. This leads to the
formulation
\begin{equation}
\label{e:p2}
\minimize{x\in C}{\frac{1}{2}d_E(x)+\frac{2}{5}h(Dx)}
+\frac{3}{4}\|H_1x-y_1\|_2^2+\frac{3}{4}\|H_2x-y_2\|_2^2,
\end{equation}
where $D\colon\RR^N\to\RR^N\times\RR^N\colon x
\mapsto(G_1x,G_2x)$, $G_1$ and $G_2$ being
horizontal and vertical discrete difference operators, 
and where 
$(\forall(y_1,y_2)=((\eta_{1,k})_{1\leq k\leq N},
(\eta_{2,k})_{1\leq k\leq N})\in\GG_2=\RR^N\times\RR^N)$
$h(y_1,y_2)=\sum_{k=1}^N
\mathfrak{h}_2(\|(\eta_{1,k},\eta_{2,k})\|_2)$, $\mathfrak{h}_2$ 
being the Huber function \eqref{e:huber}.
We derive from \eqref{e:p2} two versions of
Problem~\ref{prob:1}.

\begin{problem}
\label{prob:50}
\rm
In Problem~\ref{prob:1}, set 
$f=\iota_C$,
$I=\{1\}$, 
$g_1=0.5d_E$, 
$L_1=\Id$, $J=\{2,3,4\}$,
$h_2=0.4h$, $L_2=D$,
$h_3=0.75\|\cdot-y_1\|_2^2$, $L_3=H_1$,
$h_4=0.75\|\cdot-y_2\|_2^2$, and $L_4=H_2$.
\end{problem}

\begin{problem}[fully proximal]
\label{prob:51}
\rm
In Problem~\ref{prob:1}, set 
$f=\iota_C$,
$I=\{1,2,3,4\}$, $J=\emp$,
$g_1=0.5d_E$, $L_1=\Id$, $g_2=0.4h$, $L_2=D$, 
$g_3=0.75\|H_1\cdot-y_1\|_2^2$, $L_3=\Id$,
$g_4=0.75\|H_2\cdot-y_2\|_2^2$, and $L_4=\Id$.
\end{problem}

We apply to these two problems Algorithms~\ref{al:2012}, 
\ref{al:fbp}, and \ref{al:kt} with all initial
vectors set to $0$. The following parameters are used, where
$\beta=\sqrt{\sum_{i\in I}\|L_i\|^2}+
\sum_{j\in J}\mu_j\|L_j\|^2$ (these parameters were found to
optimize the performance of each algorithm):
\begin{itemize}
\setlength{\itemsep}{-1pt}
\item
Algorithm~\ref{al:2012}-S (with smooth terms for
Problem~\ref{prob:50}): $\gamma_n\equiv 0.99/\beta$.
\item
Algorithm~\ref{al:2012}-P (fully proximal for 
Problem~\ref{prob:51}): $\gamma_n\equiv 0.99/\beta$.
\item
Algorithm~\ref{al:fbp}-S (with smooth terms for
Problem~\ref{prob:50}):
$\sigma_{1,n}\equiv 8/(5\beta)$ and $\tau_n\equiv 8/(5\beta)$.
\item
Algorithm~\ref{al:fbp}-P (fully proximal for Problem~\ref{prob:51}):
$\sigma_{1,n}\equiv 1/(2\beta)$, 
$\sigma_{2,n}\equiv 1/(2\beta)$, 
$\sigma_{3,n}\equiv 3/\beta$,
$\sigma_{4,n}\equiv 3/\beta$, and 
$\tau_n\equiv 1/\beta$.
\item
Algorithm~\ref{al:kt}-S (with smooth terms for
Problem~\ref{prob:50}):
$\gamma_n\equiv 0.4$, 
$\mu_{1,n}\equiv 1.0$, 
$\mu_{2,n}\equiv 2.49$,
$\mu_{3,n}\equiv 0.65$,
$\mu_{4,n}\equiv 0.65$, and 
$\lambda_n\equiv 1.99$.
\item
Algorithm~\ref{al:kt}-P (fully proximal for Problem~\ref{prob:51}):
$\gamma_n\equiv 0.25$, 
$\mu_{1,n}\equiv 1.0$, 
$\mu_{2,n}\equiv 1.5$,
$\mu_{3,n}\equiv 1.0$,
$\mu_{4,n}\equiv 1.0$, and
$\lambda_n\equiv 1.99$.
\end{itemize}
The proximity operators and the gradients used in these
experiments follow from \cite[Example~24.28]{Livre1} for $g_1$, 
Example~\ref{ex:npsi2} (with $M=2$ and $\bK=\{1,\ldots,N\}$, 
and using \cite[Example~24.9]{Livre1} to get the proximity operator
of $\mathfrak{h}_{2}$) for $h_2$, 
and \eqref{e:ex1} for $g_3$ and $g_4$.
On the other hand $\prox_{\gamma f}=\proj_C$.
We have $\|D\|^2=8$ and $\|H_1\|^2=\|H_2\|^2=1$. On the
other hand, the functions $h_2$, $h_3$, and $h_4$ are 
differentiable with a Lipschitzian gradient, and their Lipschitz
constants are respectively $0.4$, $0.75$, and $0.75$. Thus, all the
assumptions required by the algorithms are satisfied. The results
shown in Figs.~\ref{fig:ex2im2} and \ref{fig:ex2costdist}
illustrate the faster convergence of the iterates $(x_n)_{n\in\NN}$
of fully proximal
algorithms to a solution $x$ compared to their gradient-based 
versions, both in terms of computation time and iterations. 
Note that these primal-dual algorithms do not guarantee 
Fej\'er monotonicity in the primal space, i.e., that 
$(\|x_n-x\|)_{n\in\NN}$ goes to $0$ monotonically. This is
confirmed by the oscillations seen in 
Fig.~\ref{fig:ex2costdist}.
Finally, Fig.~\ref{fig:ex2im}(d) shows that the formulation
\eqref{e:p2} provides a faithful recovery of $\overline{x}$.

\subsection{Image interpolation}

\begin{figure}[t]
\centering
\begin{tabular}{c@{}c@{}}
\includegraphics[width=4.7cm]{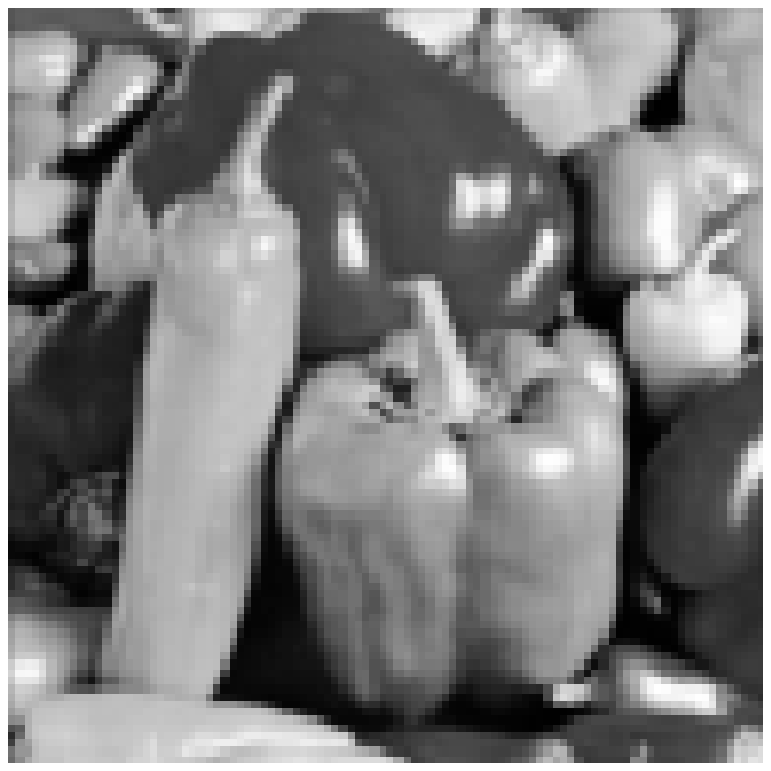}&
\hspace{0.4cm}
\includegraphics[width=4.7cm]{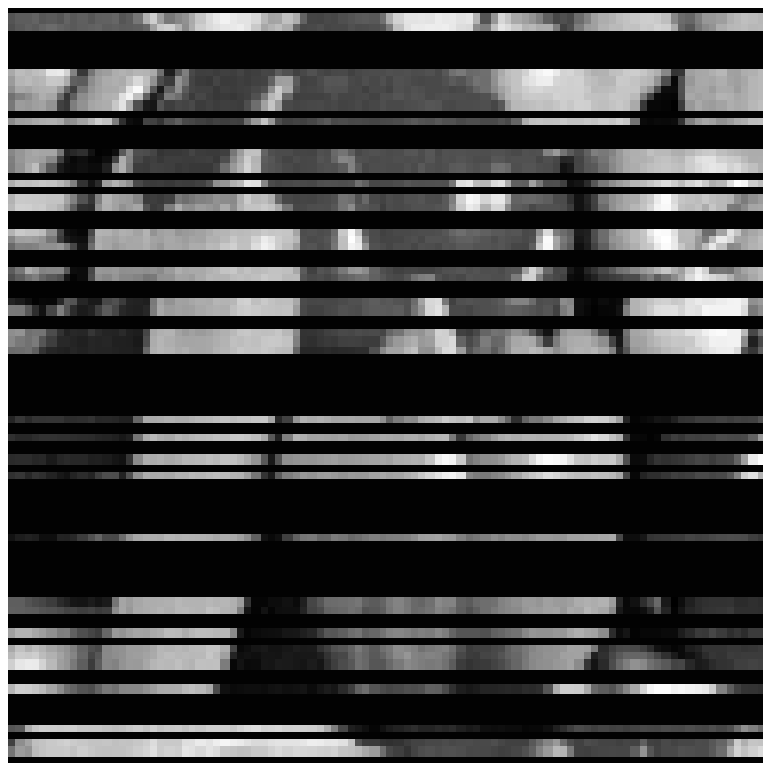}\\
\small{(a)} & \small{(b)}\\
\includegraphics[width=4.7cm]{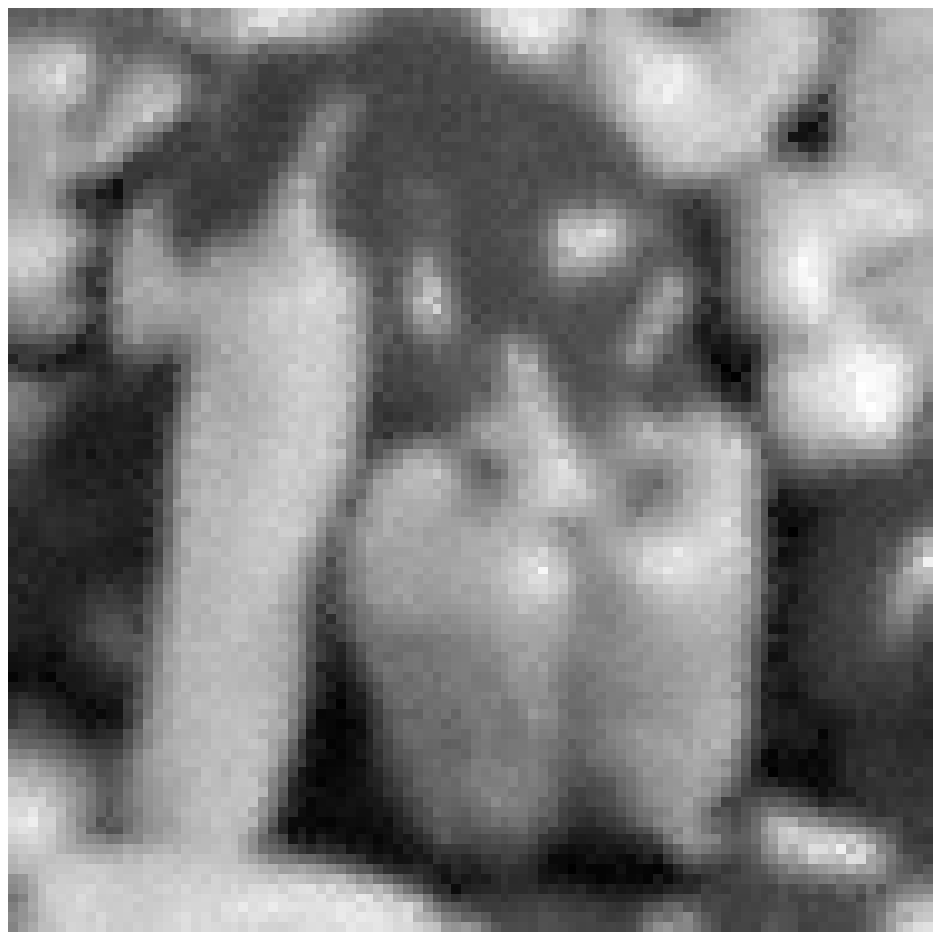}&
\hspace{0.4cm}
\includegraphics[width=4.7cm]{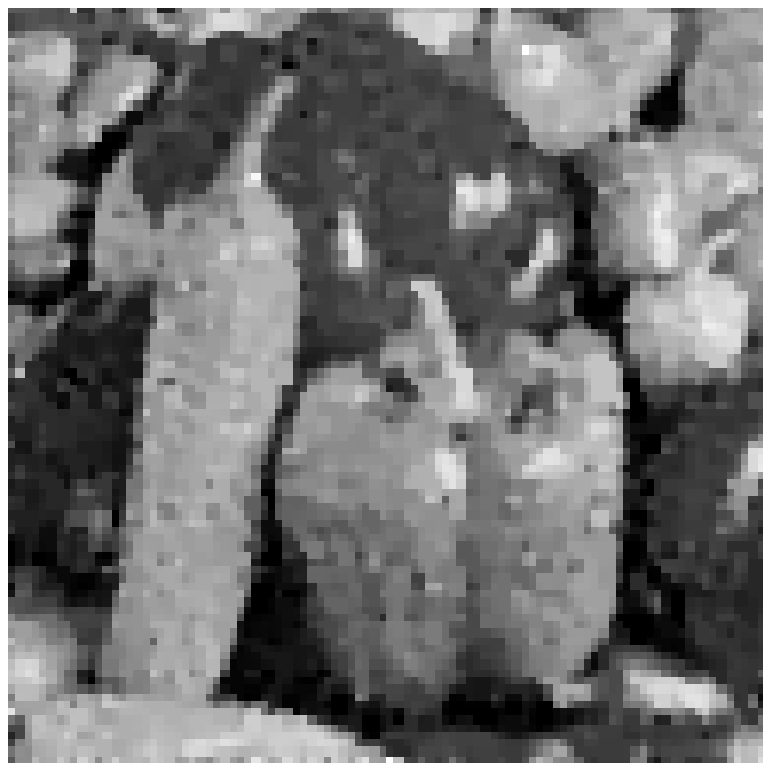}\\
\small{(c)} & \small{(d)}\\
\end{tabular} 
\caption{(a) Original image $\overline{x}$. 
(b) Degraded image $y_1$.
(c) Degraded image $y_2$. 
(d) Reconstructed image (all algorithms yield visually equivalent
images).}
\label{fig:ex3im}
\end{figure}

We consider the problem of reconstructing the 
$96\times 96$ original image $\overline{x}$ shown in 
Fig.~\ref{fig:ex3im}(a) from a noisy occulted version $y_1$ and a
blurred and noisy measurement $y_2$ 
(see Fig.~\ref{fig:ex3im}(b)--(c)). The occulted version is missing
57 lines and the observed line numbers are indexed by 
$R\subset\{1,\ldots,N\}$. The observations are 
$y_1=M\overline{x}+w_1$ and $y_2=H\overline{x}+w_2$,
where $M$ is the masking operator zeroing the lines not indexed by
$R$ and where $H$ is a blurring operator. By contrast with the 
previous experiments, the blurring is nonstationary and therefore
does not correspond to a convolution operation. More precisely, 
the action of the
blurring operator $H$ on a given pixel $(i,j)$ is to replace it by
an average of the neighboring pixels weighted by an isotropic
Gaussian kernel centered at $(i,j)$ and with random standard 
deviation $\sigma_{i,j}\in [0,1]$.
Finally $\|H\|=1$, while $w_1$ and $w_2$ are realizations of 
Gaussian white noises such that the image-to-noise-ratio for
$y_1$ is 25.9 dB and the blurred image-to-noise-ratio for $y_2$ is
31.0 dB. We denote by $(y^{(i)}_1)_{i\in R}$ the nonzero lines of
$y_1$ corresponding to the observed lines of $\overline{x}$.

To model this interpolation problem, we use 
$C=\left[0,255\right]^N$ as a hard constraint as well as 
the total variation penalty. In addition, we fit the observed lines 
via the penalty $x\mapsto\sum_{i\in R}10\|x^{(i)}-y^{(i)}_1\|_2$
and the degraded image via the penalty
$x\mapsto 5\|Hx-y_2\|_2^2$. This leads to the formulation
\begin{equation}
\label{e:p3}
\minimize{x\in C}{\|Dx\|_{1,2}+
10\sum_{i\in R}\|x^{(i)}-y^{(i)}_1\|_2}+5\|Hx-y_2\|_2^2,
\end{equation}
where $D$ is as in \eqref{e:p2} and
$(\forall(y_1,y_2)\in\GG_2=\RR^N\times\RR^N)$
$\|(y_1,y_2)\|_{1,2}=\sum_{k=1}^N\|(\eta_{1,k},\eta_{2,k})\|_2$. 
Two versions of Problem~\ref{prob:1} are employed.

\begin{problem}
\label{prob:52}
\rm
In Problem~\ref{prob:1}, set 
$f=\iota_C$, $I=\{1,2\}$, $J=\{3\}$, 
$g_1=\|\cdot\|_{1,2}$, $L_1=D$, 
$g_2=10\sum_{i\in R}\|x^{(i)}-y^{(i)}_1\|_2$, 
$L_2=\Id$,
$h_3=5\|\cdot-y_2\|_2^2$, and $L_3=H$.
\end{problem}

\begin{problem}[fully proximal]
\label{prob:53}
\rm
In Problem~\ref{prob:1}, set 
$f=\iota_C$, $I=\{1,2,3\}$, $J=\emp$,
$g_1=\|\cdot\|_{1,2}$, $L_1=D$, 
$g_2=10\sum_{i\in R}\|x^{(i)}-y^{(i)}_1\|_2$, 
$L_2=\Id$,
$g_3=5\|H\cdot-y_2\|_2^2$, and $L_3=\Id$.
\end{problem}

\begin{figure}[ht!]
\begin{tikzpicture}[scale=0.55]
\hskip -4mm
\begin{axis}[height=10cm,width=14cm,legend style={
legend cell align={left}}, xmin =0, xmax=80, ymax=0]
\addplot [dashed, very thick, mark=none, color=blue] 
table[x={time_FBF_Expl}, y={FBF_Expl}]
{figures/ex3/dist_time_log.txt};
\addlegendentry{Pb.~\ref{prob:52}/Alg.~\ref{al:2012}-S}
\addplot [very thick, mark=none, color=blue] table[x={time_FBF_Imp}, 
y={FBF_Imp}]{figures/ex3/dist_time_log.txt};
\addlegendentry{Pb.~\ref{prob:53}/Alg.~\ref{al:2012}-P}
\addplot [dashed, very thick, mark=none, color=dgreen] 
table[x={time_FB_Expl}, y={FB_Expl}] {figures/ex3/dist_time_log.txt};
\addlegendentry{Pb.~\ref{prob:52}/Alg.~\ref{al:fbp}-S}
\addplot [very thick, mark=none, color=dgreen] 
table[x={time_FB_Imp}, 
y={FB_Imp}]{figures/ex3/dist_time_log.txt};
\addlegendentry{Pb.~\ref{prob:53}/Alg.~\ref{al:fbp}-P}
\addplot [dashed, very thick, mark=none, color=red] 
table[x={time_KT_Expl}, y={KT_Expl}]{figures/ex3/dist_time_log.txt};
\addlegendentry{Pb.~\ref{prob:52}/Alg.~\ref{al:kt}-S}
\addplot [very thick, mark=none, color=red] table[x={time_KT_Imp}, 
y={KT_Imp}] 
{figures/ex3/dist_time_log.txt};
\addlegendentry{Pb.~\ref{prob:53}/Alg.~\ref{al:kt}-P}
\end{axis}
\end{tikzpicture}
\hfill
\begin{tikzpicture}[scale=0.55]
\begin{axis}[height=10cm,width=14cm,legend style={
legend cell align={left}}, xmin =0, xmax=900, ymax=0]
\addplot [dashed, very thick, mark=none, color=blue] 
table[x={n}, y={FBF_Expl}]{figures/ex3/dist_log.txt};
\addlegendentry{Pb.~\ref{prob:52}/Alg.~\ref{al:2012}-S}
\addplot [very thick, mark=none, color=blue] table[x={n}, 
y={FBF_Imp}]{figures/ex3/dist_log.txt};
\addlegendentry{Pb.~\ref{prob:53}/Alg.~\ref{al:2012}-P}
\addplot [dashed, very thick, mark=none, color=dgreen] 
table[x={n}, y={FB_Expl}] {figures/ex3/dist_log.txt};
\addlegendentry{Pb.~\ref{prob:52}/Alg.~\ref{al:fbp}-S}
\addplot [very thick, mark=none, color=dgreen] table[x={n}, 
y={FB_Imp}]{figures/ex3/dist_log.txt};
\addlegendentry{Pb.~\ref{prob:53}/Alg.~\ref{al:fbp}-P}
\addplot [dashed, very thick, mark=none, color=red] 
table[x={n}, y={KT_Expl}]{figures/ex3/dist_log.txt};
\addlegendentry{Pb.~\ref{prob:52}/Alg.~\ref{al:kt}-S}
\addplot [very thick, mark=none, color=red] table[x={n}, y={KT_Imp}] 
{figures/ex3/dist_log.txt};
\addlegendentry{Pb.~\ref{prob:53}/Alg.~\ref{al:kt}-P}
\end{axis}
\end{tikzpicture}
\centering
\caption{Left: Normalized distance in dB to the asymptotic image
produced by each algorithm versus execution time in seconds.
Right: Normalized distance in dB to the asymptotic image
produced by each algorithm versus iteration number.}
\label{fig:eX2costdist}
\end{figure}
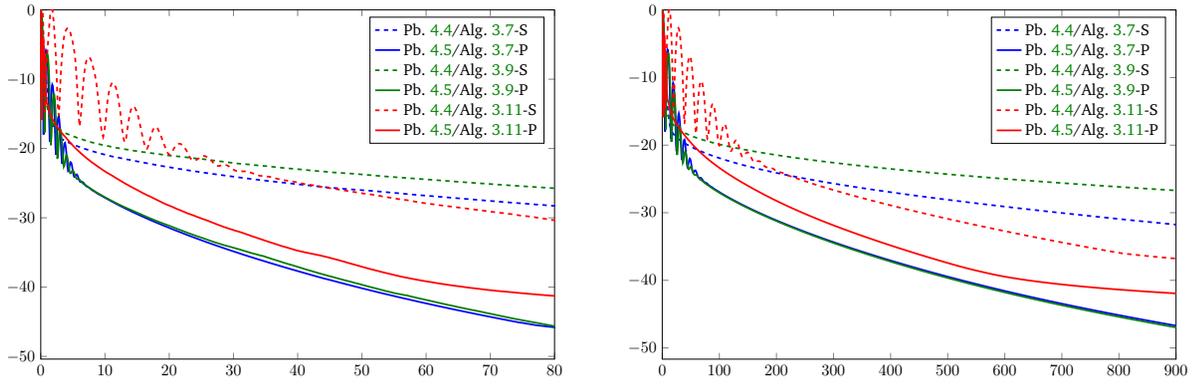

We apply to Problems~\ref{prob:52} and \ref{prob:53} 
Algorithms~\ref{al:2012}, 
\ref{al:fbp}, and \ref{al:kt} with all initial
vectors set to $0$. The following parameters are used, where
$\beta=\sqrt{\sum_{i\in I}\|L_i\|^2}+
\sum_{j\in J}\mu_j\|L_j\|^2$ (these parameters were found to
optimize the performance of each algorithm):
\begin{itemize}
\setlength{\itemsep}{-1pt}
\item
Algorithm~\ref{al:2012}-S (with smooth terms for
Problem~\ref{prob:52}): $\gamma_n\equiv 0.99/\beta$.
\item
Algorithm~\ref{al:2012}-P (fully proximal for 
Problem~\ref{prob:53}): $\gamma_n\equiv 0.99/\beta$.
\item
Algorithm~\ref{al:fbp}-S (with smooth terms for
Problem~\ref{prob:52}):
$\sigma_{1,n}\equiv 2/(5\beta)$ and 
$\tau_n\equiv 0.1/\beta$.
\item
Algorithm~\ref{al:fbp}-P (fully proximal for Problem~\ref{prob:53}):
$\sigma_{1,n}\equiv 1/\beta$, 
$\sigma_{2,n}\equiv 1/\beta$, 
$\sigma_{3,n}\equiv 1/\beta$,
and 
$\tau_n\equiv 1/\beta$.
\item
Algorithm~\ref{al:kt}-S (with smooth terms for
Problem~\ref{prob:52}):
$\gamma_n\equiv 1$, 
$\mu_{1,n}\equiv 0.1$,
$\mu_{2,n}\equiv 0.1$,
$\mu_{3,n}\equiv 0.01$, and 
$\lambda_n\equiv 1.9$. 
\item
Algorithm~\ref{al:kt}-P (fully proximal for Problem~\ref{prob:53}):
$\gamma_n\equiv 0.5$, 
$\mu_{1,n}\equiv 1.0$, 
$\mu_{2,n}\equiv 0.1$,
$\mu_{3,n}\equiv 0.01$, and
$\lambda_n\equiv 1.9$.
\end{itemize}
The proximity operators of $f$ and $g_3$ are
discussed in Section~\ref{sec:42}. The proximity operator of $g_1$
is computed similarly to that of $g_2$ in Problem~\ref{prob:51}
since, in view of Lemma~\ref{l:3}, we can apply \eqref{e:ex140}
with $\phi=|\cdot|$. It follows from \cite[Propositions~24.8(ii)
and 24.11]{Livre1} and \cite[Example~24.20]{Livre1}
that the proximity operator of $g_2$ for index $\gamma\in\RPP$ is
\begin{multline}
\prox_{\gamma g_2}\colon(x^{(i)})_{1\leq i\leq N}\mapsto
(p^{(i)})_{1\leq i\leq N},\quad\text{where}\quad
(\forall i\in\{1,\ldots,N\})\\
p^{(i)}=
\begin{cases}
x^{(i)},&\text{if}\;\:i\notin R;\\
y_1^{(i)}+\Bigg(1-\dfrac{\gamma}
{\max\big\{\big\|x^{(i)}-y_1^{(i)}\big\|_2,\gamma\big\}}\Bigg)
\big(x^{(i)}-y_1^{(i)}\big),&\text{if}\;\:i\in R.
\end{cases}
\end{multline}
As seen in Fig.~\ref{fig:ex3im}(d), the missing lines are
satisfactorily reconstructed. On the other hand, the convergence
profiles displayed in Fig.~\ref{fig:eX2costdist} indicate that 
the fully proximal algorithms behave better than their
gradient-based counterparts.

\subsection{Inconsistent convex feasibility problems}
\label{sec:44}

\subsubsection{Mathematical model}
As mentioned in the Introduction, the convex feasibility formalism 
first proposed in \cite{Youl82} has enjoyed continued interest 
from the image recovery community
\cite{Byrn14,Aiep96,Liuy14,Tofi16}.
A structured formulation of this problem is the following.

\begin{problem}
\label{prob:2}
\rm
Let $\HH$ be a real Hilbert space, let $E$ be a nonempty closed
convex subset of $\HH$, let $K$ be
a nonempty finite subset of $\NN$, and let 
$(\GG_k)_{k\in K}$ be a family of real Hilbert spaces. For every
$k\in K$, suppose that $0\neq L_k\in\BL(\HH,\GG_k)$
and let $C_k$ be a nonempty closed convex subset
of $\GG_k$. The goal is to 
\begin{equation}
\label{e:cfp1}
\text{find}\;\;x\in E\;\;\text{such that}\;\;
(\forall k\in K)\;\; L_kx\in C_k.
\end{equation}
\end{problem}

In applications, because of possible inaccuracies in 
\emph{a priori} knowledge, unmodeled dynamics, or too aggressive
confidence bounds on stochastic constraints, the above convex 
feasibility problem may turn out to be inconsistent 
\cite{Cens18,Sign94,Aiep96,Gold85,Youl86}, i.e.,
$E\cap\bigcap_{k\in K}L_k^{-1}(C_k)=\emp$. To deal with this
situation, 
we propose the following variational formulation as a relaxation of
\eqref{e:cfp1}. 

\begin{problem}
\label{prob:3}
\rm
Consider the setting of Problem~\ref{prob:2}.
Let $(I,J)$ be a partition of $K$ such that, for every $i\in I$, 
$\phi_i\in\Gamma_0(\RR)$ is an even function that vanishes only
at $0$ and, for every $j\in J$, $\psi_j\colon\RR\to\RR$ is an
even differentiable convex function that vanishes only at $0$ with
a Lipschitzian derivative. The problem is to 
\begin{equation}
\label{e:prob2}
\minimize{x\in E}{
\sum_{i\in I}\,\phi_i\big(d_{C_i}(L_ix)\big)+
\sum_{j\in J}\,\psi_j\big(d_{C_j}(L_jx)\big)}.
\end{equation}
\end{problem}

Problem~\ref{prob:3} unifies several formulations that have been
proposed in the literature as surrogates to the possibly
inconsistent Problem~\ref{prob:2}:
\begin{itemize}
\setlength{\itemsep}{-1pt}
\item
If \eqref{e:cfp1} happens to have solutions, they are the same 
as those of \eqref{e:prob2}.
\item
In \eqref{e:prob2}, $E$ plays the role of a hard
constraint; if no such constraint is present, one can set $E=\HH$.
Further hard constraints can be modeled by
taking $\phi_i=\iota_{\{0\}}$ for certain $i\in I$.
\item
Suppose that $J=\emp$ and $(\forall i\in I)$
$\phi_i=\iota_{\{0\}}$. Then \eqref{e:prob2} reverts to
\eqref{e:cfp1}. 
\item
Let $(\omega_j)_{j\in J}$ be real numbers in $\rzeroun$ such that
$\sum_{j\in J}\omega_j=1$. Suppose that $I=\emp$, $E=\HH$, and 
$(\forall j\in J)$ $\GG_j=\HH$, $L_j=\Id$, and
$\psi_j=\omega_j|\cdot|^2/2$. Then we recover
the least-squares formulation of \cite{Sign94}, namely
the problem of minimizing 
$\sum_{j\in J}\,\omega_jd_{C_j}^2$ over $\HH$.
\item
Let $(\omega_j)_{j\in J}$ be real numbers in $\rzeroun$ such that
$\sum_{j\in J}\omega_j=1$. Suppose that $I=\emp$ and 
$(\forall j\in J)$ $\GG_j=\HH$, $L_j=\Id$, and
$\psi_j=\omega_j|\cdot|^2/2$. Then we recover
the hard-constrained least-squares formulation of 
\cite{Sign99}, namely
\begin{equation}
\label{e:1999}
\minimize{x\in E}{\frac{1}{2}\sum_{j\in J}\,\omega_j
d_{C_j}^2(x)}.
\end{equation}
For $J=\{1\}$, this framework reduces to the formulation
proposed in \cite{Gold85,Youl86}.
\item
Let $(\omega_j)_{j\in J}$ be real numbers in $\rzeroun$ such that
$\sum_{j\in J}\omega_j=1$. Suppose that $I=\emp$ and 
$(J_1,J_2)$ is a partition of $J$. Suppose further that
$(\forall j\in J_1)$ $\GG_j=\HH$, $L_j=\Id$, and
$\psi_j=\omega_j|\cdot|^2/2$, and that
$(\forall j\in J_2)$ $\psi_j=\omega_j|\cdot|^2/2$.
Then we recover the formulation of \cite{Cens05}, namely
\begin{equation}
\label{e:2005}
\minimize{x\in E}{
\frac{1}{2}\sum_{j\in J_1}\,\omega_jd_{C_j}^2(x)+
\frac{1}{2}\sum_{j\in J_2}\,\omega_jd_{C_j}^2(L_jx)}.
\end{equation}
\item
Let $(\HH_\ell)_{1\leq \ell\leq m}$ and $(\KK_k)_{1\leq k\leq p}$ be 
real Hilbert spaces. For every $\ell\in\{1,\ldots,m\}$ and every
$k\in\{1,\ldots,p\}$, let $C_\ell$ be a nonempty closed convex
subset of $\HH_\ell$, let $E_k$ be a nonempty closed convex subset
of $\KK_k$, and let 
$M_{k\ell}\in\BL(\HH_\ell,\KK_k)$. Consider the multivariate
convex feasibility problem
\begin{equation}
\label{e:9-21}
\text{find}\;\;x_1\in C_1,\ldots, x_m\in C_m\;\;\text{such that}\;\;
\sum_{\ell=1}^mM_{1\ell}x_\ell\in
E_1,\ldots,\,\sum_{\ell=1}^mM_{p\ell}x_\ell\in E_p.
\end{equation}
Now let $(I,J_0)$ be a partition of $\{1,\ldots,m\}$, let
$(\phi_i)_{i\in I}$ be even functions in $\Gamma_0(\RR)$ vanishing
only at $0$, and, for every $j\in J_0$, let $\psi_j\colon\RR\to\RR$
be an even differentiable convex function that vanishes only at $0$
with a Lipschitzian derivative. A relaxation of \eqref{e:9-21}
proposed in \cite[Section~3.3]{Nmtm09} is
\begin{equation}
\label{e:2009}
\minimize{x_1\in\HH_1,\ldots,\,x_m\in\HH_m}{
\sum_{i\in I}\phi_i\big(d_{C_i}(x_i)\big)+
\sum_{j\in J_0}\psi_j\big(d_{C_j}(x_j)\big)+
\frac{1}{2p}\sum_{k=1}^pd_{E_k}^2
\Bigg(\sum_{\ell=1}^mM_{k\ell}x_\ell\Bigg)}.
\end{equation}
Now let 
\begin{equation}
\begin{cases}
\HH=\bigoplus_{\ell=1}^m\HH_\ell,\;
\GG_0=\bigoplus_{k=1}^p\KK_k\\
L_0\colon\HH\to\GG_0\colon(x_1,\ldots,x_m)\mapsto
\big(\sum_{\ell=1}^mM_{1\ell}x_\ell,\ldots,\sum_{\ell=1}^m
M_{p\ell}x_\ell\big)\\ 
C_0=E_1\times\cdots\times E_p,\;\psi_0=|\cdot|^2/(2p)\\
(\forall k\in I\cup J_0)\;\GG_k=\HH_k,\;
L_k\colon\HH\to\GG_k\colon(x_1,\ldots,x_m)\mapsto x_k\\
J=\{0\}\cup J_0,\; E=\HH.
\end{cases}
\end{equation}
Then \eqref{e:9-21} reduces to an instance of
\eqref{e:cfp1}, and \eqref{e:2009} of \eqref{e:prob2}.
\end{itemize}

\begin{remark}
\label{r:47}
\rm
Problem~\ref{prob:3} corresponds to the special
case of Problem~\ref{prob:1} in which
$f=\iota_E$, $(\forall i\in I)$ $g_i=\phi_i\circ d_{C_i}$,
and $(\forall j\in J)$ $h_j=\psi_j\circ d_{C_j}$.
To solve it with a fully proximal algorithm, one can use
Lemma~\ref{l:3} for the nonsmooth terms, and 
Example~\ref{ex:2} or Example~\ref{ex:7} for the smooth ones.
\end{remark}

\begin{remark}
\label{r:48}
\rm
The proposed variational formulation \eqref{e:prob2} is also of 
interest beyond the field of image recovery. For instance, the
set-theoretic Fermat-Weber (a.k.a.\ Heron) problem arising in
location theory is to \cite{Mord12}
\begin{equation}
\minimize{x\in\HH}{\sum_{i\in I}d_{C_i}(x)}.
\end{equation}
Problem~\ref{prob:3} therefore provides variants and
generalizations of this formulation. 
\end{remark}

\subsubsection{Application to image reconstruction from phase}
\begin{figure}[t]
\centering
\begin{tabular}{c@{}c@{}c@{}}
\includegraphics[width=4.7cm]{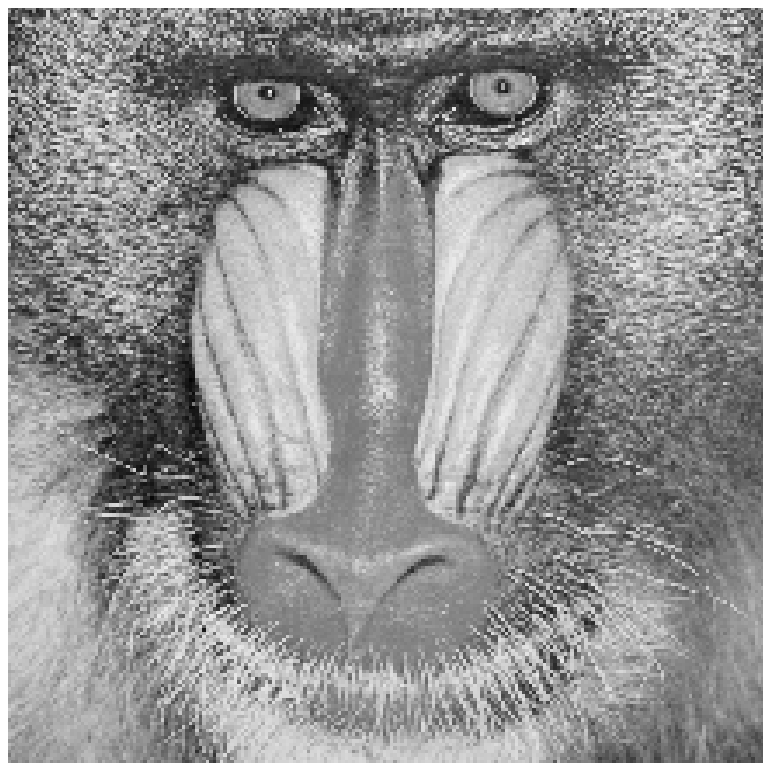}&
\hspace{0.4cm}
\includegraphics[width=4.7cm]{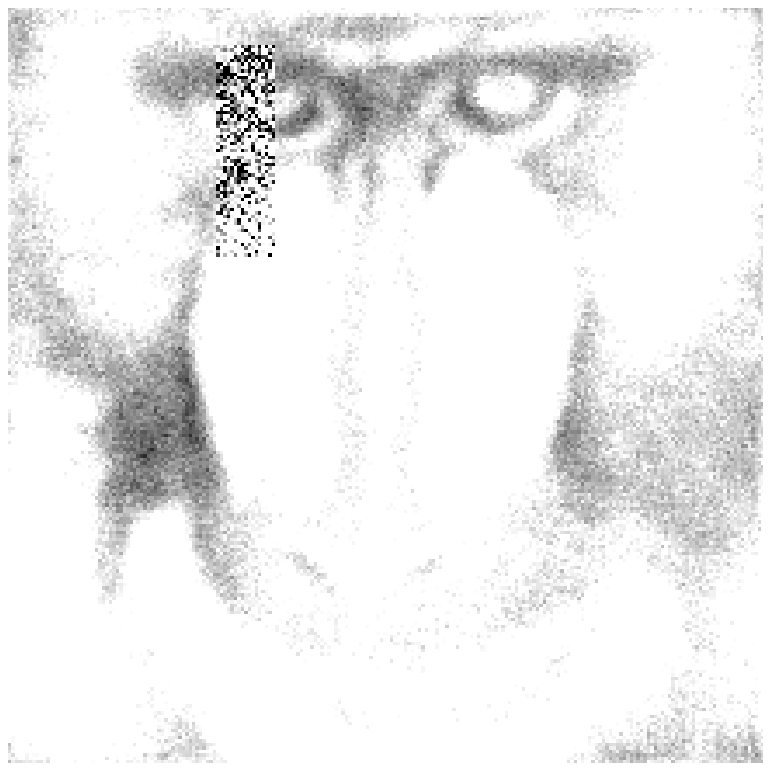}&
\hspace{0.4cm}
\includegraphics[width=4.7cm]{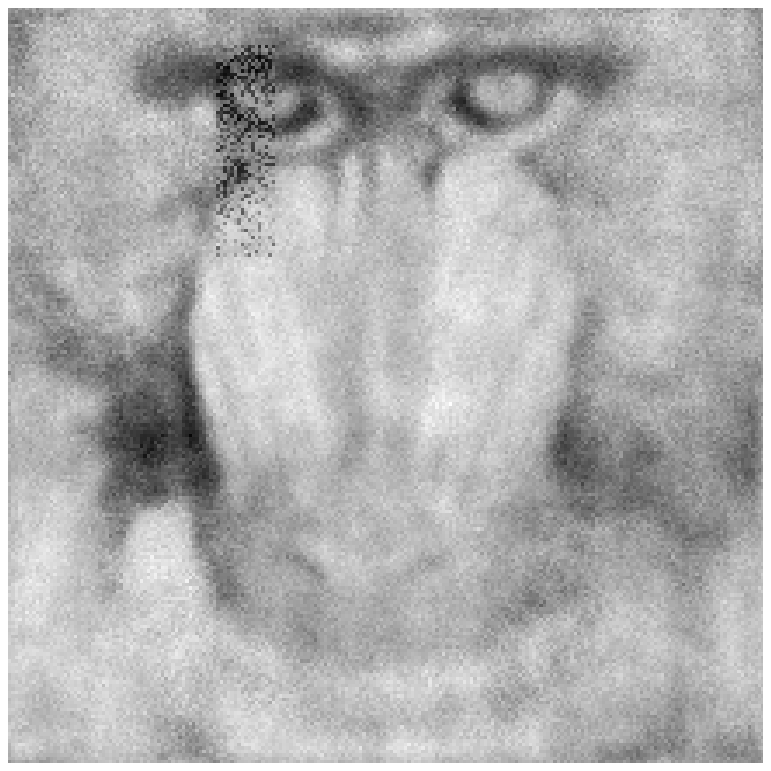}\\
\small{(a)} & \small{(b)} & \small{(c)}\\
\hspace{0.4cm}
\end{tabular} 
\caption{(a) Original image $\overline{x}$. 
(b) Reference image $r$.
(c) Reconstructed image (all algorithms yield visually equivalent
images).}
\label{fig:ex4im}
\end{figure}

This numerical example revolves around the classical problem of
recovering an
image $\overline{x}$ from the observation of its Fourier phase 
$\theta=\angle\,\widehat{\overline{x}}$ \cite{Kerm70,Levi83}. 
The original $512\times 512$ image $\overline{x}$ is shown in
Fig.~\ref{fig:ex4im}(a). The problem is modeled by
Problem~\ref{prob:2} as a convex feasibility problem with the
following constraint sets.
\begin{itemize}
\setlength{\itemsep}{-1pt}
\item
Mean pixel value: 
$C_1=\menge{x\in\RR^N}{\scal{x}{1}=\mu}$.
\item
Upper bound on the norm of the gradient: $Dx\in C_2$, with
\begin{equation}
C_2=\menge{y\in\RR^N\times\RR^N}{\|y\|_2\leq\eta}, 
\end{equation}
where $D$ is
as in \eqref{e:p2}.
\item
Phase: 
$C_3=\menge{x\in\RR^N}{\angle\widehat{x}=\theta}$.
\item
Proximity to the reference image $r$ of Fig.~\ref{fig:ex4im}(b):
$C_4=\menge{x\in\RR^N}{\|x-r\|_2\leq\xi}$.
The image $r$ is a blurred and noise corrupted version of
$\overline{x}$, which is further degraded by saturation (the pixel
values beyond 130 are clipped to 130) and the addition of a local
high intensity noise on a rectangular area around the right eye.
The only information available to the user is the bound $\xi$ 
on the distance of $r$ to the true image. 
\item
Pixel range: $C_5=[0,255]^N$.
\end{itemize}
Now set $E=C_5$, $K=\{1,2,3,4\}$, $\HH=\GG_1=\GG_3=\GG_4=\RR^N$,
$\GG_2=\RR^N\times\RR^N$, and $L_2=D$ in 
Problem~\ref{prob:2}. Then the 
feasibility problem \eqref{e:cfp1} amounts to finding an image 
$x\in C_1\cap D^{-1}(C_2)\cap C_3\cap C_4\cap C_5$. Because of 
inaccuracies in the values of $\mu$, $\xi$, $\theta$, and
$\eta$, this problem turns out to be inconsistent and we 
therefore turn to the formulation of
Problem~\ref{prob:3}.  
To ensure the robustness of the model to possible outliers, we 
adopt a constrained
Huber framework, namely
\begin{equation}
\minimize{x\in C_5}{\mathfrak{h}_{\rho_1}(d_{C_1}(x))+
\mathfrak{h}_{\rho_2}(d_{C_2}(Dx))+\mathfrak{h}_{\rho_3}
(d_{C_3}(x))+\mathfrak{h}_{\rho_4}(d_{C_4}(x))},
\end{equation}
where the functions $(\mathfrak{h}_{\rho_i})_{1\leq i\leq 4}$ are
defined in \eqref{e:huber} and $\rho_1=\rho_2=\rho_3=1000$ and 
$\rho_4=5000$. Two versions of Problem~\ref{prob:1} are employed.
\begin{problem}
\label{prob:54}
\rm
In Problem~\ref{prob:1}, set 
$f=\iota_{C_5}$, $I=\emp$, $J=\{1,2,3,4\}$, 
$h_1=\mathfrak{h}_{\rho_1}\circ d_{C_1}$, $L_1=\Id$,
$h_2=\mathfrak{h}_{\rho_2}\circ d_{C_2}$, $L_2=D$,
$h_3=\mathfrak{h}_{\rho_3}\circ d_{C_3}$, $L_3=\Id$,
$h_4=\mathfrak{h}_{\rho_4}\circ d_{C_4}$, and $L_4=\Id$.
\end{problem}

\begin{problem}[fully proximal]
\label{prob:55}
\rm
In Problem~\ref{prob:1}, set 
$f=\iota_{C_5}$, $I=\{1,2,3,4\}$, $J=\emp$,
$g_1=\mathfrak{h}_{\rho_1}\circ d_{C_1}$, $L_1=\Id$, 
$g_2=\mathfrak{h}_{\rho_2}\circ d_{C_2}$, $L_2=D$,
$g_3=\mathfrak{h}_{\rho_3}\circ d_{C_3}$, $L_3=\Id$, 
$g_4=\mathfrak{h}_{\rho_4}\circ d_{C_4}$, and $L_4=\Id$.
\end{problem}

We apply to Problems~\ref{prob:54} and \ref{prob:55} 
Algorithms~\ref{al:2012}, 
\ref{al:fbp}, and \ref{al:kt} with all initial
vectors set to $0$. The following parameters are used, where
$\beta=\sqrt{\sum_{i\in I}\|L_i\|^2}+
\sum_{j\in J}\mu_j\|L_j\|^2$ (these parameters were found to
optimize the performance of each algorithm):
\begin{itemize}
\setlength{\itemsep}{-1pt}
\item
Algorithm~\ref{al:2012}-S (with smooth terms for
Problem~\ref{prob:54}): $\gamma_n\equiv 0.99/\beta$.
\item
Algorithm~\ref{al:2012}-P (fully proximal for 
Problem~\ref{prob:55}): $\gamma_n\equiv 0.99/\beta$.
\item
Algorithm~\ref{al:fbp}-S (with smooth terms for
Problem~\ref{prob:54}): $\tau_n\equiv 1.99/\beta$.
\item
Algorithm~\ref{al:fbp}-P (fully proximal for Problem~\ref{prob:55}):
$\sigma_{1,n}\equiv 1/(1.1\beta)$, 
$\sigma_{2,n}\equiv 1/(1.1\beta)$, 
$\sigma_{3,n}\equiv 1/(1.1\beta)$,
$\sigma_{4,n}\equiv 1/(1.1\beta)$,
and 
$\tau_n\equiv 1/\beta$.
\item
Algorithm~\ref{al:kt}-S (with smooth terms for
Problem~\ref{prob:54}):
$\gamma_n\equiv 0.5$, 
$\mu_{1,n}\equiv 0.99$,
$\mu_{2,n}\equiv 0.99$,
$\mu_{3,n}\equiv 0.99$,
$\mu_{4,n}\equiv 0.99$, and 
$\lambda_n\equiv 1.9$. 
\item
Algorithm~\ref{al:kt}-P (fully proximal for Problem~\ref{prob:55}):
$\gamma_n\equiv 0.25$, 
$\mu_{1,n}\equiv 2.0$, 
$\mu_{2,n}\equiv 2.0$,
$\mu_{3,n}\equiv 0.5$,
$\mu_{4,n}\equiv 2.0$, and
$\lambda_n\equiv 1.9$.
\end{itemize}
\begin{figure}[ht!]
\begin{tikzpicture}[scale=0.55]
\hskip -4mm
\begin{axis}[height=10cm,width=14cm, legend style={at={(0.02,0.02)},
legend cell align={left}, anchor=south west}, 
xmin =0, xmax=120, ymax=0]
\addplot [dashed, very thick, mark=none, color=blue] 
table[x={time_FBF_Expl}, y={FBF_Expl}]
{figures/ex4/dist_time_log.txt};
\addlegendentry{Pb.~\ref{prob:54}/Alg.~\ref{al:2012}-S}
\addplot [very thick, mark=none, color=blue] table[x={time_FBF_Imp}, 
y={FBF_Imp}]{figures/ex4/dist_time_log.txt};
\addlegendentry{Pb.~\ref{prob:55}/Alg.~\ref{al:2012}-P}
\addplot [dashed, very thick, mark=none, color=dgreen] 
table[x={time_FB_Expl}, y={FB_Expl}] {figures/ex4/dist_time_log.txt};
\addlegendentry{Pb.~\ref{prob:54}/Alg.~\ref{al:fbp}-S}
\addplot [very thick, mark=none, color=dgreen] 
table[x={time_FB_Imp}, 
y={FB_Imp}]{figures/ex4/dist_time_log.txt};
\addlegendentry{Pb.~\ref{prob:55}/Alg.~\ref{al:fbp}-P}
\addplot [dashed, very thick, mark=none, color=red] 
table[x={time_KT_Expl}, y={KT_Expl}]{figures/ex4/dist_time_log.txt};
\addlegendentry{Pb.~\ref{prob:54}/Alg.~\ref{al:kt}-S}
\addplot [very thick, mark=none, color=red] table[x={time_KT_Imp}, 
y={KT_Imp}] 
{figures/ex4/dist_time_log.txt};
\addlegendentry{Pb.~\ref{prob:55}/Alg.~\ref{al:kt}-P}
\end{axis}
\end{tikzpicture}
\begin{tikzpicture}[scale=0.55]
\begin{axis}[height=10cm,width=14cm,legend style={at={(0.02,0.02)},
legend cell align={left}, anchor=south west}, 
xmin =0, xmax=370, ymax=0]
\addplot [dashed, very thick, mark=none, color=blue] 
table[x={n}, y={FBF_Expl}]{figures/ex4/dist_log.txt};
\addlegendentry{Pb.~\ref{prob:54}/Alg.~\ref{al:2012}-S}
\addplot [very thick, mark=none, color=blue] table[x={n}, 
y={FBF_Imp}]{figures/ex4/dist_log.txt};
\addlegendentry{Pb.~\ref{prob:55}/Alg.~\ref{al:2012}-P}
\addplot [dashed, very thick, mark=none, color=dgreen] 
table[x={n}, y={FB_Expl}] {figures/ex4/dist_log.txt};
\addlegendentry{Pb.~\ref{prob:54}/Alg.~\ref{al:fbp}-S}
\addplot [very thick, mark=none, color=dgreen] table[x={n}, 
y={FB_Imp}]{figures/ex4/dist_log.txt};
\addlegendentry{Pb.~\ref{prob:55}/Alg.~\ref{al:fbp}-P}
\addplot [dashed, very thick, mark=none, color=red] 
table[x={n}, y={KT_Expl}]{figures/ex4/dist_log.txt};
\addlegendentry{Pb.~\ref{prob:54}/Alg.~\ref{al:kt}-S}
\addplot [very thick, mark=none, color=red] table[x={n}, y={KT_Imp}] 
{figures/ex4/dist_log.txt};
\addlegendentry{Pb.~\ref{prob:55}/Alg.~\ref{al:kt}-P}
\end{axis}
\end{tikzpicture}
\centering
\caption{Left: 
Normalized distance in dB to the asymptotic image
produced by each algorithm versus execution time in seconds. 
Right: 
Normalized distance in dB to the asymptotic image
produced by each algorithm versus iteration number.}
\label{fig:eX4dist}
\end{figure}
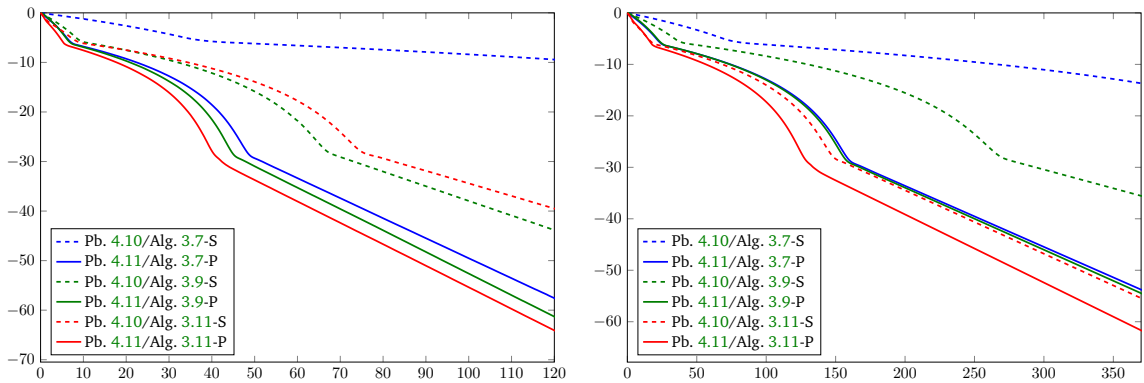
The gradient and proximity operators of the functions 
$(g_i)_{1\leq i\leq 4}$ are derived directly from
Example~\ref{ex:2}. They involve the projection operators 
$(\proj_{C_i})_{1\leq i\leq 4}$, which can be found
in \cite{Aiep96,Youl82}, as well as the proximity operator of 
$\mathfrak{h}_\rho$, which can be found in 
\cite[Example~24.9]{Livre1}.
As seen in Fig.~\ref{fig:ex4im}(c), despite the inconsistencies in
the a priori knowledge, the reconstructed image captures
important features of the original image. 
The results of Fig.~\ref{fig:eX4dist} show the faster
convergence of the fully proximal algorithms compared to
the gradient-based ones.

\end{document}